\documentclass[10pt]{amsart}

\usepackage{amsmath}
\usepackage{amsfonts}
\usepackage{amssymb}
\usepackage{amsthm}
\usepackage{hyperref}
\usepackage{color}
\usepackage{url}
\usepackage{graphicx}
\usepackage{tikz-cd}
\newtheorem{theorem}{Theorem}[section]
\newtheorem{defn}[theorem]{Definition}
\newtheorem{lem}[theorem]{Lemma}
\newtheorem{cor}[theorem]{Corollary}
\newtheorem{prop}[theorem]{Proposition}
\theoremstyle{definition}
\newtheorem{definition}{Definition}[section]
\newtheorem{rem}[definition]{Remark}
\newtheorem{qes}[definition]{Question}
\newtheorem{exm}[definition]{Example}

\newcommand{\git}{\mathbin{
		\mathchoice{/\mkern-6mu/}% \displaystyle
		{/\mkern-6mu/}% \textstyle
		{/\mkern-5mu/}% \scriptstyle
		{/\mkern-5mu/}}}% \scriptscriptstyle

\DeclareMathOperator{\coker}{coker}

\DeclareMathOperator{\sHom}{\textit{Hom}}
\DeclareMathOperator{\rHom}{\mathit{RHom}}
\DeclareMathOperator{\Hom}{Hom}
\DeclareMathOperator{\Spec}{\mathrm{Spec}}
\DeclareMathOperator{\Sym}{\mathrm{Sym}}
\DeclareMathOperator{\Hglob}{H}
\DeclareMathOperator{\nbun}{\mathcal{N}^{vir}}

\DeclareMathOperator{\Cvir}{C_{\epsilon}}
\DeclareMathOperator{\om}{\sigma}

\DeclareMathOperator{\Mor}{\text{Mor}}
\newcommand{\coh}{Wedge system }
\DeclareMathOperator{\cS}{\mathcal{S}}
\DeclareMathOperator{\cQ}{\mathcal{Q}}
\DeclareMathOperator{\wstack}{\mathbf{WS}}
\DeclareMathOperator{\cO}{\mathcal{O}}
\DeclareMathOperator{\bun}{\mathbf{Bun}}
\DeclareMathOperator{\quot}{\mathsf{Quot}}
\DeclareMathOperator{\IQ}{\mathsf{IQ}}
\DeclareMathOperator{\fix}{\mathrm{F}}
\DeclareMathOperator{\vir}{\mathrm{vir}}
\DeclareMathOperator{\rank}{\text{rank}}
\DeclareMathOperator{\vd}{\mathrm{vd}}
\DeclareMathOperator{\OG}{OG}
\DeclareMathOperator{\SG}{SG}
\DeclareMathOperator{\G}{G}
\DeclareMathOperator{\quadric}{\mathrm{Q}}
\author{Shubham Sinha}
\address{Shubham Sinha, Department of Mathematics, University of California San Diego, 9500 Gilman Drive, La Jolla, CA 92093, USA}
\email{shs074@ucsd.edu}
\begin{document}
	\title[The virtual intersection theory of isotropic Quot Schemes]{The virtual intersection theory of isotropic Quot Schemes}
	\begin{abstract}
		Isotropic Quot schemes parameterize rank $r$ isotropic subsheaves of a vector bundle equipped with symplectic or symmetric quadratic form. We define a virtual fundamental class for isotropic Quot schemes over smooth projective curves. Using torus localization, we prescribe a way to calculate top intersection numbers of tautological classes, and obtain explicit formulas when $r=2$. These include and generalize the Vafa-Intriligator formula. In this setting, we compare the Quot scheme invariants with the invariants obtained via the stable map compactification.	
	\end{abstract}
	\maketitle
	\section{Introduction}
	
	The isotropic Grassmannian $\SG(r,\mathbb{C}^N)$ (or $\OG(r,\mathbb{C}^N)$) is the variety parameterizing $r$ dimensional isotropic subspaces of a vector space $\mathbb{C}^N$ endowed with symplectic (or symmetric) non-degenerate bilinear form. The classical intersection theory of the Grassmannian $\G(r,\mathbb{C}^N)$ and isotropic Grassmannians has been an important subject connecting many areas of mathematics.
	
	The Quot scheme is a natural generalization of Grassmannian. Fix a smooth projective curve $C$ of genus $g$. The Quot scheme $\quot_d(V,r,C)$ (for short $\quot_d$) parameterizes degree $-d$, rank $r$ sub-sheaves of a fixed vector bundle $V$ over $C$. 
	
	Let $L$ be a line bundle over $C$ and let $\sigma$ be a symplectic or symmetric non-degenerate $L$-valued form on $V$: \[\sigma: V\otimes V\to L.  \]
	A subsheaf $S\subset V$ is isotropic if the restriction $\sigma|_{S\otimes S}=0$. The isotropic Quot scheme $\IQ_d(V,\sigma,r,C)$ (for short $\IQ_d$) is the closed subscheme of $\quot_d$  consisting of isotropic subsheaves.
	
	When $V$ is the trival rank $N$ bundle, $\quot_d$ provides a natural compactification of $\Mor_d(C,\G(r,\mathbb{C}^N))$, the scheme parameterizing degree $d$ maps from $C$ to the Grassmannian $\G(r,\mathbb{C}^N)$. Moreover, when $L$ is trivial and $\sigma$ is induced by a symplectic or symmetric form on $\mathbb{C}^N$ (we call such $\sigma $ standard), $\IQ_d$ gives a natural compactification for the space of maps $\Mor_d(C,\SG(r,\mathbb{C}^N))$ and $\Mor_d(C,\OG(r,\mathbb{C}^N))$ respectively.
	
	Another way to compactify the morphism space is via stable maps. This compactification is important for defining quantum cohomology (see \cite{RuTi_quantum_coh}). A geometric comparison between the Quot scheme and the stable map compactification was done in \cite{Popa_Roth}.
	
	A presentation for the quantum cohomology of $\G(r,\mathbb{C}^N)$ was derived in \cite{Q_C_Fano}, and a formula for Gromov-Ruan-Witten (GRW) invariants was proven. The presentations for the quantum cohomology rings of the isotropic Grassmannians were obtained in \cite{QC_LG}, \cite{Quantum_pieri_OG} and \cite{QC_isotropic_grassmannians}. 
	
	The  intersection theory of the Quot scheme was studied extensively in \cite{1994alg.geom..3007B}, \cite{Bert_Dask_Went}, \cite{Bertram_1997} and \cite{marian2007}. In particular, GRW invariants were recovered and new calculations were performed in \cite{marian2007}. The isotropic analogue of the Quot scheme first appeared as the Lagrangian Quot scheme over $\mathbb{P}^1$ (parameterizing maximal rank isotropic subsheaves) in \cite{QC_LG}. The Lagrangian Quot schemes have been recently studied in all genera in \cite{cheong2020irreducibility},\cite{cheong2019counting}.

	In this paper, we construct a virtual fundamental class for $\IQ_d$ for all $V$, all ranks $r$, all degrees and all genera. When $V$ is trivial and $\sigma $ is standard, we use virtual localization \cite{Graber1997LocalizationOV} to study the virtual intersection theory of $\IQ_d$. We prescribe a way to calculate top intersection numbers of tautological classes, and obtain explicit formulas when $r=2$. We further compute the Gromov-Ruan-Witten invariants obtained via the stable map compactification for the corresponding isotropic Grassmannians and compare the answers.
	
	We will now describe the results in detail. 
	
	\subsection{The Virtual Fundamental Class} Isotropic Quot schemes are, in most cases, not smooth. To define invariants, we first construct a virtual fundamental class on the isotropic Quot scheme.
	
	In \cite{marian2007}, Marian and Oprea constructed a virtual fundamental class for the Quot schemes $\quot_d$; see also \cite{Virt_Euler_CFK}. The virtual fundamental class on the isotropic Quot scheme is not a direct consequence of their construction.
	
	Let us assume $\sigma$ is symplectic. We may replace $\wedge^2$ with $\Sym^2$ when $\sigma$ is symmetric to obtain the following results.
	
	The best scenario occurs when $V$ is the trivial vector bundle over $\mathbb{P}^1$. In this case, $\quot_d$ is a smooth scheme and $\IQ_d$ is the zero locus of a section of the vector bundle $\pi_*(\wedge^2\cS^\vee)$. Here, we consider the universal exact sequence over $C\times \IQ_d  $, \[0\to \cS\to p^*V\to \cQ\to 0,  \] where $p$ and $\pi$ are the projection maps to $C$ and $\IQ_d$ respectively.
	
	Unfortunately, for an arbitrary vector bundle $V$ over a higher genus curve $C$, $\quot_d$ is not smooth and $\pi_*(\wedge^2\cS^\vee)$ is not a vector bundle.
	
	Our first main result is
	\begin{theorem}\label{thm:POT}
		There is a morphism in the derived category
		\begin{equation}
			\mathbf{R}\pi_*(J^\bullet)^\vee \to \tau_{[-1,0]}\mathbb{L}_{\IQ_d}
		\end{equation}
		where $J^\bullet=[\rHom(\cS,\cQ)\to \sHom(\wedge^2\cS,p^*L)]$, which induces a $2$-term perfect obstruction theory and hence a virtual fundamental, $[\IQ_d]^{\vir }$, on the isotropic Quot scheme.
	\end{theorem}
	
	We prove Theorem \ref{thm:POT} in Section \ref{sec:Vir_Fund_class}.
	
	Over a closed point $[0\to S\to V\to Q\to 0]$ in $\IQ_d$, the tangent space and the obstruction space are given by the hypercohomology of the complex of sheaves $[\sHom(S,Q)\to \sHom(\wedge^2S,L)]$. The virtual dimension is 
	\begin{align*}
		\vd&=\begin{cases}
			\chi(S^\vee \otimes Q)-\chi(\wedge^2S^\vee\otimes L) & \text{when $\sigma $ is symplectic} \\
			\chi(S^\vee \otimes Q)-\chi(\Sym^2S^\vee\otimes L) & \text{when $\sigma $ is symmetric}
		\end{cases},
	\end{align*}
	where $\chi(E)$ denotes the Euler characteristic of a sheaf $E$. These are easy to calculate as an application of the Riemann-Roch formula.
	\begin{rem}
		When $2r=N$ and $\sigma$ is symplectic, the isotropic Quot scheme is irreducible and generically smooth \cite{cheong2020irreducibility} for $d>>0$ and its  dimension equals the virtual dimension obtained above. In this case, the virtual fundamental class agrees with the fundamental class.
	\end{rem}
	\begin{rem}
		Our method can also be extended to obtain a virtual fundamental class for the closed subscheme of $\quot_d$ parameterizing  subsheaves $S\to V$ isotropic with respect to higher order forms $\sigma:\wedge^{k}V\to L$ and $\sigma:\Sym^{k}V\to L$.
	\end{rem}
	For the rest of the introduction, we will assume that $V$ is a trivial vector bundle of even rank $N$. We will also assume that the line bundle $L$ is trivial and the non-degenerate symplectic or symmetric form $\sigma$ is standard.
	
	\subsection{Compatibility of virtual fundamental classes} 
	The group $G= Sp(N)$ (or $G=SO(N)$) acts on the isotropic Quot scheme with $\sigma $ symplectic (resp. symmetric). The perfect obstruction theory we construct is equivariant under any one-parameter subgroup $\mathbb{C}^*\subset G$. In this case, we use the virtual localization theorem \cite{Graber1997LocalizationOV} to study the virtual intersection theory of $\IQ_d$. This has been done extensively for  $\quot_d$ in \cite{marian2007}. 
	
	We first show a compatibility result for the virtual fundamental classes. Fix a point $q\in C$. There is a natural embedding 
	\begin{equation*}
		i_q:\IQ_d\to \IQ_{d+r}
	\end{equation*}
	which sends a subsheaf $S\subset \mathbb{C}^N\otimes\cO$ to the composition \[S(-q)\to S \to \mathbb{C}^N\otimes\cO, \] which is also an isotropic subsheaf of degree $-(d+r)$. 
	
	\begin{theorem}\label{compatibiliy_theorem}
		We have the following identity in the homology $H_*(\IQ_{d+r})$ :
		\begin{equation}
			{i_q}_*(c_{\text{top}}(\wedge^2\cS^\vee_q)^2 \cap [\IQ_d]^{\vir})=c_{\text{top}}({\cS_q^{\vee}})^N\cap  [\IQ_{d+r}]^{\vir}
		\end{equation} 
		where we assume that $\sigma$ is symplectic. The corresponding identity for symmetric form is obtained by replacing $\wedge ^2$ with $\Sym^2$.
	\end{theorem}
	
	This means that the virtual fundamental classes we construct, $[\IQ_d]^{\vir}$, are related as we vary the degree $d$ by a multiple of $r$. An analogous result was proven in the case of the Quot scheme in \cite{marian2007}.
	\subsection{Virtual Invariants}
	
	Let $\{1,\delta_1,\dots \delta_{2g}, \omega\}$ be a symplectic basis for the cohomology of $C$. Let the K\"unneth decomposition of $\cS^\vee$ over $C\times \IQ_d$ be 
	\begin{align*}
		c_i(\cS^\vee)= a_i\otimes 1 + \sum_{k=1}^{2g}b_i^k\otimes\delta_k+ f_i\otimes \omega,
	\end{align*}
	where $a_i\in H^{2i}(\IQ_d)$, $b_i^{k}\in H^{2i-1}(\IQ_d)$ and $f_i\in H^{2i-2}(\IQ_d)$. 
	
	The classes $a_i$ and $f_i$ have natural algebro-geometric descriptions. 
	For any point $q\in\IQ_d$, let $\cS_q$ be the restriction of $\cS$ to $\IQ_d\times \{q\}$. Then 
	\begin{align*}
		a_i=c_i(\cS^{\vee}_q),\hspace*{1.4cm} f_i=\pi_*c_i(\cS^{\vee}).
	\end{align*}
	
	The top intersections of the corresponding $a$-classes over $\quot_d$ match the GRW invariants for Grassmannians. The explicit answers were first obtained in the physics literature by Vafa and Intriligator \cite{INTRILIGATOR_1991}. In the mathematics literature, these formulas appeared in \cite{1994alg.geom..3007B}, \cite{Q_C_Fano} and \cite{marian2007}. 
	
	We are interested in understanding the intersection products of the above two kinds of classes evaluated on the virtual fundamental cycle. The virtual localization theorem \cite{Graber1997LocalizationOV} allows us to evaluate all monomials in $a_i$ and $f_i$ on the virtual fundamental class $[\IQ_d]^{\vir}$. However, closed form expressions are harder to write down due to the fact that the combinatorics becomes very involved. 
	
	When $r=2$, we prove a Vafa-Intriligator type formula for such intersection numbers. We achieve this by developing combinatorial techniques in Section \ref{sec:prelim_hilb} to evaluate and sum the fixed loci contributions. In the process, we simplify some of the combinatorics in \cite{marian2007}.
	
	At this point, we will have to distinguish the two cases depending on $\sigma$ being a symplectic or symmetric form.
	
	\subsection{When $\sigma$ is symplectic}
	When $\sigma$ is symplectic and $r=2$, the virtual dimension is \[\vd=(N-1)d-(2N-5)\bar{g},\] where we use the convention \[\bar{g}=g-1.\] We further define
	\[T_{d,g}(N)= \sum_{i=0}^{d}\binom{g}{i}(-N)^{-i}. \]
	The above expression equals $(1-1/N)^g$ when $d\ge g$. Note that the non-negativity of the virtual dimension implies that $d< g$ if and only if $\vd=0$ and $N=4$ or $\vd=0$ and $g=1$.
	
	\begin{theorem}\label{thm:r=2,sympl}
		Let $\sigma$ be a symplectic form and $m_1+2m_2=\vd\ge 0$. Then
		\begin{equation}\label{eq:a_1a_2formula}
			\int_{[\IQ_d]^{\vir}}^{}a_1^{m_1}a_2^{m_2} =u\frac{N}{2}T_{d,g}(N) \sum_{\zeta\ne\pm1}^{}(1+\zeta)^{m_1+d}\zeta^{m_2}J(1,\zeta)^{\bar{g}},
		\end{equation}
		where the sum is taken over $N^{th}$ roots of unity $\zeta\ne \pm 1$. Here $u=(-1)^{\bar{g}+d}$ and
		\begin{align*}
			J(z_1,z_2)&=N^2z_1^{-1}z_2^{-1}(z_1-z_2)^{-2}(z_1+z_2)^{-1}.
		\end{align*}
	\end{theorem}
	\begin{exm}\label{exm:N=4}
		When $N=4$, the virtual dimension $\vd=3d-3\bar{g}$. The above theorem specializes to 
		\begin{align*}
			\int_{[\IQ_d]^{\vir}}^{}a_1^{m_1}a_2^{m_2} =\begin{cases}
				2^{2d-m_2-\bar{g}}3^g & \vd >0\\
				2^{\bar{g}}(3^g+(-1)^{\bar{g}})& \vd=0.
			\end{cases}
		\end{align*}
		When $\vd=0$, the resulting invariant can be interpreted as a `virtual' count of isotropic subsheaves of $V$. This virtual count matches the enumerative count \cite{cheong2019counting} of the rank two maximal degree isotropic subbundle of a general rank 4 stable bundle endowed with an $\cO$-valued symplectic form.
	\end{exm}

\begin{exm}\label{exm:g=1}
When $g=1$, the virtual dimension $\vd=(N-1)d$. Then
\begin{align*}
	\int_{[\IQ_d]^{\vir}}^{}a_1^{\vd}=\begin{cases}
(-1)^d\frac{N-1}{2}[q^{Nd}]\Big(  \frac{N(1-q)^{N-1}}{(1-q)^N-q^N}-\frac{1}{1+2q}\Big)& d>0\\
\frac{N(N-2)}{2}& d=0
\end{cases}.
\end{align*}
\end{exm}
	
	We have the following results involving $f$-classes; the latter are typically intractable by other methods.
	\begin{theorem}\label{thm:f_classes_d>g}
		Let $m_1+m_2+1=\vd$ and $d>g$, then 
		\begin{align*}
			\int_{[\IQ_d]^{\vir}}f_2a_1^{m_1}a_2^{m_2}=&\bigg(1-\frac{1}{N}\bigg)^{g}\sum_{\zeta\ne \pm1}^{}\bigg( D\circ B(1,\zeta)-\frac{\zeta B(1,\zeta)}{(1+\zeta)}\bigg).
		\end{align*}
		where \[D \circ R(z_1,z_2)=\frac{z_1z_2}{2}\bigg(\frac{\partial}{\partial z_1}+\frac{\partial}{\partial z_2}\bigg)R(z_1,z_2)\] is a differential operator and
		\begin{equation*}
			B(z_1,z_2)=  u(z_1+z_2)^{m_1}(z_1z_2)^{m_2}\frac{(z_1+z_2)^{d-\bar{g}}}{(z_1-z_2)^{2\bar{g}}}\prod_{i=1}^{2} (Nz_i^{N-1})^{\bar{g}}.
		\end{equation*}
	\end{theorem}
	In Section \ref{sec:f_intersection}, we provide a complete answer for the intersection numbers of the form $ f_2^{\ell}a_1^{m_1}a_2^{m_2}\cap [\IQ_d]^{\vir}$ at the cost of making the formula more cumbersome. The answer involves higher degree differential operators. We remark here that our method can also be applied to obtain virtual intersection numbers involving higher powers of $f_2$ over the Quot schemes as well (for which closed form expressions were not known).
	
	\subsection{When $\sigma$ is symmetric}
	When $r=1$, every rank $r$ subsheaf of a symplectic vector bundle is isotropic. In this case $\IQ_d=\quot_d$. However, when $\sigma$ is a symmetric form, this is not the case. 
	\begin{prop}
		Let $r=1$, let $N$ be even and let $\sigma$ be a symmetric form. Then
		\begin{align*}
			\int_{[\IQ_d]^{\vir}}^{}a_1^{\vd}= (N-2)^g2^{2d-\bar{g}},
		\end{align*}
		where $\vd=(N-2)(d-\bar{g})$ is the virtual dimension and $d\ge g$.
	\end{prop}
	When $r=2$, the virtual dimension of $\IQ_d$ is \[\vd = (N-3)d-\bar{g}(2N-7).\]
	
	\begin{theorem}\label{thm:r=2_symmetric}
		Let $m_1+2m_2=\vd$ and $N=2n+2$.
		\begin{itemize}
			\item[(i)] When $m_2>0$, then
			\begin{equation*}
				\int_{[\IQ_d]^{\vir}}^{}a_1^{m_1}a_2^{m_2} =
				c\sum_{\zeta\ne \pm 1}(1+\zeta)^{m_1+d}\zeta^{m_2}J(1,\zeta)^{\bar{g}}
			\end{equation*}
			\item[(ii)] When $m_2=0$,
			\begin{equation*}
				\int_{[\IQ_d]^{\vir}}^{}a_1^{m_1} =c\bigg(4(-n)^{\bar{g}}+\sum_{\zeta\ne \pm 1}(1+\zeta)^{m_1+d}J(1,\zeta)^{\bar{g}}\bigg),
			\end{equation*}
		\end{itemize}
		where the sum is taken over $2n^{th}$ roots of unity $\zeta\ne \pm 1$. Here $u=(-1)^{\bar{g}+d}$,
		\begin{align*}
			c=u4^dnT_{d,g}(2n), \quad 
			J(z_1,z_2)=n^2(z_1+z_2)^{-1}(z_1-z_2)^{-2}.
		\end{align*}		
	\end{theorem}
	In the above theorem, there are two differences from Theorem \ref{thm:r=2,sympl} which make the proof more difficult. First, the case $m_2=0$ requires extra care. Second, in the sum above $\zeta$ is $(N-2)^{th}$ root of unity. This arises from picking a non-standard $\mathbb{C}^*$ action in the localization formula. In particular, the fixed loci thus obtained come equipped with a non-standard virtual structure. 
	
	We observe a surprising duality in the $a$-class intersection numbers over the symmetric and symplectic isotropic Quot schemes. We will later observe the same phenomenon for GRW invariants.
	\begin{cor}
		Let $\IQ_d$ (and $\widetilde{\IQ}_d$) be symplectic (respectively symmetric) isotropic Quot scheme parameterizing  rank $2$ degree $d$ isotropic subsheaves of $\mathbb{C}^N\otimes\cO$ (and $\mathbb{C}^{N+2}\otimes\cO$ respectively). Then, for integers $m_1,m_2$ such that $m_1+2m_2=(N-1)d-\bar{g}(2N-5)$ and $m_2-\bar{g}>0$, we have
		\begin{align*}
			\int_{[\widetilde{\IQ}_d]^{\vir}}^{}a_1^{m_1}a_2^{m_2-\bar{g}}=4^{d-2\bar{g}}\int_{[\IQ_d]^{\vir}}a_1^{m_1}a_2^{m_2}.
		\end{align*}
	\end{cor}

	\subsection{Gromov-Ruan-Witten Invariants}
	
	In the previous sections, we considered the Quot scheme compactification of the morphism space $\Mor_d(C,\SG(2,N))$ and $\Mor_d(C,\OG(2,N))$. 
	
	Let $(M,\omega)$ be a compact symplectic manifold with a generic almost complex structure $J$ tamed by $\omega$ (i.e. $\omega(v,Jv)>0$ for all non-zero $v\in TM $). We will further assume that $H_2(M,\mathbb{Z})\cong \mathbb{Z}$ and $M$ is positive in the sense that $c_1(TM,J) \cdot f_*[\mathbb{P}^1]>0$ for all non-constant $J$-holomorphic maps $f:\mathbb{P}^1\to M$. 

	The morphism space of $J$-holomorphic maps from $C$ to $(M,\omega)$ can be compactified by letting the curve $C$ `bubble' \cite{RuTi_quantum_coh}. The boundary of this compactification includes $C$ with finitely many trees of rational curves. This leads to the definition of quantum cohomology and Gromov-Ruan-Witten (GRW) invariants. We briefly describe these terms, but readers are suggested to see \cite{Q_C_Fano}, \cite{J-holo_curves_SalamonDuff} for more details.
	
	Let $\alpha\in H^2(M,\mathbb{Z})$ be a positive generator. Define the index $e$ of $M$ by $c_1(M)=e\alpha$. Let $d\in H^2(M,\mathbb{Z})$ and $\alpha_1,\dots,\alpha_s$ be cohomology classes in $H^*(M,\mathbb{Z})$ satisfying 
	\begin{align}\label{eq:exp_dim}
		\frac{1}{2}\sum_{i=1}^{s}\deg \alpha_i= ed+\dim(M)(1-g).
	\end{align}
	The right side of the above expression is the expected dimension of the moduli space of maps $f:C\to M$ with $f_*(C)=d\in H_2(M,\mathbb{Z})$.
	
	Let $B_1,\dots,B_s$ be a generic choice of the  Poincar\'e dual homology classes of $\alpha_1,\dots,\alpha_s$. Then for $s$ generic points $p_1,\dots,p_s\in C$, the GRW invariants \[\Phi_{g,d}(\alpha_1,\dots ,\alpha_s)\]
	is the algebraic count (considering sign and multiplicities) of $J$-holomorphic curves $f:C\to X$ such that $f(p_i)\in B_i$ and $f_*([C])=d$. The GRW invariants depend on the genus but not the complex structure of the curve.
	
	Quantum cohomology packages the information of $3$-point genus zero GRW invariants giving a deformation of the usual cohomology ring (see \cite{J-holo_curves_SalamonDuff} for more details). A presentation of quantum cohomology of $\SG(r,N)$ and $\OG(r,N)$ was described in \cite{QC_isotropic_grassmannians} and \cite{Quantum_pieri_OG}. In \cite{Q_C_IG}, the authors gave a simpler presentation for $\SG(2,N)$. We extend their result obtaining a similar presentation for $\OG(2,N)$.
	
	Let $N=2n+2$. We have the universal exact sequence \[0\to \cS\to \mathbb{C}^N\otimes\cO\to \cQ\to 0\] over $\OG(2,N)$. Let $\cS^{\perp}\subset \mathbb{C}^N\otimes\cO$ be the rank $N-2$ orthogonal complement.
	
	We have the following cohomology classes :
	\begin{itemize}
		\item The Chern classes $a_i= c_i(\cS^\vee)$ for $i\in\{1,2\}$.
		\item Let $b_i=c_{2i}(\cS^{\perp}/\cS)$ for $i\in \{1,\dots ,n-1 \}$. The bundle $\cS^{\perp}/\cS$ is self dual, hence all the odd Chern classes vanish.
		\item Let $\xi$ be the Edidin-Graham square root class \cite{Char_classes_quadric} of the bundle $\cS^{\perp}/\cS$. In particular, it satisfies \[ (-1)^{n-1}\xi^2=b_{n-1}. \] 
	\end{itemize}
	\begin{prop}
		The quantum cohomology ring $QH^*(\OG(2,2n+2),\mathbb{C})$ is isomorphic to the quotient of the ring $\mathbb{C}[a_1,a_2,b_1,\dots,b_{n-2}, \xi,q]$ by the ideal generated by the relations \[\xi a_2=0\] and 
		\begin{equation*}
			(1+(2a_2-a_1^2)x^2+a_2^2x^4)(1+b_1x^2+\cdots + b_{n-2}x^{2n-4}+(-1)^{n-1}\xi^2x^{2n-2})=1+4 qa_1x^{2n},
		\end{equation*}
		where $x$ is a formal variable.
	\end{prop}
	
Define the GRW invariant \[\langle a_1^{m_1}a_2^{m_2}\rangle_g=\Phi_{g,d}(a_1,\dots,a_1,a_2,\dots ,a_2),\] where $a_1$ and $a_2$ appear $m_1$ and $m_2$ times respectively; and $d$ is chosen (if possible) such that it satisfies \eqref{eq:exp_dim}.
	
	In \cite{Q_C_Fano}, Siebert and Tian gave a remarkable technique to compute the higher genus GRW invariants using a given presentation for the quantum cohomology. We explicitly calculate the GRW invariants for $\SG(2,N)$ and $\OG(2,N)$  in Theorems \ref{thm:GRW_invariant_SG} and \ref{thm:GRW_invariants_OG} respectively. We verify the slogan below for $r=2$.
	\[\text{``GRW Invariants} = \text{Virtual $a$-class intersections''}\]
	In particular, we prove the following theorem.
	\begin{theorem}
		Let $d$, $m_1$ and $m_2$ be non-negative integers such that $\vd=m_1+2m_2$ is the expected dimension. The GRW invariants for $\SG(2,N)$ (and $\OG(2,N)$)
		\begin{equation*}
			\langle a_1^{m_1}a_2^{m_2}\rangle_g=	\int_{[\IQ_d]^{\vir}}^{}a_1^{m_1}a_2^{m_2},
		\end{equation*}
where $\IQ_{d}$ is the symplectic (respectively symmetric) isotropic Quot scheme.
	\end{theorem}
	\begin{qes}
		In the large degree regime, we expect that $\IQ_{d}$ and the corresponding stable map compactification are irreducible and the above invariants are enumerative. The irreducibility of the Lagrangian Quot schemes for $d>>0$ is proven in \cite{cheong2020irreducibility}. 
	\end{qes}
	\subsection{Virtual Euler Characteristic}
	The topological Euler characteristics of schemes $\IQ_{d}$ is given by
	\begin{equation*}
		\sum_{d=0}^{\infty}e(\IQ_d)q^d=2^r\binom{n}{r}(1-q)^{r(2g-2)},
	\end{equation*}
where $N=2n$.
	
	Let $X$ be a scheme admitting a 2-term perfect obstruction theory. The virtual Euler characteristic is defined \cite{Virt_Euler_FG}, \cite{Virt_Euler_CFK}
	\begin{equation*}
		e^{\vir}(X)=\int_{[X]^{\vir}}^{}c(T_X^{\vir}).
	\end{equation*}
	The virtual Euler characteristic of Quot scheme parameterizing zero dimensional quotients over surfaces were calculated in \cite{oprea2021quot}.
	
	When $X$ is smooth and the obstruction bundle vanishes, the virtual Euler characteristic $e^{\vir}(X)$ matches the topological Euler characteristic of $X$. The isotropic Quot schemes, $\IQ_1$, are smooth for $C=\mathbb{P}^1$ and all values of $N=2n$ and $r$. 
	By contrast, the isotropic Quot schemes $\IQ_{d}$ are not smooth for $d>1$ even when $C=\mathbb{P}^1$. Thus the virtual Euler characteristics, $e^{\vir}(\IQ_{d})$, are new invariants. While we do not a have a closed form expression for these power series, nonetheless we find a finite number of values using Sagemath \cite{sagemath}. We provide a small list of these invariants in Section \ref{sec:vir_Euler_char}. 
	
\begin{center}
		\begin{figure}\label{fig:plot}
			\includegraphics*[width=10cm]{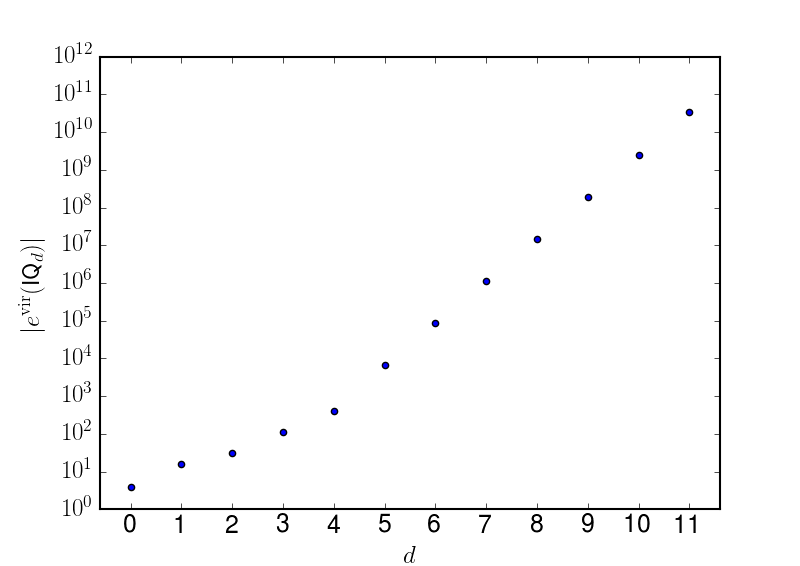}
			\caption{The absolute value of the virtual Euler characteristic of $\IQ_{d}$ in log scale, where $r=2$ and $\sigma$ is the standard symplectic form on $\mathbb{C}^4\otimes \cO$ over $\mathbb{P}^1$.}
		\end{figure}
	\end{center}

	When $r=2$, $N=4$ and $\sigma$ is symplectic, we plot a log scale graph for the absolute value of $e^{\vir}(\IQ_{d})$. The plot (see Figure \ref{fig:plot}) indicates an exponential growth in contrast with the polynomial expression for the topological Euler characteristics. 

	\begin{qes}
		Find a closed form expression for the virtual Euler characteristic of $\quot_d$ and $\IQ_{d}$ for all genus $g$ and all ranks $r$ and $N$.
	\end{qes}

	\subsection{Plan of the paper} We construct the virtual fundamental class over $\IQ_d$ in Section \ref{sec:Vir_Fund_class}, thus proving Theorem \ref{thm:POT}. In Sections \ref{sec:Nvir_sympl} and \ref{sec:Nvir_symm}, we will describe the torus action on $\IQ_d$ and find an expression for the equivariant virtual normal bundles over the fixed loci. Section \ref{sec:prelim_hilb} is technical, and it contains calculations on the product of symmetric powers of curves. These will be used in Sections  \ref{sec:Int_a_classes} and \ref{sec:f_intersection} to prove Theorems \ref{thm:r=2,sympl} and \ref{thm:r=2_symmetric}. The quantum cohomology and GRW invariants are calculated in Section \ref{sec:GRW_Invariants}; this section is technically disjoint from all the other sections.
	
	\subsection{Acknowledgements}
	I would like to thank Professor Dragos Oprea for suggesting this problem and for numerous useful conversations. The computational exploration and verification was done using the open source mathematical software Sage \cite{sagemath}. This work was partially funded by NSF grant DMS 1802228.
	
	\section{Virtual Fundamental Class}\label{sec:Vir_Fund_class}
	
	We will construct a natural 2-term perfect obstruction theory for the isotropic Quot Schemes $\IQ_{d}$ over smooth projective curves. This yields a virtual fundamental class using the results in \cite{Behrend_1997} and \cite{li1996virtual}. The argument can be slightly simplified for trivial vector bundles $V=\mathbb{C}^N\otimes\cO$ over $\mathbb{P}^1$ and we will explain this case first.
	
	We will assume that $\sigma$ is a symplectic non-degenerate bilinear form on a vector bundle $V$, although similar results hold for symmetric bilinear forms and can be proved verbatim replacing $\wedge^2$ with $\Sym^2$.
	
	\subsection{Background}
	We will briefly describe the results pertaining to the construction of virtual fundamental classes in \cite{Behrend_1997}.
	
	Let $X$ be a scheme (or a stack) over a scheme (or a stack) $S$ and $\mathbb{L}_{X/S}$ be the relative cotangent complex.
	\begin{defn}
		A 2-term relative perfect obstruction theory is a morphism in the derived category\[ \phi : E^{\bullet}\to \tau_{[-1,0]}\mathbb{L}_{X/S}, \] where $E^\bullet= [E^{-1}\to E^{0}]$ is a complex of vector bundles over $X$ of amplitude contained in $[-1,0]$ and satisfies:
		\begin{itemize}
			\item $h^0$ is an isomorphism and
			\item $h^{-1}$ is a surjection.
		\end{itemize}
	\end{defn} 
	
	Let $[E_0\to E_1]$ be the dual of $E^\bullet$. Given a 2-term perfect obstruction theory, \cite{Behrend_1997} and \cite{li1996virtual} define a cone inside $E_{1}$. The virtual fundamental class is then defined to be an element in $H_{2e}(X)$ given by the refined intersection of the cone with the zero section of $E_1$. Here $e= \rank E_0-\rank E_1$ is called the virtual dimension of $X$.
	
	For practical purposes, we only need the description of the virtual tangent (or cotangent) bundle, which is an element in the $K$-theory \[T^{\vir}_X=[E_0]-[E_1]\in K^0(X).\]
	
	The simplest case is when $X$ is a closed subscheme of a smooth scheme $Y$ cut out by a section $s$ of a vector bundle $V$ over $Y$. In this case, there is a natural 2-term perfect obstruction theory given by $[V^{\vee}|_X\to \Omega_Y|_X]$. Note that when $s$ is a regular section, we get the usual fundamental class.
	
	For the remainder of this section, we provide a 2-term perfect obstruction theory for $\IQ_d$.
	\subsection{Genus 0}
	Over $\mathbb{P}^1$, the Quot scheme $\quot_d(\mathbb{C}^{N},r,\mathbb{P}^1)$ is smooth for any choice of $N,r$ and $d$. The isotropic Quot scheme $\IQ_{d}$ is smooth for $d=0,1$ for all $r$ and $N$, but it is singular for higher values of $d$. 
	
	The isotropic Quot schemes $\IQ_d$ can be described as the zero locus of a section of a vector bundle over $\quot_d$. Therefore, the virtual fundamental class exists and is given by the Euler class of the vector bundle. The following well-known Propositions explain the details.
	\begin{prop}\label{prop:quot_smooth}
		For any choice of $N$, $r$ and $d$, $\quot_d(\mathbb{C}^{N},r,\mathbb{P}^1)$ is smooth.
	\end{prop}
	\begin{proof}
		The deformation theory for Quot schemes is given by $\text{Ext}^\bullet(S,Q) $. Since we work over curves it is enough to show that $\text{Ext}^1(S,Q)=0$. Using Serre duality, $\text{Ext}^1(S,Q)=\text{Ext}^0(Q,S(-2))^\vee$. Since $\mathbb{C}^{N}\otimes\cO\to Q$ is a surjection and $S\to \mathbb{C}^{N}\times\cO$ is an injection, it is enough to show that $\Hom(\mathbb{C}^{N}\otimes\cO, \mathbb{C}^{N}\otimes\cO(-2))=0$, which is clear. 
	\end{proof}
	
	\begin{prop}
		Let $\pi:\quot_d\times \mathbb{P}^1\to \quot_d$ be the projection. Then $\pi_* (\wedge^2\cS^{\vee})$ is a locally free sheaf.
	\end{prop}
	\begin{proof}
		Note that for any point $q=[0\to S\to \cO^{N}\to Q\to 0]$ in the Quot scheme, $\mathbb{C}^N\otimes\cO\to S^{\vee}$ is generically surjective and so is \[\phi:\wedge^2( \mathbb{C}^{N}\otimes \cO)\to \wedge^2 S^{\vee}.\] Observe that $\wedge^2(\mathbb{C}^N\otimes \cO)= \mathbb{C}^{\binom{N}{2}}\otimes \cO$.
		We have the following exact sequences of sheaves
		\begin{center}
			\begin{tikzcd}[column sep=0cm, row sep=0cm]
				0&\to&\ker \phi &\to&  \mathbb{C}^{\binom{N}{2}}\otimes \cO&\to& \text{im}\,\phi&\to& 0\\
				0&\to&\text{im}\,\phi &\to& \wedge^2S^{\vee}&\to& \coker\phi&\to& 0
			\end{tikzcd}.
		\end{center}
		
		Since $\coker(\phi)$ is zero dimensional and $ \mathbb{C}^{\binom{N}{2}}\otimes \cO$ is a trivial vector bundle over $\mathbb{P}^1$, their first sheaf cohomology groups vanish.	
		The first exact sequence implies $H^1(\mathbb{P}^1,\text{im}\,\phi)=0$. The second exact sequence gives us  $H^1(\mathbb{P}^1,\wedge^2(S^{\vee}))=0$, hence $h^0(\wedge^2 S^{\vee})=\chi(\wedge^2 S^{\vee})$ is constant. Using Grauert's theorem we conclude that $\pi_*(\wedge^2(\cS^{\vee}))$ is locally free.
	\end{proof}
	
	The symplectic form $\sigma: \wedge^2(\mathbb{C}^N\otimes\cO)\to \cO$ induces an element of $H^0(\mathbb{P}^1,\wedge^2S^\vee)$ given as the composition \[\wedge^2 S\to\wedge^2\mathbb{C}^N\otimes\cO \xrightarrow{\sigma} \cO \] for any subsheaf $S$ of $\mathbb{C}^N\otimes\cO$. This induces a section, denoted as $\tilde{\sigma}$, of $\pi_*(\wedge^2\cS^\vee)$ over $\quot_d$.
	
	Recall that $\IQ_d$ is the subscheme of $\quot_d$ consisting of subsheaves $S$ of $\mathbb{C}^N\otimes\cO$ such that the above composition is zero, hence $\IQ_d=\text{Zero}(\tilde{\sigma})$. Therefore, we have a natural perfect obstruction theory and a virtual fundamental class proving Theorem \ref{thm:POT} in this case.

	\subsection{The Perfect Obstruction theory in general}
	
	In the general case, the two main aspects of the above proof break down, namely $\quot_d$ is not always smooth and the sheaf $\pi_*(\wedge^2\cS^\vee)$ may not be locally free. To construct a perfect obstruction theory, we will have to make a few auxiliary constructions.
	
	Fix $V,L, r$ and $d$. Let $\bun$ be the moduli stack of rank $r$ and degree $d$ vector bundles over $C$. There is a natural forgetful map $\mu: \quot_d\to \bun $ sending the exact sequence $0\to S\to V\to Q\to 0$ to $[S^\vee]\in \bun$. 
	
	We define another stack $\wstack$ which parameterizes pairs $(S,\phi)$, where $S$ is a vector bundle with $S^\vee\in \bun$ and $\phi:\wedge^2S\to L$ is a morphism of sheaves. This also comes equipped with a natural map $\eta:\wstack\to \bun$ sending the pair $(S,\phi)$ to $[S^\vee]$.
	
	We have tabulated the situation in the following commutative diagram 
	\begin{center}
	\begin{tikzcd}
		\IQ_d \arrow[hookrightarrow,r,"i"]\arrow[d, "\mu\circ i"]&\quot_d\arrow[d, "\tilde{\sigma }"]\arrow[ld, "\mu"]\\
		\bun \arrow[hookrightarrow,r , "z"]\arrow[leftarrow,r, bend right, "\eta"]&\wstack 
	\end{tikzcd}.
\end{center}
	Here $\tilde{\sigma} $ is the map sending the short exact sequence $0\to S\to V\to Q\to 0 \in \quot_d$ to the pair $(S,\phi)$, where $\phi$ is the composition $\wedge^2S\to \wedge^2V\xrightarrow{\sigma}L$. 
	
	Recall $\IQ_d$ is precisely the closed locus in $\quot_d$ which is sent to $(S,0)$ under the map $\tilde{\sigma}$. There is a zero section $z:\bun \to \wstack$ sending $[S^\vee]$ to $(S,0)$, and we see that $\IQ_d$ is the fiber product of the maps $\tilde{\sigma}$ and $z$. 
	
	The advantage of the above description is that we understand the cotangent complex of $\quot_d$ and $\bun$, and the new stack $\wstack$ is an abelian cone over $\bun$. We will first describe relative perfect obstruction theory for the maps $\mu$ and $\eta$, and use it to obtain a relative perfect obstruction theory for $\IQ_d$ relative to $\bun$. Since $\bun$ is a smooth Artin stack, this standardly  yields a global perfect obstruction theory for $\IQ_d$, by \cite[Appendix B]{Graber1997LocalizationOV}.

	\subsection{A perfect obstruction theory for $\wstack$}
	We will first carefully define the stack $\wstack$ and show that it is an abelian cone over $\bun$. We will use the results in \cite{scattareggia2018perfect} and \cite{scattareggia_2017} to obtain perfect obstruction theory of $\wstack$ over $\bun$.
	\begin{defn}
		A \coh is a pair $(S,\phi)$ where $S$ is a locally free sheaf on $C$ and $\phi$ is a morphism of sheaves $\phi:\wedge^2S\to L$ over $C$. A family of Wedge systems over a scheme $T$ is $(\pi :C\times T\to T, \cS,\phi : \wedge^2\cS \to p^*L )$	where $p:C\times T\to C$ is the first projection and $\cS $ is a locally free sheaf over $C\times T$.
	\end{defn}
	An isomorphism of two families of \coh $(\pi :C\times T\to T, \cS,\phi : \wedge^2\cS \to p^*L )$ and $(\pi :C\times T\to T, \cS',\phi' : \wedge^2\cS' \to p^*L )$ over $T$ is an isomorphism $\alpha:\cS \to \cS'$ over $C\times T$ such that $ \phi =\phi'\circ \wedge^2\alpha$.
	\begin{defn}
		Let $\wstack$ be the category fibered in groupoids defined by $\wstack (T)$ being the families of Wedge systems over T. Let $\eta:\wstack\to \bun $ be the forgetful morphism.
	\end{defn}
	\begin{prop}\label{prop:S_abelian_cone}
		There is a natural isomorphism of $\bun$-stacks
		\begin{equation}
			\wstack \to \Spec\Sym (\mathbf{R}^1\pi_*(\wedge^2\cS \otimes p^*L^\vee\otimes \omega_\pi)) 
		\end{equation} 
		where $\omega_\pi$ is the relative dualising sheaf of $\pi:\wstack\times C \to\wstack$. In particular $\wstack$ is an abelian cone over $\bun $. Thus $\wstack$ is an algebraic stack.
		\begin{proof}
			The proof is almost same as the proof of Prop 1.8 in \cite{scattareggia2018perfect}. Let $T$ be a scheme, then $\wstack (T)=\{ {t:T\to \bun , \phi :\bar{t}^*\wedge^2\cS\to p^*L }  \}$, where $\bar{t}$ is the induced map from $C\times T\to C\times \bun$. Using Grothendieck duality and base change there is a canonical bijection between $\Hom(\bar{t}^*\wedge^2\cS, p^*L)$ and $\Hom(t^*\mathbf{R}^1\pi_*(\wedge^2\cS \otimes p^*L^\vee\otimes \omega_\pi),\cO_T)$ which is compatible with pull backs.
		\end{proof}
	\end{prop}

	\begin{cor}\label{cor:POT_eta}
		There is a relative perfect obstruction theory for $\eta$ induced by  
		\begin{equation*}
			\mathbf{R}\pi_*(\sHom(\wedge^2\cS,p^*L))^\vee
			\to \tau_{[-1,0]}\mathbb{L}_\eta.
		\end{equation*}
	\end{cor}
	\begin{proof}
		The corollary follows using Lemma \ref{lem:relative_obs_theory} by observing that \[\rHom (\mathbf{R}{\pi}_*(\wedge^2\cS \otimes p^*L^\vee\otimes \omega_\pi[1]),\cO_{\wstack} )\] is isomorphic to $\mathbf{R}\pi_*(\sHom(\wedge^2\cS,p^*L))$ in the derived category.
	\end{proof}

	\begin{lem}\label{lem:relative_obs_theory}
		Let $\pi:Y'\to Y$ be a relative dimension one, flat, projective morphism of algebraic stacks and let $F\in Coh(Y')$ be flat over $Y$, then the abelian cone $\wstack:=\Spec\Sym (\mathbf{R}^1\pi_* F)\xrightarrow{\eta} Y$ has a relative perfect obstruction theory induced by the canonical morphism
		\begin{equation}
			\mathbf{R}\bar{\pi}_*(\bar{F}[1])\to \tau_{[-1,0]}\mathbb{L}_\eta
		\end{equation}
		where $\bar{\pi}: Y'\times_Y \wstack\to \wstack$ and $\bar{F}$ is the induced sheaf on $Y'\times_Y \wstack$.
	\end{lem}
	\begin{proof}
		We will briefly explain the argument assuming $Y$ is a scheme. The complete proof is exactly the same as the proof of Proposition 2.4 in \cite{scattareggia2018perfect}.
		
		Under the given conditions, $F$ can be shown to admit a resolution
		\[0\to K\to M\to F\to 0 \] where $M$ is locally free, $\pi_*K=\pi_*M=0$ and the first derived pushforwards $\mathbf{R}^1\pi_*M$ and $\mathbf{R}^1\pi_*K$ are locally free. Then $\eta$ admits a factorization 
		\[ \wstack \xrightarrow{i} \Spec\Sym(\mathbf{R}^1\pi_*M)\xrightarrow{q} Y  \]
		where $\eta=q\circ i$, $q$ is a smooth morphism and $i$ is a closed embedding. Then $\tau_{[-1,0]}\mathbb{L}_\eta\cong [I|_{\wstack}\to \Omega_q |_{\wstack}]$, where $I$ is the ideal sheaf of $i$. There is a natural isomorphism $\eta^*\mathbf{R}^1\pi_*M\to \Omega_q |_{\wstack} $ and surjection $\eta^*\mathbf{R}^1\pi_*K \to I|_{\wstack}  $. 
		
		Therefore, it remains to show that $[\eta^*\mathbf{R}^1\pi_*K \to \eta^*\mathbf{R}^1\pi_*M ]$ is quasi-isomorphic to $	\mathbf{R}\bar{\pi}_*(\bar{F}[1])$. By cohomology and base-change, $[\eta^*\mathbf{R}^1\pi_*K \to \eta^*\mathbf{R}^1\pi_*M ]$ is isomorphic to $[\mathbf{R}^1\bar{\pi}_*\bar{\eta}^*\bar{K} \to \mathbf{R}^1\bar{\pi}_*\bar{\eta}^*\bar{M} ]$, where \[0\to\bar{ K}\to \bar{M}\to \bar{F}\to 0\] is the induced resolution on $Y'\times_Y \wstack$. The required statement is obtained by the distinguished triangle of the above short exact sequence.
	\end{proof}

	\subsection{Perfect Obstruction theory}
	Recall that we have a map $\tilde{\sigma}:\quot_d \to \wstack$ which takes a subsheaf $[0\to S\to V\to Q\to 0]$ to the point $(S,\phi)$ in $\wstack$ where $\phi$  is the composition of $\wedge^2S\to \wedge^2V\to L$. This can be defined as a morphism of $\bun$-stacks. 

	Consider the morphisms
	\begin{center}
		\begin{tikzcd}
			\quot  \arrow[r, "\tilde{\sigma}"]
			& \wstack\arrow[r,"\eta "] &\bun.
		\end{tikzcd}
	\end{center}
	Let $\mu =\eta \circ \tilde{\sigma}$. There exists a distinguished triangle 
	\begin{equation}\label{eq:dist_trianle1}
		\tilde{\sigma}^*\mathbb{L}_\eta \to \mathbb{L}_\mu \to \mathbb{L}_{\tilde{\sigma}}\to \tilde{\sigma}^*\mathbb{L}_\eta[1].
	\end{equation}
	
	Note that the Quot schemes over smooth curves have perfect obstruction theories as described in \cite{marian2007}. In order to obtain the relative perfect obstruction theory over $\bun$, we consider $\quot_d $ as an open substack of the abelian cone \[\Spec\Sym (\mathbf{R}^1\pi_*(\cS \otimes p^*V^\vee \otimes\omega_\pi)).\] Therefore Lemma \ref{lem:relative_obs_theory} and relative duality implies that the morphism 
	\begin{equation*}
		\mathbf{R}\pi_*(\sHom(\cS,p^*V))^\vee
		\to \tau_{[-1,0]}\mathbb{L}_\mu
	\end{equation*}
	induces a perfect obstruction theory for $\mu : \quot_d \to \bun$. We also recall Corollary \ref{cor:POT_eta}. Thus we get a map of distinguished triangles completing \eqref{eq:dist_trianle1} by the axioms of derived category:

		\begin{equation}
			\begin{tikzcd}\label{eq:dist_triangle}
			\mathbf{R}\pi_*(\sHom(\wedge^2\cS,p^*L))^\vee\arrow[r,"d\tilde{\sigma}"]\arrow[d]& \mathbf{R}\pi_*(\sHom(\cS,p^*V))^\vee \arrow[r]\arrow[d]&\mathbf{R}\pi_*(D^{\bullet})^\vee\arrow[d] 
			\\
			\tilde{\sigma}^*\mathbb{L}_\eta \arrow[r]& \mathbb{L}_\mu \arrow[r]& \mathbb{L}_{\tilde{\sigma}}
	\end{tikzcd}
		\end{equation}
	where $D^\bullet=[\sHom(\cS,p^*V)\xrightarrow{d\sigma} \sHom(\wedge^2\cS,p^*L)]$. The description of $d\sigma$, given below, is important for proving Lemma \ref{lem:generic_surjectivity}.
	
	Fix a vector bundle $S$ in $\bun$, then the map $\tilde{\sigma}$ restricts to a quadratic map $\Hom (S,V)\to \Hom (\wedge^2S,L)$ sending $f$ to $\sigma\circ\wedge^2f$. Vanishing of this map is precisely the locus of the fiber of $\IQ_d$ over $S$. 
	Hence the tangent space at a point $f=[0\to S\xrightarrow{f} V\to Q\to 0]$ in $\IQ_d$ relative to $\bun$  is given as kernel of the linear map $d\tilde{\sigma}:\Hom (S,V)\to \Hom (\wedge^2S,L)$ sending $g$ to the map $[u\wedge v\to \sigma(f(u)\wedge g(v)+g(u)\wedge f(v))]$. The corresponding map of sheaves $d\sigma:\sHom(S,V)\to \sHom(\wedge^2S,L) $ over the fiber $C\times \{f\}$ is given by the same expression over each open sets of $C$.
	
	Over $C\times \IQ_d$ we have the universal section $f$ of the vector bundle $\sHom(\cS,p^*V)$. The above description induces a morphism of locally free sheaves
	\begin{equation*}
		d\sigma : \sHom(\cS,p^*V)\to \sHom(\wedge^2\cS,p^*L).
	\end{equation*}
	
	We have seen in Proposition \ref{prop:S_abelian_cone} that $\wstack $ is an abelian cone, therefore it comes equipped with the zero section $z:\bun \to \wstack$ which is a closed immersion. Recall that $\IQ_d$ sits inside the commutative diagram
	\begin{center}
		\begin{tikzcd}
			\IQ_d \arrow[hookrightarrow,r,"i"]\arrow[d, "\mu\circ i"]&\quot_d\arrow[d, "\tilde{\sigma }"]\arrow[ld, "\mu"]\\
			\bun \arrow[hookrightarrow,r , "z"]\arrow[leftarrow,r, bend right, "\eta"]&\wstack 
		\end{tikzcd}.
	\end{center}
	
	Observe that $\IQ_d$ is the inverse image $\tilde{\sigma}^{-1}(z(\bun)) $. The perfect obstruction theory $\mathbf{R}\pi_*(D^{\bullet})^{\vee}$ of $\sigma$ induces a perfect obstruction theory of $\IQ_{d}$ relative to $\bun$ using the map of cotangent complex
	\begin{equation}\label{eq:Cotangent_cmlx_IQ/Bun}
		i^*\mathbb{L}_{\tilde{\sigma}}\to \mathbb{L}_{\IQ_d/\bun}.
	\end{equation}
		\begin{lem}\label{lem:rel_obs_theory}
		There is a perfect obstruction theory of $\IQ_d$ relative to $\bun$ induced by
		\begin{equation}\label{eq:IQ_rel_Obst}
			\mathbf{R}\pi_*(D^\bullet)^\vee \to \tau_{[-1,0]}\mathbb{L}_{\IQ_d/\bun}.
		\end{equation}
		where $D^\bullet=[\sHom(\cS,p^*V)\xrightarrow{d\sigma} \sHom(\wedge^2\cS,p^*L)]$ is the two term complex over vector bundles with amplitude in [0,1] over $C\times \IQ_d$.
	\end{lem}
\begin{proof}
	We obtain the perfect obstruction theory in \eqref{eq:IQ_rel_Obst} by restricting the perfect obstruction theory of $\tilde{\sigma}$ in \eqref{eq:dist_triangle} to $\IQ_d$ using \eqref{eq:Cotangent_cmlx_IQ/Bun}.
	
Let $D^\bullet|_C=[Hom(S,V)\xrightarrow{d\sigma} Hom(\wedge^2S,L)]$ be the restriction to a fibers, denoted as $C$, of $\pi:C\times \IQ_d\to \IQ_d$. Consider the hypercohomology long exact sequence
\begin{equation*}
\cdots \to\Hglob^1(Hom(S,V))\to \Hglob^1(Hom(\wedge^2S,L))\to \mathbb{H}^2(D^\bullet|_C)\to \Hglob^2(Hom(S,V))=0.
\end{equation*}
Since $d\sigma$ is generically surjective (see Lemma \ref{lem:generic_surjectivity}) and $C$ is one dimensional, $\Hglob^1(Hom(S,V))\to \Hglob^1(Hom(\wedge^2S,L))$ is surjective. Thus we conclude that $\mathbb{H}^2(D^\bullet|_C)$ vanishes.
\end{proof}

	\begin{lem}\label{lem:generic_surjectivity}
		The restriction of $d\sigma$ to each fiber $C=C\times \{f\}$, where $[0\to S\xrightarrow{f} V\to Q\to 0]$ is an element in  $\IQ_{d}$, is generically surjective.
	\end{lem}
\begin{proof}
	Note that $f$ is morphism of vector bundle over $C\backslash A$ where $A$ is finite set of points in $C$. We will show that the linear map of vector spaces
	\begin{align*}
		\phi:Hom(S_x\to V_x)&\to Hom(\wedge^2 S_x,L_x) \\
		g&\to [u\wedge v\to \sigma(f(u)\wedge g(v)+g(u)\wedge f(v))]
	\end{align*} 
	is surjective for all $x\in C\backslash A$. This is now an exercise in linear algebra.
	
	Let $N=2n$. We can choose symplectic coordinates $\{e_1,\dots e_N\}$ of $V_x$ such that $\sigma(e_i,e_{n+i})=1$ and $f$ identifies the isotropic subspace $S_x$ with $span \{e_1,\dots ,e_r\}$. An element $g\in Hom(S_x\to V_x)$ can be identified with an $N\times r$ matrix $( B_{i,j} )$. A simple calculation shows that $g\in \ker \phi$ if and only if $B_{i,n+k}=B_{k,n+i}$ for all $1\le i,k\le r$. Thus the rank of $\ker \phi$ is $Nr-\binom{r}{2}$, hence $\phi$ is surjective.
\end{proof}

	\begin{proof}[Proof of Theorem \ref{thm:POT}]
		In Lemma \ref{lem:rel_obs_theory}, we constructed a relative perfect obstruction theory.	We follow the arguments in \cite[Appendix B]{Graber1997LocalizationOV} verbatim to obtain an absolute perfect obstruction theory. Here we use the fact that $\bun$ is a smooth Artin stack with obstruction theory given by $\mathbf{R}\pi_*(\sHom(\cS,\cS))^\vee[-1]\to \mathbb{L}_{\bun}$.
	\end{proof}
	
	\begin{rem}
		We note that when $V$ and $L$ are trivial and $\sigma$ is induced from a standard symplectic or symmetric form on $\mathbb{C}^N$, there is another way to construct the virtual fundamental class for $\IQ_d$ using the theory of quasi-maps to GIT quotients as discussed in \cite{ciocanfontanine2011stable}.
		
		Indeed, $\IQ_d$ can be considered as the moduli space of quasi maps from $C$ to $\SG(r,N)$ (or $\OG(r,N)$). The isotropic Grassmannian can be realized as a GIT quotient of $W\git_{\theta} G$, where $\theta = \det^{-1}$ is the multiplicative character of $G=GL_r$ and $W=\{f\in \Hom(\mathbb{C}^r,\mathbb{C}^N) : \sigma(f(u),f(v))=0 \,\forall \,u,v\in \mathbb{C}^r\}$ is a closed subscheme of the affine space $\Hom(\mathbb{C}^r,\mathbb{C}^N)$.
	\end{rem}

	\section{Symplectic isotropic Quot schemes}\label{sec:Nvir_sympl}
	Throughout this section we will assume that $\om$ is the standard symplectic form on $\mathbb{C}^N\otimes\cO$; i.e., it is induced by the block matrix
	\begin{align*}
		\om=\begin{bmatrix}
			0&I_n\\
			-I_n&0
		\end{bmatrix}
	\end{align*}
	where $N=2n$. 
	
	There is a natural action of $Sp(2n)$ on $\IQ_d$ induced by the respective action on $ \mathbb{C}^{2n}$. We consider the subtorus $G=\mathbb{C}^*\subseteq Sp(2n) $ given by $(t^{-w_1},\dots , t^{-w_N})$ where $w_i=-w_{i+n}$ for $1\le i\le n$. The weights $w_i$ are assumed to be distinct, unless stated otherwise.
	
	\subsection{Fixed Loci} \label{sec:fixed_loci}
	
	Each summand $\cO$ of $\mathbb{C}^N\otimes\cO$ is acted upon with different weights. A point $[0\to S\to \mathbb{C}^N\otimes\cO\to Q\to 0]$ in $\IQ_d$ is fixed under the action of $G$ if and only if :
	\begin{itemize}
		\item[(i)] $S$ splits as a direct sum of line bundles \[S=\oplus_{j=1}^rL_j,\] where $L_j$ is subsheaf of one of the $N$ copies of $\cO$ of $\mathbb{C}^N\otimes\cO$. Denote $k_j$ by the position of this copy of $\cO$. 
		\item[(ii)] $k_j-k_i\not \equiv 0 \mod n$ for any $1\le i<j\le r$ : This ensures that $S$ is isotropic.
	\end{itemize}
	
	Let $\underline{k}=\{k_1,\dots,k_r\}$ and $\vec{d}=(d_1,\dots,d_r)$ where $d_j=\deg L_j$ and \[d_1+\cdots +d_r=d. \] We require $\{i,i+n\}\not \subset  \underline{k}$ for any $1\le i \le n$. Let $\fix_{\vec{d},\underline{k}}$ be the set of fixed points with the numerical data $\vec{d}$ and $\underline{k}$. Note that there are $2^r\binom{n}{r}$ possible values of $\underline{k}$ and $\binom{d+r-1}{r-1}$ choices of $\vec{d}$.
	
	Denote $\cO_{k_i}$ be the $k_i$'th copy of $\cO$ in $\mathbb{C}^N\otimes\cO$. The short exact sequence \[0\to L_i\to \cO_{k_i}\to T_i\to 0\] defines an element of $C^{[d_i]}$, the Hilbert scheme of $d_i$ points on $C$. Therefore we have 
	\begin{align*}
		\fix_{\vec{d},\underline{k}}=C^{[d_1]}\times C^{[d_2]}\times \dots \times C^{[d_r]}.
	\end{align*}
	
	\subsection{The Equivariant Normal bundle}
	
	Let $0\to \mathcal{L}_{i}\to\cO_{k_i}\to \mathcal{T}_{i} \to0$ be the universal exact sequence over $C\times C^{[d_i]}$. We use the same notation for the pull-back exact sequence over $C\times \fix_{\vec{d}, \underline{k}}$.
	
	Let $0\to \cS\to \mathbb{C}^{N}\otimes\cO\to \cQ\to 0$ be the universal exact sequence over $C\times \IQ_d$. This restricts to 
	\begin{align*}
		0\to \mathcal{L}_{1}\oplus \dots \oplus \mathcal{L}_{r}\to \mathbb{C}^{N} \otimes \cO\to \mathcal{T}_{1}\oplus \dots \oplus \mathcal{T}_{r} \oplus \mathbb{C}^{N-r} \otimes \cO\to 0  
	\end{align*} 
	on $C\times \fix_{\vec{d},\underline{k}}$.
	
	 Let $\pi_!$ be the derived pushforward  $\mathbf{R}^0\pi_*-\mathbf{R}^1\pi_*$ in the K-theory. Recall that in Theorem \ref{thm:POT}, we provided a perfect obstruction theory for the isotropic Quot scheme. In the $K$-theory of $\IQ_d$, the corresponding virtual tangent bundle is given by 
	\begin{align*}
		T^{\vir}=\pi_![(\rHom(\cS,\cQ))]-\pi_![(\sHom(\wedge^2{\cS},\cO))].
	\end{align*}
	
	The restriction of the virtual tangent bundle in the $\mathbb{C}^*$-equivariant $K$-theory of $\fix_{\vec{d},\underline{k}}$ is given  by the following formula
	\begin{align*}
		\pi_!\bigg(\sum_{i,j\in [r]}^{} [\mathcal{L}_i^\vee\otimes\mathcal{T}_j ]+\sum_{i\in [r],k\in \underline{k}^c}[\mathcal{L}^\vee_i]-\sum_{1\le i<j\le r}^{} [\mathcal{L}^\vee_i\otimes\mathcal{L}^\vee_j]\bigg),
	\end{align*}
	where the above three groups of elements have $\mathbb{C}^*$ weights $(w_{k_i}-w_{k_j})$, $(w_{k_i}-w_{k})$ and $(w_{k_i}+w_{k_j})$ respectively. 
	
	Note that the fixed part of the restriction of $T^{\vir}$ to $\fix_{\vec{d},\underline{k}}$ is \[\sum_{i\in [r]}^{}\pi_![\mathcal{L}^\vee_i\otimes\mathcal{T}_i],
	\]
	which matches the tangent bundle of $\fix_{\vec{d},\underline{k}}$. The induced virtual class $[\fix_{\vec{d},\underline{k}}]^{\vir}=[\fix_{\vec{d},\underline{k}}]$ agrees with the usual fundamental class.
	
	The virtual equivariant normal bundle $\nbun$ is given by the moving part of the restriction of $T^{\vir}$.  Using the identity in $K$-theory, \[[\mathcal{L}_i^\vee\otimes\mathcal{T}_j ]=[\mathcal{L}_i^\vee\otimes\cO_{k_j} ]-[\mathcal{L}_i^\vee\otimes\mathcal{L}_j ],\]
	we obtain the following equality
	\begin{align}\label{eq:Nvir_sympl}
		\nbun=	\pi_!\bigg(\sum_{\substack{i\in [r],k\in [N]\\ k\ne k_i } }[\mathcal{L}^\vee_i]-\sum_{\substack{i,j\in [r] \\ i\ne j}}^{} [\mathcal{L}_i^\vee\otimes\mathcal{L}_j ]-\sum_{1\le i<j\le r}^{} [\mathcal{L}^\vee_i\otimes\mathcal{L}^\vee_j]\bigg),
	\end{align}
	where the terms are acted on with wights $(w_{k_i}-w_{k})$, $(w_{k_i}-w_{k_j})$ and $(w_{k_i}+w_{k_j})$ respectively.
	
	\subsection{Chern polynomials} 
	In the subsection we briefly describe certain Grothendieck-Riemann-Roch calculations for the map $\pi:C\times X\to X$, where \[X=C^{[d_1]}\times C^{[d_2]}\times \dots \times C^{[d_r]}.\]
	
	Let $\{1,\delta_1,\dots, \delta_{2g}, \omega \}$ be the symplectic basis for the cohomology ring of $C$ with the relations $\delta_i\delta_{i+g}=\omega=-\delta_{i+g}\delta_i$ for all $1\le i\le g$. Consider the K\"unneth decomposition of the cohomology classes $c_1(\mathcal{L}^\vee)$ in $C\times C^{[d_i]}$ with respect to a chosen symplectic basis of $\Hglob^*(C)$,
	\begin{equation}\label{eq:def_x,y_classes}
		c_1(\mathcal{L}_i^\vee)=x_i\otimes 1 + \sum_{k=1}^{2g}y_i^k\otimes\delta_k+ d_i\otimes \omega.
	\end{equation}
	The theta class, $\theta_i \in \Hglob^*(C^{[d_i]})$, is the pullback of the usual theta class under the map
	\begin{align*}
		C^{[d_i]}\to \text{Pic}^{d_i}.
	\end{align*}
	We have the following relation (explained in \cite{ACGH}) \[\bigg(\sum_{k=1}^{2g}(y_i^k\otimes \delta_k) \bigg)^2=-2\theta_i\otimes \omega. \]
	We will use the same notation for the pullback of $x_i, y_i^{k}$ and $\theta_i$ under the map  \begin{equation*}
		pr_{i}: X \to  C^{d_i}.
	\end{equation*}
	
	Let $E$ be a vector bundle of rank $m$ and let $c_t(E)=1+c_1(E)t+\cdots+c_m(E)t^m$ be its Chern polynomial. We extend the definition of $c_t$ to the $K$-theory in the usual way. We can use Grothendieck-Riemann-Roch to obtain expression for the Chern polynomials $c_t(\pi_![\mathcal{L}_i^\vee])$, $c_t(\pi_![\mathcal{L}_i^\vee\otimes\mathcal{L}_j])$  and $c_t(\pi_![\mathcal{L}_i^\vee\otimes\mathcal{L}_j^\vee])$:
	\begin{align}\label{GRR}
		c_t(\pi_![\mathcal{L}_i^\vee])&= (1+tx_i)^{d_i-\bar{g}}e^{-\frac{t\theta_i}{(1+tx_i)}}   \\
		c_t(\pi_![\mathcal{L}_i^\vee\otimes\mathcal{L}_j])&=(1+t(x_i-x_j))^{d_i-d_j-\bar{g}} e^{-\frac{t(\theta_i+\theta_j+\phi_{ij})}{1+t(x_i-x_j)}}\nonumber \\
		c_t(\pi_![\mathcal{L}_i^\vee\otimes\mathcal{L}_j^\vee])&=(1+t(x_i+x_j))^{d_i+d_j-\bar{g}} e^{-\frac{t(\theta_i+\theta_j-\phi_{ij})}{1+t(x_i+x_j)}} \nonumber\\
		c_t(\pi_![\mathcal{L}_i^\vee\otimes\mathcal{L}_i^{\vee}])&=(1+2tx_i)^{2d_i-\bar{g}} e^{-\frac{4t\theta_i}{1+2tx_i}}\nonumber
	\end{align}
	where $\phi_{ij}=-\sum_{k=1}^{g}(y_i^ky_j^{k+g}+y_j^ky_i^{k+g})$. The detailed calculation for the first two expression can be found in \cite{ACGH} and \cite{marian2007}. The other two expressions are obtained in a similar way. We will briefly explain the last one for completeness: The first Chern class is $c_1(\mathcal{L}^\vee \otimes \mathcal{L}^\vee)= 2c_1(\mathcal{L}^\vee)$, therefore the Chern character  \[ch(\mathcal{L}^\vee\otimes\mathcal{L}^\vee) = e^{2x}\otimes 1 + e^{2x}(2d- 4\theta)\otimes \omega + 2\sum_{k}^{}y^k\otimes \delta_k. \]  
	We may further apply Grothendieck Riemann Roch to obtain the Chern characters of $\pi_![\mathcal{L}^\vee\otimes \mathcal{L}^\vee]$ and then covert it into Chern polynomials to obtain the required result. The Chern character is 
	\begin{align*}
		ch(\pi_![\mathcal{L}^\vee\otimes \mathcal{L}^\vee]) &= \pi_* (ch(\mathcal{L}^\vee\otimes \mathcal{L}^\vee)(1+(1-g)\omega))\\
		&= e^{2x}(2d+(1-g)-4\theta).
	\end{align*}

	\subsection{The Euler class of virtual normal bundle}\label{sec:Euler_class_sympl}
	Next we would like to find the equivariant Euler class of $\nbun$ in the equivariant cohomology ring $H^*(\fix_{\vec{d},\underline{k}})[[t,t^{-1}]]$. This will be useful in the virtual localization formula.
	
	Let $E$ be one of the line bundles appearing in the formula for $\nbun$ in \eqref{eq:Nvir_sympl}. We evaluated the formula for the total Chern classes $c_q(\pi_!E)$ in \eqref{GRR}. Let $\pi_!E$ be acted on with weight $w$, then the equivariant Euler class is a homogeneous element in $H^*(\fix_{\vec{d},\underline{k}})[t,t^{-1}]$ and is given by
	\begin{equation*}
		e_{\mathbb{C}^*}(\pi_!E)=(wt)^mc_{\frac{1}{wt}}(\pi_!E)
	\end{equation*}
	where $m=\chi(\pi_!E)$ is the virtual rank.
	
	Consider the polynomial $P(X)=\prod_{i=1}^{N}(X-w_it)$. Let $Y_i=x_i+w_{k_i}t$ be a change of variable over $\mathbb{C}[[t]]$. Then 
	\begin{align}\label{eq:euler_class_computation_1}
		\prod_{\substack{i\in [r],k\in [N]\\ k\ne k_i } }\frac{1}{e_{\mathbb{C}^*}(\pi_![\mathcal{L}^\vee_i])}
		&= \prod_{\substack{i\in [r],k\in [N]\\ k\ne k_i } }(Y_i-w_kt)^{-d_i+\bar{g}}e^{\frac{\theta_i}{(Y_i-w_kt)}} \\
		&= \prod_{i\in [r]}^{}\bigg(\frac{P(Y_i)}{x_i}\bigg)^{-d_i+\bar{g}}e^{\theta_i\big(\frac{P'(Y_i)}{P(Y_i)}-\frac{1}{x_i} \big) }\nonumber 
	\end{align}
	Here we are using the elementary identity \[\frac{P'(X)}{P(X)}=\sum_{k=1}^{N}\frac{1}{X-w_kt}. \] For the remaining classes, we obtain
	\begin{align}\label{eq:euler_class_computation_2}
		\prod_{\substack{i,j\in [r] \\ i\ne j}}^{} e_{\mathbb{C}^*}(\pi_![\mathcal{L}_i^\vee\otimes\mathcal{L}_j ])
		&= \prod_{\substack{i,j\in [r] \\ i\ne j}}^{} (Y_i-Y_j)^{d_i-d_j-\bar{g}} e^{-\frac{(\theta_i+\theta_j+\phi_{ij})}{Y_i-Y_j}}\\
		&= (-1)^{\bar{g}\binom{r}{2}+d(r-1)} \prod_{i< j}^{}(Y_i-Y_j)^{-2\bar{g}} \nonumber 
	\end{align}
	\begin{align}\label{eq:euler_class_computation_3}
		\prod_{\substack{i,j\in [r] \\ i< j}}^{} e_{\mathbb{C}^*}(\pi_![\mathcal{L}^\vee_i\otimes\mathcal{L}^\vee_j])	&= \prod_{ i<j}^{}(Y_i+Y_j)^{d_i+d_j-\bar{g}} e^{-\frac{(\theta_i+\theta_j-\phi_{ij})}{Y_i+Y_j}}
	\end{align}
	Using the multiplicative property for the Euler classes, we have the the following expression for the equivariant Euler class of the virtual normal bundle :
	\begin{align}\label{eq:e(Nvir)spl}
		\frac{1}{e_{\mathbb{C}^*}(\nbun)}=u\prod_{i}^{}h_i^{d_i-\bar{g}} e^{\theta_iz_i}
		\cdot \prod_{i< j}\frac{(Y_i+Y_j)^{d_i+d_j-\bar{g}}}{(Y_i-Y_j)^{2\bar{g}}}e^{-\frac{\theta_i+\theta_j-\phi_{ij}}{Y_i+Y_j}}
	\end{align}
	where $u=(-1)^{\bar{g}\binom{r}{2}+d(r-1)}$, $h_i=\frac{x_i}{P(Y_i)}$ and 
	\begin{equation}\label{eq:def_z_i}
		z_i=\bigg(\frac{P'(Y_i)}{P(Y_i)}-\frac{1}{x_i} \bigg).
	\end{equation}
	
	\section{Symmetric isotropic Quot scheme}\label{sec:Nvir_symm}
	
	Throughout this section we will assume $N=2n$, $V=\mathbb{C}^N\otimes{\cO}$ is the trivial vector bundle over $C$ and $\sigma$ is induced by a non-degenerate symmetric form on $\mathbb{C}^N$. We may assume that the symmetric form $\sigma$ is given by the block matrix 
	\begin{align*}
		\sigma =\begin{bmatrix}
			0&I_n\\
			I_n&0
		\end{bmatrix}.
	\end{align*}

	There is a natural action of $SO(N)$ on the $\IQ_d$ induced by the respective action on $ \mathbb{C}^N$. The subtorus $G=\mathbb{C}^*\subset SO(N) $ given by $(t^{-w_1},\dots , t^{-w_N})$ also acts on $\IQ_{d}$ where the weights $w_i=-w_{i+n}$ for $1\le i\le n$.
	
	\subsection{Fixed Loci} 	
	When the weights are distinct, we get the same description of fixed loci as in the case of $\sigma$ symplectic. Thus the fixed loci of the $\mathbb{C}^*$ action are isomorphic to a disjoint union of 
	\[\fix_{\vec{d},\underline{k}}=C^{[d_1]}\times C^{[d_2]}\times \dots \times C^{[d_r]}\]
	for each possible tuple of positive integers $\vec{d}=(d_1,d_2,\dots,d_r)$ such that $d_1+d_2+\dots +d_r=d$ and $\underline{k}=\{k_1,\dots , k_r\}\subset \{1,\dots , N\}$ such that $\{i,i+n\}\not \subset  \underline{k}$ for any $1\le i \le n$. 
	
	We will use the localization formula with distinct weights to show compatibility of the virtual fundamental classes in Theorem \ref{compatibiliy_theorem}. We will use non-distinct weights to obtain the Vafa-Intriligator type formula in Theorem \ref{thm:r=2_symmetric}. In the latter case, we will obtain different fixed loci; we will describe it in Section \ref{sec:Int_a_classes}. The description of the equivariant normal bundle will be crucial in proving both the theorems.

	\subsection{Equivariant Normal bundle}\label{sec:Euler_class_symm}
	Let $0\to \cS\to \mathbb{C}^{N}\otimes\cO\to \cQ\to 0$ be the universal exact sequence over $C\times \IQ_d$. This restricts to 
	\begin{align*}
		0\to \mathcal{L}_{1}\oplus \dots \oplus \mathcal{L}_{r}\to\mathbb{C}^{N} \otimes \cO \to \mathcal{T}_{1}\oplus \dots \oplus \mathcal{T}_{r} \oplus \mathbb{C}^{N-r} \otimes\cO\to 0  
	\end{align*} 
	on $C\times \fix_{\vec{d},\underline{k}}$, where $0\to \mathcal{L}_{i}\to\cO\to \mathcal{T}_{i} \to0$ is the universal exact sequence over $C\times C^{[d_i]}$ at the position $k_i$.
	
	Recall that in Theorem \ref{thm:POT}, we provided a perfect obstruction theory for the isotropic Quot scheme. In the $K$-theory of $\IQ_d$, the corresponding virtual tangent bundle is given by 
	\begin{align*}
		T^{\vir}=\pi_![(\rHom(\cS,\cQ))]-\pi_![(\sHom(\Sym^2{\cS},\cO))].
	\end{align*}
	
	The restriction of the virtual tangent bundle in the $\mathbb{C}^*$ equivariant $K$-theory of $\fix_{\vec{d},\underline{k}}$ is given by 
	\begin{align*}
		\pi_!\bigg(\sum_{i,j\in [r]}^{} [\mathcal{L}_i^\vee\otimes\mathcal{T}_j ]+\sum_{i \in [r],k\in \underline{k}^c}[\mathcal{L}^\vee_i]-\sum_{1 \le i\le j\le r}^{} [\mathcal{L}^\vee_i\otimes\mathcal{L}^\vee_j]\bigg).
	\end{align*}
	where the above three summands have $\mathbb{C}^*$ weights $(w_{k_i}-w_{k_j})$, $(w_{k_i}-w_{k})$ and $(w_{k_i}+w_{k_j})$ respectively. 
	
	The fixed part of the restriction of $T^{\vir}$ to $\fix_{\vec{d},\underline{k}}$ is \[\sum_{i\in \underline{k}}^{}\pi_![\mathcal{L}^\vee_i\otimes\mathcal{T}_i]
	\]
	which matches with the $K$-theory class of the tangent bundle of $\fix_{\vec{d},\underline{k}}$. 
	
	The virtual normal bundle $\nbun$ is given by the moving part of the restriction of $T^{\vir}$. In the $K$-theory of $\fix_{\vec{d}}$, 
		\begin{equation}
		\nbun=\pi_!\bigg(\sum_{\substack{i\in[r]\\ k\ne k_i}}[\mathcal{L}^\vee_i]-\sum_{\substack{i,j\in [r]\\i\ne j}}^{} [\mathcal{L}_i^\vee\otimes\mathcal{L}_j ]-\sum_{1\le i\le j\le r}^{} [\mathcal{L}^\vee_i\otimes\mathcal{L}^\vee_j]\bigg).
	\end{equation}
	 Next we would like to determine the equivariant Euler class of $\nbun$ in the equivariant cohomology ring $H^*(\fix_{\vec{d},\underline{k}})[t,t^{-1}]$.  
	 
	 Let $P(X)=\prod_{k=1}^{N}(X-w_{k}t)$ and $Y_i=x_i+w_{k_i}t$. Using \eqref{eq:euler_class_computation_1}, \eqref{eq:euler_class_computation_2} and \eqref{eq:euler_class_computation_3}and the identity
	\begin{align}
		\prod_{i\in [r]}^{} e_{\mathbb{C}^*}(\pi_![\mathcal{L}^\vee_i\otimes\mathcal{L}^\vee_i])	&= \prod_{ i\in [r]}^{}(2Y_i)^{2d_i-\bar{g}} e^{-\frac{2\theta_i}{Y_i}},
	\end{align}
	we obtain the expression for the equivariant Euler class of $\nbun$:
	\begin{align}\label{eq:e(N)_symmetric}
		\frac{1}{e_{\mathbb{C}^*}(\nbun)}=&u2^{2d-r\bar{g}} \prod_{i=1}^{r}h_i^{d_i-\bar{g}}Y_i^{2d_i-\bar{g}}e^{\theta_iz_i}\prod_{i< j}\frac{(Y_i+Y_j)^{d_i+d_j-\bar{g}}}{(Y_i-Y_j)^{2\bar{g}}}e^{-\frac{\theta_i+\theta_j-\phi_{ij}}{Y_i+Y_j}}
	\end{align}
	where $u=(-1)^{d(r-1)+\binom{r}{2}\bar{g}}$ and
\begin{equation}
	\begin{aligned}
		h_i&=\frac{x_i}{P(Y_i)}\\
		z_i&=\frac{P'(Y_i)}{P(Y_i)}-\frac{2}{Y_i}-\frac{1}{x_i}.
	\end{aligned}
\end{equation}

	\section{Compatibility of virtual fundamental classes}
	In this section we only consider $\IQ_d$ with $V, L$ trivial and $N$ even. Fix a point $q\in C$. Then there is a natural embedding 
	\begin{equation}
		i_q:\IQ_d\to \IQ_{d+r}
	\end{equation}
	which sends a subsheaf $S\subset \mathbb{C}^N\otimes \cO$ to the composition $S(-q)\to S \to \mathbb{C}^N\otimes \cO $. Observe that $S(-q)$ is an isotropic subsheaf because the composition
	\begin{equation*}
		S(-q)\to S\to \mathbb{C}^N\otimes \cO \xrightarrow{\sigma} \mathbb{C}^N\otimes \cO\to S^\vee\to S(-q)^{\vee}
	\end{equation*}
	is zero. 
	
	\begin{proof}[Proof of Theorem \ref{compatibiliy_theorem}]	
		We work with the symmetric isotropic Quot scheme. The argument in the symplectic case is similar. 
		
		Let $j$ be the inclusion of the fixed loci into $\IQ_{d}$. The virtual localization formula \cite{Graber1997LocalizationOV} asserts that
		\begin{align*}
			[\IQ_{d}]^{\vir}=j_*\sum_{\vec{d},\underline{k}}\frac{[\fix_{\vec{d},\underline{k}}]^{\vir}}{e_{\mathbb{C}^*}(\nbun)}
		\end{align*}
		in $A_{*}^{\mathbb{C}^*}(\IQ_{d})\otimes \mathbb{Q}[t,t^{-1}]$ where $t$ is the generator of the equivariant ring of $\mathbb{C}^*$. Note that $[\fix_{\vec{d}, \underline{k}}]^{\vir}=[\fix_{\vec{d}, \underline{k}}]$ in our case. We will show the compatibility of the virtual fundamental classes by equating the fixed loci contributions.
		
		We denote $\bar{\fix}=\fix_{\vec{d}+(1,\dots,1),\underline{k}}$ and $\fix= \fix_{\vec{d},\underline{k}}$ for notational convenience. These are fixed loci on $\IQ_{d}$ and $\IQ_{d+r}$ respectively. 
		
		The map $i_q$ restricts to the natural map over the fixed locus $\tilde{i}_q: \fix \to\bar{\fix}$. This sends the fixed point $L_1\oplus \cdots \oplus L_r\subset \mathbb{C}^N\otimes\cO$ to $L_1(-q)\oplus \cdots \oplus L_r(-q)\subset \mathbb{C}^N\otimes\cO$. We have the identity (see \cite{marian2007} for more details)
		\begin{align*}
			\tilde{i}_{q*}[\fix]=\prod_{\ell=1}^{r}\bar{x}_i\cap [\bar{\fix}],
		\end{align*}
		where $\bar{x}_i$ are the cohomology classes on $\bar{\fix}$ defined in \eqref{eq:def_x,y_classes}. 
		
		In the equivariant cohomology of the fixed loci $\fix$, 
		\begin{align*}
			c_{\text{top}}(\Sym^2\cS^{\vee}_q)|_{\fix}& = \prod_{ 1\le i\le j\le r}(Y_i+Y_j)
		\end{align*}
		where $Y_i=x_i+w_{k_i}t$, and over $\bar{\fix}$ we have
		\begin{equation}\label{eq:compatibility_proof}
			c_{\text{top}}(\sHom(\cS_q,\mathbb{C}^N\otimes\cO))|_{\bar{\fix}}= \prod_{ i=1}^{r}\bar{x}_i\cdot \prod_{ i=1}^{r}\bar{h}_i^{-1}.
		\end{equation}

		Using the description of the Euler class of the equivariant normal bundle in \eqref{eq:e(N)_symmetric}, we have
		
		\begin{align*}
			\prod_{ 1\le i\le j\le r}(Y_i+Y_j)^2\cdot \frac{1}{e_{\mathbb{C}^*}(\nbun_{\fix/\IQ_d})}=\tilde{i}_q^*\prod_{ i=1}^{r}h_i^{-1}\cdot \tilde{i}_q^*\frac{1}{e_{\mathbb{C}^*}(\nbun_{\bar{\fix}/\IQ_{d+r}})}.
		\end{align*}
		Hence the fixed loci contribution matches in the application of equivariant virtual localization in \cite{Graber1997LocalizationOV} to $\IQ_{d+r}$ for the fixed loci of the kind $\bar{\fix}=\fix_{\vec{d},\underline{k}}$ with $d_i>0$ for any $1\le i\le r$ with the corresponding contribution over $\IQ_{d}$. When $d_i=0$ for some $i$, the fixed point contribution vanishes since $\bar{x}_i$ appears in \eqref{eq:compatibility_proof}.
	\end{proof}

	\section{Symmetric powers of curves}\label{sec:prelim_hilb}
	In this section we will describe the intersection theory of the products of symmetric powers of curves  \[X_{\vec{d}}= C^{[{d_1}]}\times\cdots \times C^{[d_r]}.\] This will be needed to obtain the Vafa-Intriligator type formula for the intersection of $a$ and $f$ classes over isotropic Quot schemes.
	
	There are two difficulties in the calculation of the virtual intersection numbers involving the above classes : knowing how to intersect $\theta$, $\phi_{ij}$ and $x$ (defined in section \ref{sec:Euler_class_sympl}), and summing over all the fixed loci. Note that the number of fixed loci increases as $d$ increases. Moreover, the expressions for the Euler class of the virtual normal bundles \eqref{eq:e(Nvir)spl} and \eqref{eq:e(N)_symmetric} over the fixed loci involve many complicated terms.  
	
	We describe techniques to evaluate intersection numbers involving the above terms. For the summation, we will use a beautiful combinatorial technique called multivariate Lagrange-B\"urman formula.
	
	\subsection{Intersection theory of $X_{\vec{d}}$}
	
	The following are some known facts about the $x$, $\theta$ and $y$ classes (see \cite{ACGH} and \cite{Th}) over $C^{[d]}$:
	
	\begin{itemize}
		\item The intersections of $x$ and $\theta $ are given by:
		\begin{equation*}
			\int_{C^{[d]}}^{}\theta^\ell x^{d-\ell}=\begin{cases}
				\frac{g!}{(g-\ell)!}&\ell\le g \\
				0& \ell>g
			\end{cases}.
		\end{equation*}
		In particular, for any polynomial $P$, and $\ell\le g$ 
		\begin{equation}\label{eq:theta_int}
			\int_{C^{[d]}}^{}\theta^\ell P(x)=\frac{g!}{(g-\ell)!}\int_{C^{[d]}}^{} x^\ell P(x).
		\end{equation}
		\item The non-zero integrals in the $y$ classes over $C^{[d]}$ satisfy 
		\begin{itemize}
			\item[(i)] $y^{k}$ appears with exponent at most 1 because these are odd classes.
			\item[(ii)] $y^{k}$ appears if and only if $y^{k+g}$ appears.
			\item[(iii)] For any choice of choice of distinct integers $k_1,\dots, k_s\in \{1,\dots g\}$ and a polynomial $P$ in two variables,
			\begin{equation}\label{y_int}
				\int_{C^{[d]}}y^{k_1}y^{k_1+g}\cdots y^{k_s}y^{k_s+g}P(x,\theta)=\frac{(g-s)!}{g!} \int_{C^{[d]}}\theta^sP(x,\theta). 
			\end{equation} 
		\end{itemize}
	\end{itemize}
	
	Fix $\vec{d}=(d_1,\dots,d_r)$ and $X_{\vec{d}}=C^{[d_1]}\times\cdots\times C^{[d_r]}$. For $1\le i\le r$, define the cohomology classes $x_i,y_i^k$ and $\theta_i$ on $X_{\vec{d}}$ obtained by pulling back the corresponding classes from $C^{[d_i]}$. 
	\begin{prop}\label{prop:phi_int}
		Let $P$ be a polynomial in $2r$ variables, then 
		\begin{align}
			\int_{X_{\vec{d}}}^{}\phi_{12}^{2\ell}P(\underline{x},\underline{\theta})= (-1)^\ell\binom{2\ell}{\ell}\binom{g}{\ell}^{-1}\int_{X_{\vec{d}}}^{}(\theta_1\theta_2)^\ell P(\underline{x},\underline{\theta})
		\end{align}
		where $\underline{x}=(x_1,\dots,x_r)$ and $\underline{\theta}=(\theta_1,\dots,\theta_r)$.
	\end{prop}
	\begin{proof}
		Recall that \[\phi_{12}=-\sum_{k=1}^{g}(y_1^ky_2^{k+g}+y_2^ky_1^{k+g}). \]
		For parity reasons, $\phi_{12}$ must appear with even exponent.
		
		Using \eqref{y_int}, $\phi_{12}^{2\ell}$ can be replaced by a constant multiple of $\theta_1^\ell\theta_2^\ell$, where the constant is $\frac{(g-\ell)!^2}{g!^2}$ times the sum of coefficients of \[y_1^{k_1}y_1^{k_1+g}\dots y_1^{k_\ell}y_1^{k_\ell+g}\cdot y_2^{k_1}y_2^{k_1+g}\dots y_2^{k_\ell}y_2^{k_\ell+g}\] in the multinomial expansion of $\phi_{12}^{2\ell}$. We observe that
		\begin{align*}
			(y_1^{k}y_2^{k+g}+y_2^ky_1^{k+g})^2&=y_1^{k}y_2^{k+g}y_2^ky_1^{k+g}+y_2^{k}y_1^{k+g}y_1^ky_2^{k+g}\\
			&=-2y_1^{k}y_1^{k+g}y_2^ky_2^{k+g}. 
		\end{align*} 
		
		Thus the required sum of coefficients is \[(-2)^{\ell}\binom{g}{\ell}\binom{2\ell}{2,\dots, 2},\] where $\binom{g}{\ell}$ is the number of choices for $\{k_{i_1},\dots, k_{i_\ell}\}$ and $\binom{2\ell}{2,\dots, 2}$ is the number of ways of picking $\ell$ pairs of factors in $\phi_{12}^{2\ell}$ each of which contributes $(-2)$. The binomial identity 
		\begin{align}
			(-2)^\ell\binom{g}{\ell}\binom{2\ell}{2,\dots, 2}\frac{(g-\ell)!^2}{(g!)^2}=(-1)^\ell\binom{2\ell}{\ell}\binom{g}{\ell}^{-1}
		\end{align}
		completes the proof. 
	\end{proof}
	
	\subsection{Summing over $|\vec{d}|=d$}
	In Section \ref{sec:Int_a_classes} and \ref{sec:f_intersection}, we will use the localization formula to calculate the tautological intersection numbers. We use the independence of the weights in the localization formula. We will describe how to sum over the fixed point contributions for a special choice of weights. The following two Propositions are crucial for our argument.
	
	Let $w_1,\dots, w_r$ be $r$ distinct $N^{\text{th}}$ roots of unity and let $P(Y)=Y^N-1$.
	\begin{prop}\label{prop:theta_sum_int}
		Let $p_1,\dots, p_r$ and $d$ be non-negative integers and $R(Y_1,\dots ,Y_r)$ be a homogeneous rational function of degree $s=Nd-r\bar{g}(N-1)-p$ where $p_1+\cdots +p_r=p$. Let $B(Y)=\frac{aY^N+b}{Y}$, $Y_i=x_i+w_i$, $h_i=\frac{x_i}{P(Y_i)}$ and \[z_i=\frac{B(Y_i)}{P(Y_i)}-\frac{1}{x_i}. \]
		Then we have the following identity
		\begin{align}\label{Eq:LB}
			&\sum_{|\vec{d}|=d}^{}\int_{X_{\vec{d}}}^{}R(Y_1,\dots,Y_r)\prod_{i=1}^{r}\frac{\theta_i^{p_i}}{p_i!}e^{\theta_iz_i}h_i^{d_i-\bar{g}}
			\\&=N^{-r}\frac{R(w_1,\dots, w_r)}{(w_1\cdots w_r)^{\bar{g}}}\prod_{i=1}^{r}\binom{g}{p_i}w_i^{p_i}[q^{d}](a+b+aq)^{rg-p}(1+q)^{d-r{g}}q^p\nonumber.
		\end{align}
	\end{prop}
	\begin{proof}
		The expression inside the integral is considered in the power series ring $\mathbb{Q}[[x_1,\dots,x_r,\theta_1,\dots,\theta_r]]$. We will first single out the terms containing $\theta_i$. We know that $\theta^k=0$ for $k>g$ thus 
		\begin{align*}
			\frac{\theta_i^{p_i}}{p_i!}e^{\theta_iz_i}&=\sum_{\ell=0}^{g-p_i} \frac{\theta_i^{p_i+\ell}}{p_i!\ell!}\bigg(\frac{B(Y_i)}{P(Y_i)}-\frac{1}{x_i}\bigg)^\ell
		\end{align*}
		We replace $\theta_i^{p_i+\ell}$ by $\frac{g!}{(g-p_i-\ell)!}x_i^{p_i+l}$ using \eqref{eq:theta_int}. We further simplify 
		\begin{align*}
			\sum_{\ell=0}^{g-p_i} \frac{g!x_i^{p_i+\ell}}{p_i!(g-p_i-\ell)!}\frac{1}{\ell!}\bigg(\frac{B(Y_i)}{P(Y_i)}-\frac{1}{x_i}\bigg)^\ell&=\binom{g}{p_i}\cdot x_i^{p_i}\cdot\bigg(\frac{x_iB(Y_i)}{P(Y_i)} \bigg)^{g-p_i}.
		\end{align*}
		
		Plugging this back in \eqref{Eq:LB}, we obtain the following integral of a power series in the variables $x_1,\dots,x_r$
		\begin{align*}
			\sum_{|\vec{d}|=d}^{}\int_{X_{\vec{d}}}^{}R(Y_1,\dots,Y_r)\prod_{i=1}^{r}\binom{g}{p_i}\cdot x_i^{p_i}\cdot\bigg(\frac{x_iB(Y_i)}{P(Y_i)} \bigg)^{g-p_i}h_i^{d_i-\bar{g}}.
		\end{align*}
		 We now have to find the coefficient of $x_1^{d_1}\dots x_r^{d_r}$ in the above expression and sum it over $|\vec{d}|=d_1+\cdots +d_r=d$. For such problems, we have a very useful result from combinatorics, the Lagrange-B\"urmann formula \cite{Lagrange_Burman}, which states
		\begin{equation}\label{eq:L_B}
			\sum_{|\vec{d}|}q_1^{d_1}\cdots q_2^{d_2}([x_1^{d_1}\cdots x_r^{d_r}]f(x_1,\dots,x_r)\prod_{i=1}^{r}h_i^{d_i})=f(x_1,\dots,x_r)\cdot \prod_{i=1}^{r}\frac{1}{h_i}\frac{d{x_i}}{d{q_i}}
		\end{equation}
		where $q_i=\frac{x_i}{h_i}$ and $h_i:=h_i(x_i)$ are power series with $h_i(0)\ne 0$. 
		
		We can apply this formula to
		\begin{align*}
			h_i&=\frac{x_i}{P(Y_i)}   \\
			f(x_1,\dots x_r)&
			= R(Y_1,\dots,Y_r)\prod_{i=1}^{r}\binom{g}{p_i}\cdot x_i^{g}\cdot\bigg(\frac{B(Y_i)}{P(Y_i)} \bigg)^{g-p_i}\bigg(\frac{x_i}{P(Y_i)}\bigg)^{-\bar{g}}\\
			&= R(Y_1,\dots,Y_r)\prod_{i=1}^{r}\binom{g}{p_i}B(Y_i)^{g-p_i}P(Y_i)^{p_i}h_i.
		\end{align*}
		
		We have the change of variable \[q_i=\frac{x_i}{h_i}=P(Y_i)=Y_i^N-1=(x_i+w_i)^N-1,\]and the inverse is given by
		\begin{align*}
			x_i=Y_i-w_i=w_i(1+q_i)^{1/N}-w_i.
		\end{align*}
		Observe that the derivative
		\begin{align*}
			\frac{dx_i}{dq_i}&=\frac{1}{P'(Y_i)}.
		\end{align*}
		
		By direct computation
		\begin{align}\label{eq:power_series_version_theta}
			f(x_1,\dots,x_r)\cdot \prod_{i=1}^{r}\frac{1}{h_i}\frac{d{x_i}}{d{q_i}}&=R(Y_1,\dots,Y_r)\prod_{i=1}^{r}\binom{g}{p_i} \frac{B(Y_i)^{g-p_i}P(Y_i)^{p_i}}{P'(Y_i)}.
		\end{align}
		
		In \eqref{Eq:LB}, we are interested in finding the sum over the coefficients of $q_1^{d_1}\cdots q_r^{d_r}$ where $d_1+\cdots +d_r=d$. To find this sum, we will substitute \[q_1=\cdots=q_r=q\] to obtain a power series in one variable $q$ and find the coefficient of $q^d$. 
		
		In this situation, 
		\begin{align*}
			&Y_i=w_i(1+q)^{1/N}, \hspace{.5cm}B(Y_i)=\frac{(aq+(a+b))}{w_i(1+q)^{1/N}},\\&P'(Y_i)=Nw_i^{-1}(1+q)^{\frac{N-1}{N}}.
		\end{align*}
		Note that $R$ is a homogeneous rational function of degree $s$, thus $R(Y_1,\dots Y_r)=R(w_1,\dots,w_r)(1+q)^{s/N}$. Substituting, the power series \eqref{eq:power_series_version_theta} becomes
		\begin{align*}
			&R(w_1,\dots,w_r)(1+q)^{\frac{s}{N}}\prod_{i=1}^{r}\binom{g}{p_i}\frac{w_i^{p_i-\bar{g}}}{N}\frac{(a+b+aq)^{g-p_i}}{(1+q)^{\frac{g-p_i}{N}+\frac{N-1}{N}}}q^{p_i}\\
			&=(a+b+aq)^{rg-p}(1+q)^{d-r{g}}q^pN^{-r}\frac{R(w_1,\dots, w_r)}{(w_1\cdots w_r)^{\bar{g}}}\prod_{i=1}^{r}\binom{g}{p_i}w_i^{p_i},
		\end{align*}
		where $p=p_1+\cdots+p_r$. 
	\end{proof}
\begin{rem}
	When $p\ge rg$ then $p_i>g$ for some $i$, thus the integral is $0$ since $\theta^p_i=0$. Therefore we may assume that the first term is a polynomial. Moreover, when $d\ge rg$ or $b=0$ and $d\ge p$ then the answer in \eqref{Eq:LB} is given by 
	\begin{align*}
		\frac{a^{rg}}{N^{r}}\frac{R(w_1,\dots,w_r)}{(w_1\cdots w_r)^{\bar{g}}}\prod_{i=1}^{r}\binom{g}{p_i}\frac{w_i^{p_i}}{a^{p_i}}.
	\end{align*}
\end{rem}
	\begin{rem}
		The above proposition, specialized to $B(Y)=P'(Y)$ and $p=0$, greatly simplifies the combinatorics used in finding the Vafa-Intriligator formula for Quot schemes in Section 4 of \cite{marian2007}. 
	\end{rem}
	The previous result does not suffice for the calculation of virtual intersection numbers over isotropic Quot schemes. When rank $r=2$, the following proposition can be used to find Vafa-Intriligator type formulas for $\IQ_d$.
	
	\begin{prop}\label{prop:fixed_loci_calculation}
		Let $R(Y_1,Y_2)$ be a homogeneous rational function of degree $s=Nd-2\bar{g}(N-1)$. We borrow the notation $X_{\vec{d}}$, $Y_i$, $P(Y)$, $B(Y)$, $h_i$ and $z_i$ from Proposition \ref{prop:theta_sum_int}. Let $T(q)=(a+b+aq)/q$. Then we have the following identity
		\begin{align*}
			\sum_{|\vec{d}|=d}^{}\int_{X_{\vec{d}}}^{}&R(Y_1,Y_2)e^{-\frac{\theta_1+\theta_2-\phi_{12}}{Y_1+Y_2}}\prod_{i=1}^{2}e^{\theta_iz_i}h_i^{d_i-\bar{g}}
			\\&=\frac{1}{N^2}\frac{R(w_1,w_2)}{(w_1w_2)^{\bar{g}}}[q^d](1+q)^{d}\bigg(\frac{qT(q)}{1+q}\bigg)^{2g}\bigg(1-\frac{1}{T(q)}\bigg)^g.
		\end{align*}
		In particular, when $d\ge2g$ the above value is 
		\begin{equation*}
			\frac{a^g(a-1)^g}{N^2}\frac{R(w_1,w_2)}{(w_1w_2)^{\bar{g}}}.
		\end{equation*}
	\end{prop}
	\begin{proof}
		We will first replace exponents of $\phi_{12}$ with the exponents of $\theta_1\theta_2$ using Proposition \ref{prop:phi_int}. For parity reasons $\phi_{12}$ must appear with an even power to obtain a non-zero number. Thus we can make following replacements:
		\begin{align*}
			e^{-\frac{\theta_1+\theta_2- \phi_{12}}{Y_1+Y_2}}
			&\to\sum_{p=0}^{\infty}\frac{(-1)^p}{p!(Y_1+Y_2)^p}\bigg(\sum_{2\ell+r+s=p}^{}\binom{p}{2\ell,r,s}\theta_1^{r}\theta_2^{s}\phi_{12}^{2\ell}  \bigg)\\
			&\to\sum_{p=0}^{\infty}\sum_{2\ell+r+s=p}^{}\frac{(-1)^{p-\ell}}{p!}\binom{p}{2\ell,r,s}\frac{\binom{2\ell}{\ell}}{\binom{g}{\ell}}\frac{\theta_1^{r+\ell}\theta_2^{s+\ell}}{(Y_1+Y_2)^p} \\
			&= \sum_{p=0}^{\infty}\sum_{2\ell+r+s=p}^{}\frac{(-1)^{p-\ell}}{(Y_1+Y_2)^p}\frac{\binom{p}{2\ell,r,s}}{\binom{p}{r+\ell}}\frac{\binom{2\ell}{\ell}}{\binom{g}{\ell}}\frac{\theta_1^{r+\ell}\theta_2^{s+\ell}}{(r+\ell)!(s+\ell)!}.
		\end{align*}
		Now we use Proposition \ref{prop:theta_sum_int} to reduce the problem to finding 
		\begin{align*}
			\sum_{\substack{2\ell+r+s=p}}^{}(-1)^{p-\ell} \frac{\binom{p}{2\ell,r,s}}{\binom{p}{r+\ell}}\frac{\binom{2\ell}{\ell}}{\binom{g}{\ell}}\cdot \frac{1}{N^2}\frac{R(w_1,w_2)w_1^{r+\ell}w_2^{s+\ell}}{(w_1+w_2)^p(w_1w_2)^{\bar{g}}}\binom{g}{r+\ell}\binom{g}{s+\ell}
			\\
			\cdot[q^{d}](1+q)^d\bigg(\frac{a+b+aq}{1+q}\bigg)^{2g}\bigg(\frac{q}{a+b+aq}\bigg)^{p}
		\end{align*}
		where the sum is taken over $r,s,\ell$ such that $r+\ell,s+\ell\le g$. Rearranging the binomial coefficients, the above expression is same as
		\begin{align*}
			[q^d](1+q)^d&\bigg(\frac{a+b+aq}{1+q}\bigg)^{2g}\frac{1}{N^2}\frac{R(w_1,w_2)}{(w_1w_2)^{\bar{g}}}\\
			&\cdot\sum_{\substack{2\ell+r+s=p}}^{} (-1)^{\ell}\binom{g}{\ell}\binom{g-\ell}{r}\binom{g-\ell}{s}\frac{(-w_1)^{r+\ell}(-w_2)^{s+\ell}}{T(q)^p(w_1+w_2)^p}.
		\end{align*}
		The summation in the above expression greatly simplifies via the following lemma.
	\end{proof}

	\begin{lem}
		Let $g$ and $d$ be integers, then 
		\begin{align*}
			\sum_{\substack{2\ell+r+s=p}}^{}(-1)^{\ell}\binom{g}{\ell}\binom{g-\ell}{r}\binom{g-\ell}{s}\frac{(-w_1)^{r+\ell}(-w_2)^{s+\ell}}{T(q)^p(w_1+w_2)^p} =\bigg(1-\frac{1}{T(q)}\bigg)^g.
		\end{align*}
	\end{lem}
	\begin{proof}
		The lemma follows by observing that the given expression simplifies as 
		\begin{align*}
			&\sum_{\ell}^{}\binom{g}{\ell}\frac{(-1)^{\ell}}{T(q)^{2\ell}}\frac{(-w_1)^{\ell}(-w_2)^{\ell}}{(w_1+w_2)^{2\ell}} \bigg(1-\frac{w_1}{T(q)(w_1+w_1)}\bigg)^{g-\ell}
			\bigg(1-\frac{w_2}{T(q)(w_1+w_1)}\bigg)^{g-\ell}
			\\
			&=\bigg(\bigg(1-\frac{w_1}{T(q)(w_1+w_1)}\bigg)
			\bigg(1-\frac{w_2}{T(q)(w_1+w_1)}\bigg)-\frac{w_1w_2}{T(q)^2(w_1+w_2)^2}\bigg)^g\\
			&=\bigg(1-\frac{1}{T(q)}\bigg)^g.
		\end{align*}
	\end{proof}

	\section{Intersection of $a$-classes}\label{sec:Int_a_classes}
	
	In this section we will prove Theorem \ref{thm:r=2,sympl} and \ref{thm:r=2_symmetric}, which are explicit expressions for the intersections of $a$-classes in the symplectic and symmetric case respectively. 
	\subsection{$a$-class intersections for $\sigma$ symplectic} Let $r=2$. In this case the virtual dimension of $\IQ_d$ is given by
	\[\vd = (N-1)d-(2N-5)\bar{g}. \]
	
	Let us define  
\begin{equation}\label{eq:T_d,g}
T_{d,g}(N)= [q^d](1+q)^{d-g}\bigg(1+\frac{N-1}{N}q\bigg)^g.
\end{equation}In particular, when $d\ge g$, we get $T_{d,g}(N)=(1-1/N)^g$. A simple usage of Lagrange inversion theorem implies \[T_{d,g}(N)=[q^d](1-q/N)^g(1-q)^{-1}\] and hence $T_{d,g}(N)$ is the sum of the first $d$ terms in the binomial expansion of $(1-1/N)^g$.
	
	\begin{theorem}\label{thm:7.1_r=2_sympl}
		Let $Q(X_1,X_2)$ be a polynomial of weighted degree $\vd$, where the variables $X_i$ have degree $i$. 
		Then, 
		\begin{equation}
			\int_{[\IQ_d]^{\vir}}^{}Q(a_1,a_2) =uT_{d,g}(N) \sum_{w_1,w_2}^{}S(w_1,w_2)J(w_1,w_2)^{\bar{g}}(w_1+w_2)^d
		\end{equation}
		where the sum is taken over all the pairs of $N^\text{th}$ roots of unity $\{w_1,w_2\}$ with $w_1\ne \pm w_2$. Here $u=(-1)^{\bar{g}+d}$ and
		\begin{align*}
			J(w_1,w_2)&=N^2w_1^{-1}w_2^{-1}(w_1-w_2)^{-2}(w_1+w_2)^{-1},
		\end{align*}  and $S(w_1,w_2)=Q(w_1+w_2,w_1w_2)$.
	\end{theorem}
	\begin{proof}
		The equivariant pull back of $a_i$ to the fixed loci is the $i$th elementary symmetric function $\sigma_i((w_1t+x_1),(w_2t+x_2))$, hence $Q(a_1,a_2)$ pulls back to $S(w_1t+x_1,w_2t+x_2)$. We are in a position to apply the equivariant virtual localization formula \cite{Graber1997LocalizationOV} which yields 
		\begin{align}\label{eq:localiz_r=2_symp}
			\int_{[\IQ_d]^{\vir}}^{}Q(a_1,a_2)=\sum_{d_1+d_2=d}^{}\sum_{w_1,w_2}^{} \int_{\fix_{\vec{d},\underline{k}}}^{} \frac{S(Y_1,Y_2)}{e_{\mathbb{C}^*}(\nbun)},
		\end{align}
		where the sum is taken over all the prescribed choices for $\{w_1,w_2\}$ and $Y_i=x_i+w_it$. 
		
		After appropriately replacing $\theta$ and $\phi_{12}$ classes with $x$ classes as described in Section \ref{sec:prelim_hilb}, the above expression can be written as a rational function in $x_1,x_2$ and $t$ of with total degree $d$ . The integral can thus be evaluated by finding coefficient of $x_1^{d_1}x_2^{d_2}$. The homogeneity and the identity $d_1+d_2=d $ ensures that resulting element in $\mathbb{C}[t,t^{-1}]$ has $t$ degree 0. Hence we can safely assume $t=1$ for the purpose of our calculation without changing the value of integral.
		
		Moreover, the localization formula is independent of the choice of the weights $(w_1,\dots w_N)$ as long as these are distinct and satisfy $w_i=-w_{i+n}$ for $1\le i\le n$. Hence we may assume these to be distinct roots of the polynomial $P(X)=X^N-1$.
		
		We substitute the expression \eqref{eq:e(Nvir)spl} of the Euler class of $\nbun$ into \eqref{eq:localiz_r=2_symp} to get 
		\begin{align*}
			\sum_{w_1,w_2}^{}\sum_{d_1+d_2=d}^{} \int_{\fix_{\vec{d},\underline{k}}}^{} R(Y_1,Y_2)e^{-\frac{\theta_1+\theta_2-\phi_{12}}{Y_1+Y_2}}\prod_{i=1}^{2}e^{\theta_iz_i}h_i^{d_i-\bar{g}},
		\end{align*}
		where by \eqref{eq:def_z_i} $z_i=\frac{P'(Y_i)}{P(Y_i)}-\frac{1}{x_i}$, $h_i=\frac{x_i}{P(Y_i)}$ and \[R(Y_1,Y_2)=uS(Y_1,Y_2)\frac{(Y_1+Y_2)^{d-\bar{g}}}{(Y_1-Y_2)^{2\bar{g}}}. \]
		The homogeneous degree of $R$ is $\vd +(d-3\bar{g})=Nd -2\bar{g}(N-1)$,  therefore Proposition \ref{prop:fixed_loci_calculation} gives the required intersection number 
		\begin{align}
			\sum_{w_1,w_2}\frac{1}{N^2}\frac{R(w_1,w_2)}{(w_1w_2)^{\bar{g}}}[q^d]N^{2g}(1+q)^{d-g}\bigg(1+\frac{N-1}{N}q\bigg)^g,
		\end{align}
	completing the proof.
	\end{proof}
	\begin{proof}[Proof of Theorem \ref{thm:r=2,sympl}]	
		In the statement of Theorem  \ref{thm:7.1_r=2_sympl}, the expression \[S(w_1,w_2)J(w_1,w_2)^{\bar{g}}(w_1+w_2)^d\] is homogeneous of degree $N(d-2\bar{g})$, hence this equals
		$S(1,\zeta)J(1,\zeta)^{\bar{g}}(1+\zeta)^d,$
		where $\zeta=w_2/w_1$.
	\end{proof}

	\subsection{$a$-class intersections for $\sigma$ symmetric}\label{sec_r=1_symmetric}
	Define \[\tilde{T}_{d,g}(N)=[q^d]\bigg(1+\frac{N-2}{N}q\bigg)^g(1+q)^{d-g}. \]
	\begin{prop}\label{prop:r=1}
		Over $\IQ_d$, where $N$ is even, $r=1$ and $\sigma$ is symmetric, the top intersection of the tautological class is given by
		\begin{align}
			\int_{[\IQ_d]^{\vir}}^{}a_1^{\vd}= N^{g}\tilde{T}_{d,g}(N)2^{2d-\bar{g}}
		\end{align}
		where $\vd=(N-2)(d-\bar{g})$ is the virtual dimension.
	\end{prop}
	\begin{proof}
		The restriction of $a_1$ to the fixed locus $\fix_{d, i}=C^{[d]}$ is $Y_i=x_i+w_it$.
		The Euler class of the equivariant normal bundle of the fixed locus is given by \eqref{eq:e(N)_symmetric}
		\begin{align*}
			\frac{1}{e^{\vir}_{\mathbb{C}^*}(\nbun)}&=2^{2d-\bar{g}}Y_i^{2d-\bar{g}}h_i^{d-\bar{g}}e^{\theta_iz_i }
		\end{align*}
		where  $z_i=(B(Y_i)/P(Y_i)-1/x_i)$ and \[\frac{B(Y)}{P(Y)}= \frac{P'(Y)}{P(Y)}-\frac{2}{Y}. \] The equivariant virtual localization formula gives 
		\begin{equation*}
			\int_{[\IQ_d]^{\vir}}^{}a_1^{\vd}=\sum_{i=1}^{N}\int_{\fix_{d, i}}\frac{Y_i^{\vd}}{e^{\vir}_{\mathbb{C}^*}(\nbun)}.
		\end{equation*}
		We choose the weight of the action to be $N^{\text{th}}$ roots of unity, thus $P(X)=X^N-1$, hence $B(Y)= \frac{(N-2)Y^N+2}{Y}$, and we obtain the integral as a special case of Proposition \ref{prop:theta_sum_int} by putting $r=1$ and $p=0$. 
	\end{proof}
	\begin{rem}
		Similar results can be obtained when $N$ is odd, $r=1$ and $\sigma$ symmetric. In particular, when the virtual dimension is non-zero, 
				\begin{align}
					\int_{[\IQ_d]^{\vir}}^{}a_1^{\vd}= (N-1)^{g}2^{2d-\bar{g}}T_{d,g}(N-1).
				\end{align}
	\end{rem}

	When $r=2$,  localizing with distinct weights makes combinatorics very difficult. However using two equal weights enable us to find a simple formula for these intersections. Using exactly two equal weights results in getting $C^{[d_1]}\times\IQ_{d_2}(\mathbb{C}^2\otimes\cO, r=1,\sigma)$ as part of the fixed loci. We will first show that \[\IQ_{d}(\mathbb{C}^2\otimes\cO, r=1,\sigma)= C^{[d]}\sqcup C^{[d]},\]and the two components $C^{[d]}$ come equipped with a non-standard virtual structure. We will use Proposition \ref{prop:r=1} to understand how to intersect over these non-standard loci.
	
	Recall that the virtual dimension of $\IQ_d$ is 
	\begin{equation*}
		\vd=(N-3)d -\bar{g}(2N-7).
	\end{equation*}

	Let $N=2n$. Let $G=\mathbb{C}^*$ act on $\IQ_d$ with weights \[(w_1,\dots,w_N)=(\zeta ,\zeta^2,\dots\zeta^{n-1},0,\zeta^n,\dots,\zeta^{2n-2},0),\] where $\zeta$ is a primitive $(N-2)$'th root of unity. A point $[0\to S\to \mathbb{C}^N\otimes\cO\to Q\to 0]$ in $\IQ_d$ is fixed under the action of $G$ if and only if one of the following is satisfied:
	\begin{itemize}
		\item[(i)] The sheaf $S$ splits as $L_1\oplus L_2$ where $L_i$ is a subsheaf of one of the $N-2$ copies of $\cO$, at position $k_i\notin\{n,2n\}$, in $\mathbb{C}^N\otimes\cO$ such that $k_1-k_2\not \equiv 0 \mod n$. The corresponding fixed locus is \[\fix_{\vec{d},\underline{k}} \cong C^{[d_1]}\times C^{[d_2]}, \] where $\deg{L_i}=d_i$ and $\underline{k}=(k_1,k_2)$.
		\item[(ii)] The sheaf $S$ splits as $L_1\oplus E$ where $L_1$ is a subsheaf of one the copies of $\cO$, at position $k\notin \{n,2n\}$, in $\mathbb{C}^N\otimes\cO$ and $E$ is an isotropic rank one subsheaf of $\cO_n\oplus\cO_{2n}$, the sum of copies of $\cO$ at positions $n$ and $2n$. Let $\fix_{\vec{d}, k}$ be the component of the fixed loci consisting of $(L_1,E)$, where $d_1=\deg L_1$, $d_2=\deg E$ and $k$ is the position mentioned above. Note that \[\fix_{\vec{d}, k}\cong C^{[d_1]}\times\IQ_{d_2}(\cO\otimes\mathbb{C}^2, r=1,\sigma). \]
	\end{itemize}
	
	\begin{theorem}
		Let $Q(X_1,X_2)$ be a polynomial of weighted degree $\vd$, where the variables $X_i$ have degree $i$. 
		Then, 
		\begin{equation*}
			\int_{[\IQ_d]^{\vir}}^{}Q(a_1,a_2) = I_1+I_2
		\end{equation*}
		where $S(X_1,X_2)=Q(X_1+X_2,X_1X_2)$,
		\begin{align*}
			I_1&=u4^dT_{d,g}(N-2)\sum_{w_1\ne\pm w_2}S(w_1,w_2)J(w_1,w_2)^{\bar{g}}(w_1+w_2)^d,\\
			I_2&=(-1)^{d}2^{2d+2-g}T_{d,g}(N-2)(N-2)^{g}\cdot Q(1,0),
		\end{align*}
		and $	J(w_1,w_2)=\frac{(N-2)^2}{4}(w_1+w_2)^{-1}(w_1-w_2)^{-2}.$ 
		
	\end{theorem}
	\begin{proof}
		Using equivariant virtual localization formula, we can write 
		\begin{align*}
			\int_{[\IQ_d]^{\vir}}Q(a_1,a_2)= I_1+I_2,
		\end{align*}
		where 
		\begin{align*}
			I_1&= \sum_{\substack{k_1,k_2\notin\{n,2n\} \\|k_1-k_2|\ne n }}^{}\sum_{d_1+d_2=d}^{}\int_{\fix_{\vec{d},\underline{k}}}^{}\frac{i^*(Q(a_1,a_2))}{e_{\mathbb{C}^*}(\nbun_{\fix_{\vec{d},\underline{k}}})}\\
			I_2&= \sum_{\substack{k\in[N]\\ k\notin\{n,2n\}}}^{}\sum_{d_1+d_2=d}^{}\int_{\fix_{\vec{d},k}}^{}\frac{i^*(Q(a_1,a_2))}{e_{\mathbb{C}^*}(\nbun_{\fix_{\vec{d},k}})}.
		\end{align*}
		Here we denote $i^*$ the restriction to the fixed loci. The next two subsections will be devoted to the calculation of $I_1$ and $I_2$ respectively.
	\end{proof}

	\subsubsection{Fixed loci of the first kind}
	$\fix_{\vec{d},\underline{k}}=C^{[d_1]}\times C^{[d_2]}$ . In Section \ref{sec:Euler_class_symm} we noted that the $\mathbb{C}^{*}$ equivariant virtual tangent bundle is given by 
	\begin{equation*}
		T^{\vir}= \pi_![(\rHom(\cS,\cQ))]-\pi_![(\sHom(\Sym^2{\cS},\cO))].
	\end{equation*}
	The non-moving part of the restriction of $T^{\vir}$ to $\fix_{\vec{d},\underline{k}}$ matches the $K$-theory class the tangent bundle of $\fix_{\vec{d},\underline{k}}$.
	The virtual normal bundle
	\begin{equation*}
		\nbun=\pi_*\bigg(\sum_{\substack{i=1,2\\1\le k\le N\\ k_i\ne k}}[\mathcal{L}^\vee_i]-\sum_{\substack{i,j\in [2]\\i\ne j}}^{} [\mathcal{L}_i^\vee\otimes\mathcal{L}_j ]-\sum_{1\le i\le j\le 2}^{} [\mathcal{L}^\vee_i\otimes\mathcal{L}^\vee_j]\bigg).
	\end{equation*}
	Therefore using \eqref{eq:e(N)_symmetric}, we have
	\begin{align}\label{eq:Nvir_symm_r=2}
		\frac{1}{e_{\mathbb{C}^*}(\nbun)}&= u2^{2d-2\bar{g}}\frac{(Y_1+Y_2)^{d-\bar{g}}}{(Y_1-Y_2)^{2\bar{g}}}(Y_1Y_2)^{\bar{g}}e^{-\frac{\theta_1+\theta_2-\phi_{12}}{Y_1+Y_2}}\prod_{i=1}^{2}h_i^{d_i-\bar{g}}e^{\theta_iz_i}
	\end{align}
	where $P_0(X)=X^{N-2}-1$ and 
	\begin{align*}
		h_i&=\frac{x_iY_i^2}{P(Y_i)}=\frac{x_i}{P_0(Y_i)}, \hspace{1cm}
		B(Y_i)=P_0'(Y_i), \\
		z_i&=\frac{P'(Y_i)}{P(Y_i)}-\frac{2}{Y_1}-\frac{1}{x_i} 
		= \frac{B(Y_i)}{P_0(Y_i)}-\frac{1}{x_i}.
	\end{align*}
	\begin{prop}
		We have 
		\begin{align}
			I_1= u4^dT_{d,g}(N-2)\sum_{w_1,w_2}S(w_1,w_2)J(w_1,w_2)^{\bar{g}}(w_1+w_2)^d
		\end{align}
		where the sum is taken over pairs of $(N-2)^{\text{th}}$ roots of unity $\{w_1,w_2\}$ with $w_1\ne \pm w_2$, and
		\begin{align*}
			J(w_1,w_2)=\frac{(N-2)^2}{4}(w_1+w_2)^{-1}(w_1-w_2)^{-2}.
		\end{align*}
		In particular when $d\ge g$, $T_{d,g}(N-2)=(N-3)^g(N-2)^{-g}$.
	\end{prop}
	\begin{proof}
		For notational convenience, we assume $\underline{k}=(1,2)$.  The classes $a_1$ and $a_2$ restrict to $Y_1+Y_2$ and $Y_1Y_2$ respectively, where $Y_i= x_i+w_it$ in the equivariant cohomology ring $H^*(\fix_{\vec{d},\underline{k}}=C^{[d_1]}\times C^{[d_2]})[[t]]$. 
		
		We are interested in evaluating the following sum
		\begin{align*}
			\sum_{d_1+d_2=d}^{}\sum_{w_1,w_2}^{} \int_{\fix_{\vec{d},\underline{k}}}^{} \frac{S(Y_1,Y_2)}{e_{\mathbb{C}^*}(\nbun_{\fix_{\vec{d},\underline{k}}})},
		\end{align*}
		where $S(Y_i,Y_i)=Q(Y_1+Y_2,Y_1Y_2)$. After replacing the classes $\theta_i$ and $\phi_{12}$ as in the proof of Theorem \ref{thm:r=2,sympl}, the above expression becomes a homogeneous degree rational function of degree $d=d_1+d_2$ in the variables $x_i$ and $t$ and a power series in $x_1$ and $x_2$ with coefficients in $\mathbb{C}[[t,t^{-1}]]$. Integrating over $C^{[d_1]}\times C^{[d_2]}$ amounts to finding the coefficient of $x_1^{d_1}x_2^{d_2}$.
		
		Using the calculation of $e(\nbun)$ in \eqref{eq:Nvir_symm_r=2}, we reduce our problem to finding 
		\begin{equation*}
			\sum_{d_1+d_2=d}^{}\sum_{w_1,w_2}^{} \int_{\fix_{\vec{d},\underline{k}}}^{} R(Y_1,Y_2)e^{-\frac{\theta_1+\theta_2-\phi_{12}}{Y_1+Y_2}}\prod_{i=1}^{2}h_i^{d_i-\bar{g}}e^{\theta_iz_i}
		\end{equation*}
		where $(w_1,w_2)$ are the prescribed pair of $(N-2)$'th roots of unity and \[R(Y_1,Y_2)=u2^{2d-2\bar{g}}S(Y_1,Y_2)(Y_1Y_2)^{\bar{g}}\frac{(Y_1+Y_2)^{d-\bar{g}}}{(Y_1-Y_2)^{2\bar{g}}}. \]
		We apply Proposition \ref{prop:fixed_loci_calculation} to find
		\begin{align*}
			I_1=\sum_{w_1,w_2}^{}\frac{1}{(N-2)^2}\frac{R(w_1,w_2)}{(w_1,w_2)^{\bar{g}}}[q^d](N-2)^{2g}(1+q)^{d-g}\bigg(1+\frac{N-3}{N-2}q\bigg)^g.
		\end{align*}
	\end{proof}

	\subsubsection{Fixed Loci of second kind} We will first understand the virtual geometry of the isotropic Quot scheme $\IQ_d^\circ=\IQ_{d}(\cO\otimes\mathbb{C}^2, r=1,\sigma)$.
	
	\begin{lem}
		The isotropic Quot scheme $\IQ_d^\circ$ is isomorphic to the disjoint union $C^{[d]}\sqcup C^{[d]}$. The virtual tangent bundle of $\IQ_d^\circ$ restricted to either copy of $C^{[d]}$ is given by
		\begin{align*}
			T^{\vir}=\pi_!([\mathcal{L}^\vee\otimes (\mathcal{T}\oplus \cO)]-[\mathcal{L}^\vee\otimes\mathcal{L}^\vee]),
		\end{align*}
		where $\pi$ is the projection $\pi : C\times C^{[d]}\to C^{[d]}$ and $0\to\mathcal{L}\to \cO\to \mathcal{T}\to 0$ is the universal exact sequence on $C\times C^{[d]}$.
	\end{lem}
	\begin{proof}
		A subsheaf $E\subset \mathbb{C}^2\otimes\cO$ is isotropic if and only if $E$ factors through a copy of $\cO$ in $\mathbb{C}^2\otimes\cO$, hence $\IQ_d^\circ\cong C^{[d]}\sqcup C^{[d]}$. The universal short exact sequence over $C\times \IQ_d^\circ$ restricts to \[0\to \mathcal{L}\to \mathbb{C}^2\otimes\cO \to \mathcal{T}\oplus \cO\to 0  \] over each copy of $C\times C^{[d]}$. The lemma follows using the description of $T^{\vir}$ of $\IQ_d^\circ$ in Theorem \ref{thm:POT}.
	\end{proof}
	Therefore we see that the virtual fundamental class $[C^{[d]}]^{\vir}$ induced over each  component $C^{[d]}$ of $\IQ_{d}^\circ$ is different from the usual fundamental class $[C^{[d]}]$.  We also observe that the virtual dimension for $C^{[d]}$ is zero. 
	\begin{lem}\label{lem:Cvir_integral}
		Let $C^{[d]}$ be equipped with the non-standard virtual structure as described above, then 
		\begin{equation*}
			\int_{[C^{[d]} ]^{\vir}}1=2^{2d}(-1)^d\binom{\bar{g}}{d}. 
		\end{equation*}
	\end{lem}
	\begin{proof}
		We have a natural automorphism obtained by swapping the copies of the $\cO$ in $\mathbb{C}^{2}\otimes\cO$. Therefore the above intersection number is independent of the copy of $C^{[d]}$ we have chosen. The Proposition \ref{prop:r=1} tells us
		\begin{align*}
			\int_{[C^{[d]} ]^{\vir}}1=\frac{1}{2}\int_{[\IQ_{d}^\circ]^{\vir}}^{}1=2^{2d}[q^d](1+q)^{d-g}.
		\end{align*}
	\end{proof}
	Now we are ready to prove
	\begin{prop}
		We have 
		\begin{align*}
			I_2=(-1)^{d}2^{2d+2-g}(N-2)^{g}T_{d,g}(N-2)\cdot Q(1,0)
		\end{align*}
	\end{prop}
	\begin{proof}
		We are working over the fixed loci $\fix_{\vec{d}, k,\epsilon}=C^{[d_1]}\times \Cvir^{[d_2]}$ where $ k\notin\{n,2n\}$ and the first factor corresponds to the copy of $\cO$ at position $k$ and the index $\epsilon$ differentiates between the two components of $\IQ_{d_2}^{0}=C^{[d_2]}\sqcup C^{[d_2]}$. Let $\mathcal{L}_1$ and $\mathcal{L}_2$ be the pullbacks of the universal subsheaves over $C^{[d_1]}$ and $\Cvir^{[d_2]}$ to the product $\fix_{\vec{d}, k,\epsilon}$. The virtual normal bundle is the moving part of the restriction of the $T^{\vir}$ and is given by
		\begin{equation*}
			\nbun=\pi_!\bigg(\sum_{j\in [N]-\{k\}}[\mathcal{L}_1^\vee] + \sum_{\substack{j\in [N]\\ j\notin\{n,2n\}}}^{}[\mathcal{L}_2^{\vee}] -
			[\mathcal{L}_1^\vee\otimes\mathcal{L}_2]-[\mathcal{L}_1\otimes\mathcal{L}_2^\vee]- [\mathcal{L}^\vee_1\otimes\mathcal{L}^\vee_2] -[\mathcal{L}^\vee_1\otimes\mathcal{L}^\vee_1]\bigg),
		\end{equation*}
		where the above terms have $\mathbb{C}^*$ weights $(w_k-w_j)$, $-w_j$, $w_k$, $-w_k$, $w_k$ and $2w_k$ respectively.
		
		We may assume $t=1$ (see the proof of Theorem \ref{thm:7.1_r=2_sympl}. Let $Y_1=x_1+w_k$, $u=(-1)^{d+\bar{g}}$ and $P(X)=X^{N-2}-1$. A careful calculation using \eqref{eq:euler_class_computation_1}, \eqref{eq:euler_class_computation_2} and \eqref{eq:euler_class_computation_3} gives 
		\begin{align*}
			\frac{1}{e_{\mathbb{C}^*}(\nbun)}=& \bigg(\frac{Y_1^2P(Y_1)}{x_1}\bigg)^{-d_1+\bar{g}} e^{\theta_1\big(\frac{P'(Y_1)}{P(Y_1)}+\frac{2}{Y_1}-\frac{1}{x_1}\big)}
			\cdot P(x_\epsilon)^{-d_2+\bar{g}}e^{\theta_\epsilon\frac{P'(x_\epsilon)}{P(x_\epsilon)}}\\
			&\cdot u(Y_1-x_\epsilon)^{-2\bar{g}}\cdot(Y_1+x_\epsilon)^{d-\bar{g}}e^{(-\frac{\theta_1+\theta_\epsilon-\phi_{12}}{(Y_1+x_\epsilon)})}\cdot (2Y_1)^{2d_1-\bar{g}}e^{-\frac{2\theta_1}{Y_1}}
		\end{align*}
		Since $\Cvir^{[d_2]}$ has virtual dimension zero, $x_\epsilon$ and $\theta_\epsilon$ yield zero when intersected with the virtual fundamental class $[\Cvir^{[d_2]}]^{\vir}$. Thus for the purpose of our calculation, we may substitute $x_\epsilon=\theta_\epsilon=\phi_{12}=0$ in the above expression to get
		\begin{align*}\label{eq:nbun_vir_r=2}
			u2^{2d_1-\bar{g}}Y_1^{d-2\bar{g}}h_1^{d_1-\bar{g}}e^{\theta_1 z_1}\cdot (-1)^{(\bar{g}-d_2)},
		\end{align*}
		where $h_1=x_1/P(Y_1)$ and $z_1=P'(Y_1)/P(Y_1)-1/Y_1-1/x_1$.
		
		Note that $a_1$ and $a_2$ restrict to $Y_1+x_\epsilon$ and $Y_1x_\epsilon$ respectively over the fixed loci. We want to calculate
		
		\begin{align*}
			I_2&= \sum_{k=1}^{N-2}\sum_{d_1+d_2=d}^{}\sum_{\epsilon=1}^{2}\int_{[\fix_{\vec{d},k,\epsilon}]^{\vir}}^{}\frac{i^*(Q(a_1,a_2))}{e_{\mathbb{C}^*}(\nbun_{\fix_{\vec{d},k}})}
		\end{align*}
		Substituting $x_\epsilon=0$, we get
		\begin{align}
			I_2&= Q(1,0)\sum_{k=1}^{N-2}\sum_{d_1+d_2=d}^{}\sum_{\epsilon=1}^{2}\int_{[\fix_{\vec{d},k,\epsilon}]^{\vir}}^{}\frac{Y_1^{\vd}}{e_{\mathbb{C}^*}(\nbun)}.
		\end{align}
		Simplifying further using Lemma \ref{lem:Cvir_integral}, we get 
		\begin{align*}
			I_2&=Q(1,0)\sum_{k=1}^{N-2}\sum_{\epsilon=1}^{2}\sum_{d_1+d_2=d}^{}u2^{2d_1-\bar{g}}(-1)^{\bar{g}-d_2} \int_{C^{[d_1]}}^{}Y_1^{\vd +d-2\bar{g}}h_1^{d_1-\bar{g}}e^{\theta_1 z_1} \int_{[C^{[d_2]}]^{\vir}}1  \\
			&=
			Q(1,0)\sum_{k=1}^{N-2}\sum_{\epsilon=1}^{2}\sum_{d_1+d_2=d}^{}u2^{2d-\bar{g}}(-1)^{\bar{g}}\binom{\bar{g}}{d_2}\int_{C^{[d_1]}}^{}Y_1^{\vd +d-2\bar{g}}h_1^{d_1-\bar{g}}e^{\theta_1 z_1}\\
			&= 	Q(1,0)\sum_{k=1}^{N-2}\sum_{\epsilon=1}^{2}u2^{2d-\bar{g}}(-1)^{\bar{g}}(N-2)^{\bar{g}}[q^d](1+q)^{d-g}\bigg(1+\frac{N-3}{N-2}q\bigg)^g.
		\end{align*}
		The last equality follows from noting that $\binom{\bar{g}}{d_2}=[q^{d_2}](1+q)^{\bar{g}}$ and the following Lemma.
	\end{proof}
	\begin{lem}
		\begin{align*}
			\int_{C^{[d_1]}}^{}Y_1^{\vd +d-2\bar{g}}h_1^{d_1-\bar{g}}e^{\theta_1 z_1} = (N-2)^{\bar{g}}[q^{d_1}](1+q)^{d-\bar{g}-g}\bigg(1+\frac{N-3}{N-2}q\bigg)^g
		\end{align*}
	\end{lem}
	\begin{proof}
		Proposition \ref{prop:theta_sum_int} does not directly apply here due to shape of $d_1$. However, we closely follow the proof of Proposition \ref{prop:theta_sum_int}. Correctly replacing $e^{\theta_1 z_1}$ yield
		\begin{align*}
			\int_{C^{[d_1]}}^{}Y_1^{\vd +d-2\bar{g}}h_1^{d_1-\bar{g}}\bigg(\frac{x_1B(Y_1)}{P'(Y_1)}\bigg)^g.
		\end{align*}
		Applying the Lagrange-B\"urmann formula, we obtain
		\begin{align*}
			[q^{d_1}]Y_1^{\vd+d-2\bar{g}}\frac{B(Y_1)^g}{P'(Y_1)}
		\end{align*}
		where $Y_1=w_1(1+q)^{\frac{1}{N-2}}$ and $Y_1B(Y_1)=(N-3)Y_1^{N-2}+1$. Therefore, it equals 
		\begin{align*}
			(N-2)^{\bar{g}}[q^{d_1}](1+q)^{d-\bar{g}-g}\bigg(1+\frac{N-3}{N-2}q\bigg)^g.
		\end{align*}
	\end{proof}

	\section{Gromov-Ruan-Witten Invariants}\label{sec:GRW_Invariants}

	In this section we will compare the sheaf theoretic invariants obtained using isotropic Quot schemes and Gromov-Ruan-Witten invariants for Isotropic Grassmannians. We will denote by $\SG(2,N)$ and $\OG (2,N)$ the symplectic Grassmannian and orthogonal Grassmannian respectively.
	
	\subsection{Quantum Cohomology}
	The small quantum cohomology of the Isotropic Grassmannian and its presentation are known (see \cite{Quantum_pieri_OG}, \cite{QC_isotropic_grassmannians}). However, the explicit expressions for the high genus and large degree Gromov-Ruan-Witten invariants require further arguments. 
	
	When the rank $r=2$, a simpler presentation for the quantum cohomology of $\SG(2,2n)$ was obtained in \cite{Q_C_IG}. We will briefly describe their result and find a similar presentation for the quantum cohomology of $\OG(2,2n+2)$. 
	
	Let $N=2n$. We have the universal exact sequence $0\to \cS\to \mathbb{C}^N\otimes\cO\to \cQ\to 0$ over $\SG(2,N)$. Let $\cS^{\perp}\subset \mathbb{C}^N\otimes\cO$ be the rank $N-2$ vector bundle consisting of vectors perpendicular to $\cS$. 

	Moreover, $\cS^\perp$ is the kernel of the composition $\mathbb{C}^N\otimes\cO\xrightarrow{\sigma}(\mathbb{C}^N)^{\vee}\otimes\cO\to \cS^\vee $ which gives us an identity for the Chern polynomial $c_t(\cS^\vee)c_t(\cS^{\perp})=1$. This implies
	\begin{align}\label{eq:QC_SG_identity}
		c_t(\cS)c_t(\cS^\vee)c_t(\cS^{\perp}/\cS)&=1.
	\end{align}
	The above identity suggests us to define the following cohomology classes :
	\begin{itemize}
		\item The Chern classes $a_i= c_i(\cS^\vee)$ for $i\in\{1,2\}$.
		\item Let $b_i=c_{2i}(\cS^{\perp}/\cS)$ for $i\in \{1,\dots ,n-2 \}$. The bundle $\cS^{\perp}/\cS$ is self dual, hence all the odd Chern classes vanish.
	\end{itemize}
	The cohomology ring $H^*(\SG(2,2n))$ is isomorphic to the quotient of the ring $\mathbb{C}[a_1,a_2,b_1,\dots ,b_{n-2}]$ by the ideal generated by 
	\begin{equation}\label{eq:H*SG(2,2n)}
		(1+(2a_2-a_1^2)x^2+a_2x^4)(1+b_1x^2+\dots + b_{n-2}x^{2n-4})=1.
	\end{equation}
	The above identity is simply a restatement of \eqref{eq:QC_SG_identity}. The quantum cohomology ring is $H^*(\SG(2,2n))\otimes \mathbb{C}[[q]]$, where the quantum products is described in the following theorem. Note that $\deg(q)=2n-1$ is the index of $\SG(2,2n)$.
	\begin{theorem}[\cite{Q_C_IG}]
		The quantum cohomology ring $QH^*(\SG(2,2n))$ is isomorphic to the quotient of the ring $\mathbb{C}[a_1,a_2,b_1,\dots ,b_{n-2},q]$ by the ideal generated by 
		\begin{equation}\label{eq:QH_Pres.}
			(1+(2a_2-a_1^2)x^2+a_2x^4)(1+b_1x^2+\dots+ b_{n-2}x^{2n-4})=1+qa_1x^{2n}
		\end{equation}
	\end{theorem}
	
	The detailed proof of the above result can be found in \cite{Q_C_IG}. Now we will describe a similar presentation for the orthogonal Grassmannian $\OG(2,N)$, where $N=2n+2$. We will assume $n\ge 3$, otherwise $H^2(\OG(2,N),\mathbb{C})$ may have rank greater than one.
	
	We have the universal exact sequence $0\to \cS\to \mathbb{C}^N\otimes\cO\to \cQ\to 0$ over $\OG(2,N)$. Let $\cS^{\perp}\subset \mathbb{C}^N\otimes\cO$ be the rank $N-2$ vector bundle consisting of vectors perpendicular to $\cS$. 
	
	Unlike the symplectic case, there is a cohomology class which is not obtained using the universal exact sequence. Let $\quadric\subset \mathbb{P}(\mathbb{C}^N)$ be the quadric of isotropic lines in  $\mathbb{C}^N$ equipped with a non-degenerate symmetric bilinear form $\sigma$. Let $\pi :\mathbb{P}(\cS)\to \OG(2,N)$ be the projective bundle. We have the natural the map $\theta : \mathbb{P}(\cS)\to \quadric$. 
	
	Note that $O(2n+2)$ acts on $\mathbb{C}^{2n+2}$. There are precisely two $SO(2n+2)$ orbits of maximal isotropic subspaces. Two maximal isotropic subspaces $E$ and $F$ lie in different orbits if and only if $\dim E\cap F$ is even. Let $e$ and $f$ be the cohomology classes corresponding to $\mathbb{P}(E)$ and $\mathbb{P}(F)$ inside the quadric $\quadric\subset \mathbb{P}(\mathbb{C}^N)$. The classes $e$ and $f$ corresponds to two rulings of $\quadric$.
	
	The cohomology ring of $\quadric$ is generated by the hyper plane class $h$ and ruling classes $e$ and $f$ (see \cite{Char_classes_quadric}).
	
	Over $\OG(2,N)$, we have the following cohomology classes :
	\begin{itemize}
		\item The Chern classes $a_i= c_i(\cS^\vee)$ for $i\in\{1,2\}$.
		\item Let $b_i=c_{2i}(\cS^{\perp}/\cS)$ for $i\in \{1,\dots ,n-1 \}$. The bundle $\cS^{\perp}/\cS$ is self dual, hence all the odd Chern classes vanish.
		\item Let $\pi:\mathbb{P}(\cS)\to \OG$ be the projection, then we define \[\xi = \pi_*\theta^*(e-f). \]
	\end{itemize}
	The above classes still satisfy the identity \eqref{eq:QC_SG_identity}, but two new identities involving $\xi$ are required. We will briefly describe these for readers convenience.
	
	\begin{lem}\label{lem:xi_relation}
		The cohomology class $\xi$ satisfy $\xi a_2=0$ and $\xi^2=(-1)^{n-1}b_{n-1}$.
	\end{lem}
	\begin{proof}
Let $h=c_1(\cO(1))$ on $\mathbb{P}(S)$, then $h\theta^*(e-f)=0$. Multiplying $\theta^*(e-f)$ to the identity  \[h^2-hc_1(\pi^*\cS^\vee)+c_2(\pi^*\cS^\vee)=0,  \]
we obtain $\theta^*(e-f)\pi ^*a_2=0$. The projection formula implies $\xi a_2=0$.

		Using the identities $c_t(\cS)c_t(\cS^{\vee})c_t(\cS^{\perp}/\cS)=1$ and $c_t(\cS)c_t(\cQ)=1$, we obtain $c_t(\cS^{\perp}/\cS)=c_t(\cQ)c_{-t}(\cQ)$. In particular, for all $1\le k\le n-1$
		\begin{equation*}
			(-1)^{k}b_k=c_k(\cQ)^2+2\sum_{i=1}^{k}(-1)^ic_{k+i}(\cQ)c_{k-i}(\cQ).
		\end{equation*}
		When $k=n-1$, the right side of the above equality is $\xi^2$ by \cite{Quantum_pieri_OG}.

	\end{proof}
	\begin{rem}
		The class $\xi$ is the Edidin-Graham characteristic square root class for the quadratic bundle $\cS^\perp/\cS$. 
	\end{rem}
	\begin{prop}
		The cohomology ring $H^*(\OG(2,2n+2))$ is isomorphic to the quotient of the ring $\mathbb{C}[a_1,a_2,b_1,\dots,b_{n-2}, \xi]$ by the ideal generated by the relations $\xi a_2=0$ and 
		\begin{equation*}
			(1+(2a_2-a_1^2)x^2+a_2^2x^4)(1+b_1x^2+\cdots + b_{n-2}x^{2n-4}+(-1)^{n-1}\xi^2x^{2n-2})=1.
		\end{equation*}
	\end{prop}
	\begin{proof}
		Note that the topological Euler characteristic of $\OG$ is the vector space dimension of $H^*(\OG)$ and is given by $2^2\binom{n+1}{2}$. This is obtained by counting the number of fixed points under $\mathbb{C}^*$ action on $\OG$. 
		
		We can unpack the relations to obtain the generators of the ideal:
		\begin{equation}\label{eq:def_fi_OG}
	\begin{aligned}
			f_0&=\xi a_2\\
			f_1&= b_1+(2a_2-a_1^2) \\
			&\vdots
			\\
			f_{n-1}&= (-1)^{n-1}\xi^2+ b_{n-2}(2a_2-a_1^2)+b_{n-3}a_2^2\\
			f_{n}&= (-1)^{n-1}\xi^2(2a_2-a_1^2)+ b_{n-2}a_2^2
	\end{aligned}
		\end{equation}
		Define $R'=\mathbb{C}[a_1,a_2,b_1,\dots,b_{n-2}, \xi]/\langle f_0,\dots,f_n\rangle $. 
		
		Using Lemma \ref{lem:xi_relation} and $c_t(\cS)c_t(\cS^{\vee})c_t(\cS^{\perp}/\cS)=1$, we know that $f_i=0$ for all $0\le i\le n$ in $H^*(\OG)$ . Moreover, the classes $a_1,a_2$ and $\xi$ generates $H^*(\OG)$ (see \cite{Quantum_pieri_OG}). Therefore we get the surjective ring homomorphism \[R'\to H^*(\OG).\] 
		
		It is enough to show that $R'$ is a vector space of dimension at most $2^2\binom{n+1}{2}$. We bound the dimension of $R'$ using the exact sequence \[0\to \langle \xi\rangle \to R'\to R'/\langle \xi\rangle\to 0.  \]
		
		Using \eqref{eq:H*SG(2,2n)}, we observe that $R'/\langle \xi\rangle = H^*(\SG(2,2n))$. Thus $R'/\langle \xi\rangle$ has dimension $2n^2-2n$, which is the Euler characteristic of $\SG(2,2n)$.
		
		Note that $b_i\in a_1^{2i} + \langle a_2\rangle$, $\xi^2\in a_1^{2n-2}+\langle a_2\rangle$ and $\xi^2a_1^2\in \langle a_2\rangle$. Hence $\dim R'/\langle a_2\rangle\le |\{1,a_1\dots,a_1^{2n-1},\xi,\dots \xi a_1^{2n-1} \}|=4n$. Consider the exact sequence
		\[0\to \ker\to R'\xrightarrow{\cdot a_2}R'\to R'/\langle a_2\rangle\to 0. \]
		Note that $\langle \xi\rangle\subset \ker$, thus \[\dim\langle \xi\rangle\le \dim \ker = \dim R'/\langle a_2\rangle\le4n.\]
	\end{proof}

	Now we will turn our attention to the small quantum cohomology.
	\begin{prop}
		Let $n>2$. The small quantum cohomology ring $QH^*(\OG(2,2n+2))$ is isomorphic to the quotient of the ring $\mathbb{C}[a_1,a_2,b_1,\dots,b_{n-2}, \xi,q]$ by the ideal generated by the relations $\xi a_2=0$ and 
		\begin{equation}\label{eq:QH_pres_OG}
			(1+(2a_2-a_1^2)x^2+a_2^2x^4)(1+\cdots + b_{n-2}x^{2n-4}+(-1)^{n-1}\xi^2x^{2n-2})=1+4 qa_1x^{2n}.
		\end{equation}
	\end{prop}
	\begin{proof}
		The degrees of the relations in the given presentation of $H^*(\OG)$ are \[ \deg  f_i = \begin{cases}
			n+1& i=0\\
			2i& 1\le i\le n.
		\end{cases} \]Since $q$ has degree $2n-1$, the quantum term can appear only in degree $2n$ in the above presentation of the cohomology.  Therefore, \[(-1)^{n-1}\xi^2(2a_2-a_1^2) +b_{n-2}a_2^2=cqa_1\] for some constant $c$.  Recall that $(-1)^{n-1}\xi^2=b_{n-1}=c_{2n-2}(\cS^\perp/\cS)$. The first term $\xi^2a_2=0$ since $\xi a_2=0$. Note that we have the following Schubert classes
		\begin{align*}
			b_{n-1}a_1& = c_{2n-1}(\cQ) \\
			b_{n-2}a_2+b_{n-1}&=  c_{2n-2}(\cQ).
		\end{align*}
		It is enough to show that the three point GRW invariants
		\begin{align*}
			\Phi_{0,1} (a_1,c_{2n-1}(\cQ), a_1^*) =2&, \hspace{1cm}	\Phi_{0,1}( a_2,c_{2n-2}(\cQ), a_1^*)=2,
		\end{align*}
		where $a_1^*$ corresponds to the class of a line. It follows by carefully applying the quantum Pieri rule stated in \cite{Quantum_pieri_OG}, which describes the three term genus zero GWR invariants (equivalently the quantum product) of the Schubert classes. 
	\end{proof}
	
\subsection{Jacobian Calculation}\label{subsec:Jacobian_cal}
	We can unpack \eqref{eq:QH_Pres.} to write that the ideal of relations is generated by
	\begin{equation}\label{eq:def_f_i}
		\begin{aligned}
		\tilde{f}_1= & b_1+(2a_2-a_1^2) \\
		\tilde{f}_2=& b_2+b_1(2a_2-a_1^2)+a_2^2\\
		&\vdots
		\\
		\tilde{f}_{n-2}=& b_{n-2}+b_{n-3}(2a_2-a_1^2)+b_{n-4}a_2^2\\
		\tilde{f}_{n-1}=&  b_{n-2}(2a_2-a_1^2)+b_{n-3}a_2^2\\
		\tilde{f}_{n}=& b_{n-2}a_2^2-qa_1.
		\end{aligned}
	\end{equation}
	Let $R=\mathbb{C}[a_1,a_2,b_1,\dots,b_{n-2},q]/\langle \tilde{f}_1,\dots,\tilde{f}_n\rangle$  be the quantum cohomology ring of $\SG(2,2n)$ over $\mathbb{C}[q]$.
	
	In order to calculate the Gromov-Ruan-Witten invariants, we are required to compute the Jacobian
	\begin{equation*}
		J=\det\begin{bmatrix}
			\frac{\partial \tilde{f}_1}{\partial a_1}&\dots &\frac{\partial \tilde{f}_n}{\partial a_1}\\
			\vdots& & \vdots\\
			\frac{\partial \tilde{f}_1}{\partial b_{n-2}}&\dots& \frac{\partial \tilde{f}_n}{\partial b_{n-2}}
		\end{bmatrix}
	\end{equation*}
	at the vanishing locus of $(\tilde{f}_1, \tilde{f}_2,\dots,\tilde{f}_n)$. Substituting $b_1=(a_1^2-2a_2)$, this determinant equals 
	\begin{equation*}
		-4a_1\det \begin{bmatrix}
			1&b_1&b_2&b_3&\dots &b_{n-2}&\frac{q}{2a_1}\\
			1&(a_2+b_1)&(a_2b_1+b_2)&(a_2b_2+b_3) &\dots &(a_2b_{n-3}+b_{n-2})&a_2b_{n-2}\\
			1&-b_1&a_2^2&0&\dots&0&0\\
			0&1&-b_1&a_2^2&\dots&0&0\\
			0&0&1&-b_1&\dots&0&0\\
			\vdots&\vdots & \vdots\\
			0&0&0&\dots &1&-b_1&a_2^2
		\end{bmatrix}.
	\end{equation*}
	After subtracting first two rows, we observe that the above equals
	\begin{equation*}
		-4a_1a_2\det
		\begin{bmatrix}
			1&b_1&b_2&b_3&\dots &b_{n-2}&\frac{q}{2a_1}\\
			0&1&b_1&b_2 &\dots &b_{n-3}&b_{n-2}-\frac{q}{2a_1a_2}\\
			1&-b_1&a_2^2&0&\dots&0&0\\
			0&1&-b_1&a_2^2&\dots&0&0\\
			0&0&1&-b_1&\dots&0&0\\
			\vdots&\vdots & \vdots\\
			0&0&0&\dots &1&-b_1&a_2^2
		\end{bmatrix}.
	\end{equation*}
	Let $v_0, v_1,\dots , v_{n-1}$ be the column vectors in the above matrix. Then the 
	\begin{align*}
		\det[v_0,\dots,v_{n-1}]= \det [V_0,\dots V_{n-1} ]
	\end{align*}
	where $V_i=v_ib_0+v_{i-1}b_1+\cdots + v_0b_i$. Using the identity, $a_2^2b_{i-2}-b_1b_{i-1}+b_i=0$, we observe that 
	\begin{align*}
		[V_0,\dots,V_{n-1}] =
		\begin{bmatrix}
			1&B_1&B_2&B_3&\dots &B_{n-2}&B_{n-1}+\frac{q}{2a_1}\\
			0&1&B_1&B_2 &\dots &B_{n-3}&B_{n-2}-\frac{q}{2a_1a_2}\\
			1&0&0&0&\dots&0&0\\
			0&1&0&0&\dots&0&0\\
			\vdots&\vdots & \vdots\\
			0&0&0&\dots &\ \ 1&0&0
		\end{bmatrix}
	\end{align*}
	where $b_{n-1}:=0$ and $B_i:=b_ib_0+b_1b_{i-1}+b_2b_{i-2}+\dots +b_0b_i$. Therefore the required Jacobian is given by
	\begin{align}\label{eq:Jacobian_SG_Bi}
		J=-4a_1a_2\det \begin{bmatrix}
			B_{n-2}&B_{n-1}+\frac{q}{2a_1}\\
			B_{n-3}&B_{n-2}-\frac{q}{2a_1a_2}.
		\end{bmatrix}.
	\end{align}

		\subsection{Residues}
	We will use the presentation of the quantum cohomology in \eqref{eq:QH_Pres.} and \eqref{eq:QH_pres_OG} to obtain the higher genus GRW invariants for $\SG(2,2n)$ and $\OG(2,2n+2)$  using the techniques in \cite{Q_C_Fano}. We will briefly describe the result we require from \cite{Q_C_Fano}.

 Let $F\in \mathbb{C}[x_1,\dots,x_n]$ be a polynomial, and $f=(f_1,\dots,f_n):\mathbb{C}^n\to\mathbb{C}^n$ be a tuple of polynomials such that $f^{-1}(0)$ is finite. For any $p\in f^{-1}(0)$, we define
	\begin{equation*}
		Res_f(p;F):=\frac{1}{(2\pi i)^n}\int_{\Gamma^{\epsilon}_p}^{}\frac{F}{f_1\cdots f_n}dx_1\dots dx_n
	\end{equation*}
	with $\Gamma^{\epsilon}_p=\{q\in U(p): |f(q)|=\epsilon  \}$, $U(p)$ small neighborhood of $a$ with $f^{-1}(0)\cap U(p)=\{p\}$ and $\Gamma^{\epsilon}_p$ relatively compact in $U(p)$. We may further define \[Res_f(F) =\sum_{p\in f^{-1}(0)}^{}Res_f(p;F).\] 
	Note that when $p$ is a regular point, i.e. the Jacobian $J=\det\big(\partial f_i/\partial x_j \big)\ne 0$ at $p$, then 
	\[Res_f(p;F)= \bigg(\frac{F}{J}\bigg)(p) .\]

Let $M$ be a Fano manifold with $h^2(M,\mathbb{C})=1$ and the cohomology ring  $H^*(M,\mathbb{C})=\mathbb{C}[x_1,\dots,x_n]/\langle f_1,\dots,f_n\rangle $, where each $x_i$ corresponds to a pure dimensional cohomology class. Let \[QH^*(M,\mathbb{C})=\mathbb{C}[x_1,\dots,x_n,q]/\langle \tilde{f}_1,\dots,\tilde{f}_n \rangle\] be the quantum cohomology as an algebra over $\mathbb{C}[q]$. 

Substitute $q$ for a complex number, and let $\tilde{f}^{q}=(\tilde{f}^q_1,\dots,\tilde{f}^q_n)$ be the corresponding tuple of polynomials in $x_1,\dots,x_n$. Let $R_q=QH_q^*(M,\mathbb{C})$ be the corresponding quantum cohomology ring. Note that $R_q$ and $H^*(M,\mathbb{C})$ are isomorphic as vector spaces. The ring $R_q$ is equipped with a quantum multiplication that matches the usual multiplication of cohomology classes when $q=0$.

	\begin{theorem}\label{thm:residue_formula}\cite{Q_C_Fano}
	Let $M$ and $\tilde{f}^q$ be defined as above. Let $F\in \mathbb{C}[x_1,\dots,x_n]$ be a weighted homogeneous polynomial satisfying the dimension condition \eqref{eq:exp_dim} for a natural number $d$. Then
		\begin{align*}
			\langle F\rangle_gq^{d}=c^{\bar{g}} \text{Res}_{\tilde{f}^q}(J_q^gF)=\lim\limits_{\substack{y\to 0} }\sum_{x\in {(\tilde{f}^q)}^{-1}(y)}^{} ((cJ_q)^{\bar{g}}F)(x) 
		\end{align*}
where the limit is taken over regular points $y$, $c$ is a constant and $J_q=\det\big(\partial \tilde{f}^q_i/\partial x_j \big)$ is the Jacobian.
	\end{theorem}

\subsection{GRW invariants for $\SG(2,2n)$}
	We will the apply Theorem \ref{thm:residue_formula} to the presentation of the quantum cohomology $R=QH^*(\SG(2,2n))$ in \eqref{eq:QH_Pres.}. To be precise, let $(x_1,x_2,x_3,\dots,x_n)=(a_1,a_2,b_1,\dots,b_{n-2})$ and let $\tilde{f}$ defined by \eqref{eq:def_f_i}. 

Fix $q=-1$ (or any non-zero number). Equation \eqref{eq:QH_Pres.} can be rephrased as 
	\begin{equation*}
		(z^2-z_1^2)(z^2-z_2^2)Q(z)=z^{2n}+q(z_1+z_2)
	\end{equation*}
	where $a_1=z_1+z_2$, $a_2=z_1z_2$ and $Q(z)=z^{2n-4}+b_1z^{2n-6}+\dots + b_{n-2}$. Observe that $b_i$ can be represented in terms of $a_1$ and $a_2$ for all $1\le i\le n-2$.
	
	Evaluating at $z_1$ and $z_2$, we obtain 
	\begin{align*}
		z_1^{2n}=-q(z_1+z_2)\\
		z_2^{2n}=-q(z_1+z_2).
	\end{align*}

The structure of $R_q$ is  described in \cite{Q_C_IG}. The set $(\tilde{f}^q)^{-1}(0) $ has two types of points:
\begin{itemize}
\item Reduced points: The points described by the unordered pair $\{z_1,z_2\}$ satisfying
	\begin{align}\label{eq:z_i_chern_roots}
	z_2&= \zeta z_1\\
		z_1&=\omega(1+\zeta)^{\frac{1}{2n-1}},\nonumber
	\end{align} 
where $\omega^{2n-1}=-q$, $\zeta^{2n}=1$ and $\zeta\ne \pm 1$. Since $\{z_1,z_2\}$ is an unordered, $(\omega,\zeta)$ and $(\omega,\zeta^{-1})$ yields the same point. Thus there are $(n-1)(2n-1)$ such points. The non-vanishing of the Jacobian computed below implies that these points are reduced.
\item Fat point : The origin is the only other point in $(\tilde{f}^q)^{-1}(0) $. Since the vector space dimension $\dim(R_q)=2n(n-1)$, the origin is a non-reduced point of order $(n-1)$ in $\text{Spec}(R_q)$.
\end{itemize}
Thus $R_q=A_1\times A_2$ where $A_1\cong\mathbb{C}[\epsilon ]/\langle \epsilon^{n-1}\rangle$ corresponds to the fat point at origin in $Spec(R_q)$ and $\text{Spec}(A_2)$ consists of $(n-1)(2n-1)$ distinct reduced points.

\begin{prop}\label{prop:Jacobian_SG}
Let $p\in A_2$ be a reduced point described using \eqref{eq:z_i_chern_roots}. The Jacobian at $p$ is
\begin{equation}
J_q(p)=2n(2n-1)\zeta^{-1}(1+\zeta)^{-1}(1-\zeta)^{-2}z_1^{4n-5}.
\end{equation}
\end{prop}
	\begin{proof}
 We  recursively calculate a concise expression for $b_1,\dots,b_{n-2}$:
	\begin{align*}
		b_i&=z_1^{2i}(1+\zeta^2+\dots +\zeta^{2i}).
	\end{align*}
	 We define $b_i$ for all $i\in \mathbb{N}$ using the above identity. Note that $b_{n-1}=0$ and $b_0=1$.

We are now going to give a simple formula for the convolution products $B_i$, and use it to find the Jacobian.

 Let $t=z_1^2$. Let $P(x)=1+b_1x+b_2x^2+\cdots $ be the power series in $x$. Then 
		\begin{align*}
			(1-\zeta^2)P(x)&= \sum_{i=0}^{\infty} (1-\zeta^{2i+2})(tx)^i\\
			&= \frac{1}{1-tx}-\frac{\zeta^2}{1-\zeta^2tx}.
		\end{align*}
		
		Observe that $P(x)^2=1+B_1x+B_2x^2+\cdots $, which can be expressed as 
		\[P(x)^2 =\frac{1}{(1-\zeta^2)^2}\bigg(\frac{1}{(1-tx)^2}+ \frac{\zeta^4}{(1-\zeta^2tx)^2}-\frac{2\zeta^2}{1-\zeta^2}\bigg(\frac{1}{1-tx}-\frac{\zeta^2}{1-\zeta^2tx} \bigg) \bigg). \]
		Extracting the coefficient of $x^i$ in the above expression gives
		\begin{align*}
			B_i&=\frac{1}{(1-\zeta^2)^2}\bigg((i+1)t^i+(i+1)\zeta^{2i+4}t^i-\frac{2\zeta^2}{1-\zeta^2}(t^i-\zeta^{2i+2}t^i ) \bigg)\\
&=\bigg(\frac{(i+1)(1+\zeta^{2i+4})}{(1-\zeta^2)^2} - \frac{2\zeta^2(1-\zeta^{2i+2})}{(1-\zeta^2)^3} \bigg)t^{i}.
		\end{align*}
	In particular, we have
	\begin{align*}
		&B_{n-1}=n\frac{1+\zeta^2}{(1-\zeta^2)^2}t^{n-1}, \hspace{.5cm}
		B_{n-2}=\frac{2n}{(1-\zeta^2)^2}t^{n-2},\\
		&B_{n-3}= \frac{n(1+\zeta^2)}{\zeta^2(1-\zeta^2)^2}t^{n-3}.
	\end{align*}
	Substituting $q=b_{n-2}a_2^2/a_1$ and using $a_1^2=t(1+\zeta)^2$, $b_{n-2}=-t^{n-2}/\zeta^2$ and $a_2=t\zeta$ we get the expression for Jacobian for $\tilde{f}^q=(\tilde{f}^q_1,\tilde{f}^{q}_{2},\dots,\tilde{f}^q_{n})$ at $p$:
	\begin{align*}
		J_q(p)&=-4a_1a_2\bigg(\det\begin{bmatrix}
			B_{n-2}&B_{n-1}\\ B_{n-3}&B_{n-2}
		\end{bmatrix}+\det \begin{bmatrix}
			B_{n-2}& \frac{b_{n-2}a_2^2}{2a_1^2}\\ B_{n-3}& -\frac{b_{n-2}a_2}{2a_1^2}
		\end{bmatrix}  \bigg)\\
		&=-4a_1a_2\bigg(-\frac{n^2}{\zeta^2(1-\zeta^2)^2} +\frac{n}{2\zeta^2(1-\zeta^2)^2}  \bigg)t^{2n-4}\\
		&=2n(2n-1)\zeta^{-1}(1+\zeta)^{-1}(1-\zeta)^{-2}z_1^{4n-5}.
	\end{align*}
	\end{proof}

	\begin{prop}\label{prop:residue_SG_reduced}
		Let $\vd=(2n-1)d-\bar{g}(4n-5)$ and $F=a_1^{m_1}a_2^{m_2}$ such that $m_1+2m_2=\vd$, then
		\begin{align}
			\sum_{p\in A_2}^{}\text{Res}_{\tilde{f}^q}(p;J_q^{g}F)= \frac{2n-1}{2} \sum_{\zeta\ne \pm 1}^{}(1+\zeta)^{m_1}\zeta^{m_2}J(\zeta)^{\bar{g}}(1+\zeta)^d(-q)^d
		\end{align}
		where $\zeta\ne \pm 1 $ is an $2n^{\text{th}}$ root of unity and $J(\zeta):=2n(2n-1)\zeta^{-1}(1+\zeta)^{-1}(1-\zeta)^{-2}$.
	\end{prop}
	\begin{proof}
		Let $p$ be given by $(\omega,\zeta)$. Using Proposition \ref{prop:Jacobian_SG}
		\begin{align*}
			Res_{\tilde{f}^q}(p;J_q^{g}F)&= (J_q^{g-1}F)(p)\\
			&= J(\zeta)^{\bar{g}}(1+\zeta)^{m_1}\zeta^{m_2}z_1^{\vd + \bar{g}(4n-5)}. 
		\end{align*}
		Observe that $z_1^{\vd + \bar{g}(4n-5)}=(1+\zeta)^d(-q)^d,$ thus 
\begin{align*}
	\sum_{p\in A_2}^{}\text{Res}_{\tilde{f}^q}(p;J_q^{g}F)= \sum_{(\omega,\zeta)}^{}(1+\zeta)^{m_1}\zeta^{m_2}J(\zeta)^{\bar{g}}(1+\zeta)^d(-q)^d
\end{align*}
where the latter is summed over pairs $(\omega,\zeta)$ such that $\omega^{2n-1}=(-q)$ and $\zeta$ is a $2n^{\text{th}}$ root of unity with strictly positive imaginary part.
		The above expression does not depend on the choice of $\omega$ and it is invariant under $\zeta\to \zeta^{-1}$. When summed over these choices the required formula is obtained.
	\end{proof}
	
	\begin{theorem}\label{thm:GRW_invariant_SG}
		Let $m_1+2m_2=\vd=(2n-1)d-(4n-5)\bar{g}$. The GRW invariants for $\SG(2,2n)$ equal the top virtual intersections of the $a$-classes on the corresponding isotropic Quot scheme:
		\begin{align}
			\langle a_1^{m_1}a_2^{m_2}\rangle_g= 	\int_{[\IQ_d]^{\vir}}^{}a_1^{m_1}a_2^{m_2}
		\end{align}
	\end{theorem}
	\begin{proof}
		 The origin $y=0:=(0,\dots, 0)$ is not necessarily a regular point for the function $\tilde{f}^q=(\tilde{f}^q_1,\dots,\tilde{f}^q_n)$. We will evaluate the limit
\begin{equation}
\lim\limits_{y\to 0}\sum_{p\in {(\tilde{f}^q)}^{-1}(y)}^{} (J^{\bar{g}}F)(p) ,
\end{equation}  where the limit  $y\to 0$ is taken over regular values of $y$. Let $ \epsilon$ be a non-zero complex number with small absolute value, and let $y_{\epsilon}=(0,\dots,0,\epsilon^{n-1},0)$. We will see that $y_{\epsilon}$ is regular for $\epsilon$ small enough.

Reduced points : Since the Jacobian for each point $p\in A_2$ is non-zero, the inverse function theorem implies that for small enough $\epsilon$, there is exactly one reduced point $p_\epsilon$ near $p$ satisfying $f(p_\epsilon)=y_\epsilon$. Thus $y_{\epsilon}$ is a regular value for all $\epsilon$ in a neighborhood of $0$.

Let $A_2^\epsilon$ be the set of unique points $p_\epsilon$ near $p\in A_2$. Observe that the residue contribution is
\begin{equation}
\lim\limits_{\epsilon\to 0}\sum_{p_\epsilon\in A_2^{\epsilon}}^{} (J^{\bar{g}}F)(p_\epsilon) =	\sum_{p\in A_2}^{}\text{Res}_f(p;J^{g}F).
\end{equation}
This has been calculated in Proposition \ref{prop:residue_SG_reduced}.

Fat point : The vanishing of $\tilde{f}^q_1,\dots,\tilde{f}^q_{n-2}$ implies that $b_1,\dots,b_{n-2}$ is a polynomial in $a_1$ and $a_2$. Observe that \[b_i=(-1)^i(i+1)a_2^i+\langle a_1^2\rangle.  \] Since $q\ne0$, the vanishing of $\tilde{f}^q_n$ implies \[a_1=q^{-1}a_2^{n}+\langle a_1^2\rangle.  \]
Therefore $a_1=a_2^{n}h_1(a_2)$ for some power series $h$ that defines a holomorphic function for an open set containing $0$. A similar argument shows that $f_{n-1}=a_2^{n-1}h_2(a_2)$ where $h_2$ is holomorphic with non-zero constant term. Observe that $a_2^{n-1}h_2(a_2)=\epsilon\ne0$ has exactly $(n-1)$ simple zeros for all $\epsilon$ lying in a neighborhood of $0$.

Note that $a_2=O(\epsilon)$, $a_1=O(\epsilon^n)$ and $b_i=O(\epsilon^i)$ as $\epsilon$ approaches $0$. Substituting the above orders in \eqref{eq:Jacobian_SG_Bi}, we get $J=O(\epsilon^{n-2})$. Thus the residue contributions of these $n-1$ points has order $O(\epsilon^{nm_1+m_2+\bar{g}(n-2)})$, which vanishes in the limit $\epsilon\to 0$ when the the exponent $nm_1+m_2+\bar{g}(n-2)$ is non-zero. 

There are exactly two cases when the above exponent is zero: (i) $\vd=0$, $d=g-1$, $N=2n=4$; and (ii) $\vd=d=0$, $g=1$. An easy calculation shows that the residue contribution are $(2q)^{d}$ and $1$ respectively. These are the only instances where $\vd \ge 0$ and $d<g$. 

We apply Theorem \ref{thm:residue_formula} to obtain the GRW invariant up to a constant $c$. When $g=d=0$, the GRW invariants are the top intersections in the cohomology ring of $\SG(2,2n)$. Note that $\IQ_{0}\cong \SG(2,2n)$ when $g=0$, thus the virtual invariants in \eqref{eq:a_1a_2formula} must match the GRW invariants. Comparing the two we obtain $c=-1$. 

Putting together all the terms, we get
\begin{align*}
	\langle a_1^{m_1}a_2^{m_2}\rangle_g=
\begin{cases}
(-1)^{d+\bar{g}}\frac{2n-1}{2} \sum_{\zeta}^{}(1+\zeta)^{m_1+d}\zeta^{m_2}J(\zeta)^{\bar{g}} &d\ge g\\
2^{\bar{g}}3^g+(-1)^{\bar{g}}2^{d} & n=2,\ d=\bar{g}\\
2n(n-1) & g=1, d=0
\end{cases}.
\end{align*}
This match the expression in Theorem \ref{thm:r=2,sympl} (also see Examples \ref{exm:N=4} and \ref{exm:g=1}) for all $d$, $g$ and $N$.
	\end{proof}
	\subsection{GRW invariants for $\OG(2,2n+2)$}
Let $n\ge 3$. Recall the definition of ${f}_0,{f}_1,\dots,{f}_n$ from \eqref{eq:def_fi_OG}. Let $\tilde{f}_i=f_i$ for $0\le i\le n-1$ and let $\tilde{f}_n=f_n-4qa_1$ as prescribed by \eqref{eq:def_fi_OG}. In particular, 
\begin{align*}
			\tilde{f}^q_0&=\xi a_2\nonumber\\
			\tilde{f}^q_1&= b_1+(2a_2-a_1^2) \nonumber\\
			&\vdots
			\\
			\tilde{f}^q_{n-1}&= (-1)^{n-1}\xi^2+ b_{n-2}(2a_2-a_1^2)+b_{n-3}a_2^2 \nonumber\\
			\tilde{f}^q_{n}&= (-1)^{n-1}\xi^2(2a_2-a_1^2)+ b_{n-2}a_2^2 - 4qa_1 \nonumber
		\end{align*}

Let $R'=\mathbb{C}[\xi,a_1,a_2,b_1,\dots,b_{n-2},q]/\langle \tilde{f}_0,\dots,\tilde{f}_n\rangle$ be the presentation for the quantum cohomology of $\OG(2,2n+2)$ (see \eqref{eq:QH_pres_OG}). The Jacobian $J'$ for $\tilde{f}=(\tilde{f}_0,\dots,\tilde{f}_n)$ is calculated in similar fashion as it was done in the symplectic case. Observe that
\begin{equation}
J'\in -4a_1a_2^2\det \begin{bmatrix}
			B_{n-2}&B_{n-1}+\frac{4q}{2a_1}\\
			B_{n-3}&B_{n-2}-\frac{4q}{2a_1a_2}
		\end{bmatrix}+\langle \xi\rangle,
\end{equation}
where $b_0=1$, $b_{n-1}:=(-1)^{n-1}\xi^{2}$ and $B_i=b_ib_0+\cdots+b_0b_i$.

Note that modulo $\langle a_2\rangle$, we have
\begin{align*}
\tilde{f}_0&=0\\
\tilde{f}_1&=b_1-a_1^2\\
\vdots&\\
\tilde{f}_{n-1}&=(-1)^{n-1}\xi^2-b_{n-2}a_1^2\\
\tilde{f}_{n}&= (-1)^{n-1}\xi^2(-a_1^2)-4qa_1
\end{align*}
An easy calculation shows that
\begin{equation*}
J'\in-2b_{n-1}(2a_1B_{n-1} +4q )+\langle a_2\rangle.
\end{equation*}
Note that $b_i\in a_1^{2i}+\langle a_2\rangle$, thus we may further write 
\begin{equation}\label{eq:Jacobian_OG_a_1}
J'\in -2a_1^{2n-2}(2na_1^{2n-1}+4q)+\langle a_2\rangle . 
\end{equation}

Fix a non-zero number $q$. Note that $f_0=0$ implies that either $\xi=0$ or $a_2=0$. The set $(\tilde{f}^{q})^{-1}(0)$ has three types of points:
\begin{itemize}
\item Reduced points ($a_2\ne0$): The reduced points with $\xi=0$ have almost the same description as that of $Spec(A_2)$ in the symplectic case. It is obtained by replacing $q\to 4q$ and letting $a_1$ and $a_2$ be described (similar to \eqref{eq:z_i_chern_roots}) using Chern roots $\{z_1,z_2\}$ in this case.
\item Reduced points ($\xi\ne0$):  Thus $a_2=0$ and hence $b_i=a_1^{2i}$. Moreover,  $\tilde{f}^q_{n-1}=\tilde{f}^q_{n}=0$ implies
\begin{align*}
(-1)^{n-1}\xi^{2}&=a_1^{2n-2}\\
a_1^{2n}&=-4qa_1.
\end{align*}
Thus there are $(4n-2)$ points given by $(\xi,a_1)=(\sqrt{-4q}\mu^{-1},\mu^2)$ where $\mu $ is a $(4n-2)^{\text{th}}$ root of $(-4q)$. We observe that the Jacobian (see \eqref{eq:Jacobian_OG_a_1}) is non-zero.
\item Fat point $A_1$: The origin is the non-reduced point of order $(n+1)$.
\end{itemize}
The Artinian ring $R'_q$ is isomorphic to $A_1\times A_2\times A_3$ where $A_1\cong \mathbb{C}[\epsilon]/\langle \epsilon^{n+1} \rangle$. The Spec of $A_2$ and $A_3$ corresponds to the distinct reduced points with $a_2\ne0$ and $\xi\ne 0$ respectively. 

Over the points $p\in Spec(A_2)$ given by a choice of  $\{z_1,z_2\}$ as defined in \eqref{eq:z_i_chern_roots} by replacing $q\to 4q$, the Jacobian
	\begin{align*}
		J_q'(p)=2n(2n-1)(1+\zeta)^{-1}(1-\zeta)^{-2}z_1^{4n-3}.
	\end{align*}
We obtain an analogue of Proposition \ref{prop:residue_SG_reduced}:
\begin{prop}
Let $\vd=(2n-1)d-\bar{g}(4n-3)$ and $F=a_1^{m_1}a_2^{m_2}$ such that $m_1+2m_2=\vd$, then
		\begin{align}
			\sum_{p\in A_2}^{}\text{Res}_{\tilde{f}^q}(p;J'^{g}F)= \frac{2n-1}{2}\sum_{\zeta\ne \pm 1}(1+\zeta)^{m_1+d}\zeta^{m_2}J'(\zeta)^{\bar{g}}(-4q)^d
		\end{align}
		where $\zeta\ne \pm 1 $ is $2n^{\text{th}}$ root of unity and $J'(\zeta):=2n(2n-1)(1+\zeta)^{-1}(1-\zeta)^{-2}$.
\end{prop}

\begin{prop}
Let $F=a_1^{m_1}a_2^{m_2}$, where $m_1+2m_2=\vd$. Then
	\begin{align}
			\sum_{p\in A_3}^{}\text{Res}_{\tilde{f}^q}(p;J'^{g}F)=\begin{cases}
(-1)^{\bar{g}}(4n-2)^{g} (-4q)^d & m_2=0\\
0& m_2>0
\end{cases}.
		\end{align}
\end{prop}
\begin{proof}
Let $p\in A_3$ be determined by $(\xi,a_1)=(\sqrt{-4q}\mu^{-1},\mu^2)$ where $\mu$ is a $(4n-2)^{\text{th}}$ root of unity.  Note that $a_2=0$, thus the residues vanish when $m_2>0$. 

We may assume $m_2=0$. Using \eqref{eq:Jacobian_OG_a_1} and the equality $a_1^{2n-1}+4q=0$, the Jacobian is $-2a_1^{4n-3}(2n-1)$. Thus 
\begin{align*}
\text{Res}_{\tilde{f}^q}(p;J'^{g}a_1^{\vd})
&=(-1)^{\bar{g}}(2(2n-1))^{\bar{g}} a_1^{(2n-1)d}\\
&=(-1)^{\bar{g}}(4n-2)^{\bar{g}} (-4q)^d.
\end{align*} 
\end{proof}
	\begin{theorem}
		Let $m_1+2m_2=(2n-1)d-(4n-3)\bar{g}$ and $n\ge 3$. The GRW invariants for $\OG(2,2n+2)$ involving $a_1$ and $a_2$  equal the top virtual intersections of the $a$-classes on the corresponding isotropic Quot schemes.

In particular, when $d\ge g$ and
\begin{itemize}\label{thm:GRW_invariants_OG}
\item[(i)] When $m_2>0$, then
\begin{align*}
			\langle a_1^{m_1}a_2^{m_2}\rangle_g=u4^{d}
\frac{2n-1}{2}\sum_{\zeta\ne \pm 1}(1+\zeta)^{m_1+d}\zeta^{m_2}\bigg(\frac{J'(\zeta)}{4}\bigg)^{\bar{g}},
		\end{align*}
where 	$u=(-1)^{\bar{g}+d}$ and $J'(\zeta)=2n(2n-1)(1+\zeta)^{-1}(1-\zeta)^{-2}$.
\item[(ii)] When $m_2=0$, then
\begin{align*}
\langle a_1^{m_1}\rangle_g=u4^{d}\bigg(\frac{(-1)^{\bar{g}}(4n-2)^{\bar{g}}}{4^{\bar{g}}}+
\frac{2n-1}{2}\sum_{\zeta\ne \pm 1}\frac{(1+\zeta)^{m_1+d}J'(\zeta)^{\bar{g}}}{4^{\bar{g}}}  \bigg).
\end{align*}
\end{itemize}
	\end{theorem}
The proof of the above theorem is similar to that of Theorem \ref{thm:GRW_invariant_SG}.

	\section{Intersection of $f$ classes}\label{sec:f_intersection}

	  We will find an explicit expression for the intersection numbers of polynomials in $a$ and $f$ classes in terms of multivariate generating functions. We obtain Theorem \ref{thm:f_classes_d>g} as a corollary. While the computations are more involved, the basic ideas are similar to those in Section \ref{sec:Int_a_classes}.

	We will only work with symplectic isotropic Quot scheme $\IQ_{d}$ with $r=2$. A similar analysis can be carried out when $\sigma$ is symmetric.

Over the fixed loci $\fix_{\vec{d},\underline{k}}$, the equivariant restriction of the $f$ classes are given by $f_1=d$ and $f_2=\phi_{12}+d_1(x_2+w_2t)+d_2(x_1+w_1t)$. The formula for the intersection of $f$ classes with a polynomial in $a$ classes involves differential operators.
	
	Let $P(X)=X^N-1$ and
	\begin{equation*}
		T_{g}(t,Y_1,Y_2)=\bigg(\prod_{i=1}^{2} (1-\eta_i)-\prod_{i=1}^{2}t^2\eta_i\bigg)^g,
	\end{equation*}
	where $\eta_i=\frac{P(Y_i)}{P'(Y_i)(Y_1+Y_2)}$. When  $Y_i=w_i(1+q_i)^{\frac{1}{N}}$, $T_{g}(t,Y_1,Y_2)$ is a power series in $q_1$ and $q_2$ over $\mathbb{C}[t]$. This should be considered as an analogue of $T_{d,g}(N)$ in \eqref{eq:T_d,g}. In particular,
	\begin{equation*}
		T_{g}(1,w_1(1+q)^{\frac{1}{N}},w_2(1+q)^{\frac{1}{N}})= \bigg(1-\frac{q}{N(1+q)}\bigg)^g.
	\end{equation*}

Let $\partial_i$ and $\partial_t$ be the partial derivatives with respect to $Y_i$ and $t$ respectively. Define the differential operators $\mathfrak{d}_t=-(Y_1+Y_2)\partial_t$, 
		\begin{align*}
			\Delta^u&:=\sum_{i=0}^{u}\binom{u}{i}(q_1\partial_1)^i(q_2\partial_2)^{u-i}Y_2^iY_1^{u-i},\\
			(\Delta+\mathfrak{d}_t)^m&:=\sum_{u=0}^{m}\binom{m}{u}\Delta^u\mathfrak{d}_t^{m-u}.
		\end{align*}
		Note that $\Delta^u$ defined above is not $u^{\text{th}}$ power of the operator $\Delta$.
	
	\begin{theorem}\label{thm:f-classes}
		Let $Q(X_1,X_2)$ be a weighted homogeneous polynomial and $m$ be a positive integer satisfying $\vd=m+\deg Q$, where $\deg Q$ is the weighted degree. Then 
		\begin{align*}
			\int_{[\IQ_d]^{\vir}}f_2^{m}Q(a_1,a_2)=\sum_{w_1,w_2}^{}[q^{d}](\Delta+\mathfrak{d}_t)^m B(Y_1,Y_2) T_{g}(t,Y_1,Y_2)\bigg|_{t=1, q=q_1=q_2}
		\end{align*}
		where the sum is taken over $N^{\text{th}}$ roots of unity $\{w_1,w_2\}$ such that $w_1\ne \pm w_2$, $u=(-1)^{\bar{g}+d}$, $Y_i=w_i(1+q_i)^{1/N}$ and
		\begin{equation*}
			B(Y_1,Y_2)=  uQ(Y_1+Y_2, Y_1Y_2)\frac{(Y_1+Y_2)^{d-\bar{g}}}{(Y_1-Y_2)^{2\bar{g}}}\prod_{i=1}^{2} P'(Y_i)^{\bar{g}}.
		\end{equation*}

	\end{theorem}
	\begin{proof}
		Using the same arguments as in the proof of Theorem \ref{thm:7.1_r=2_sympl}, we see that the required intersection number equals
		\begin{align*}
			\sum_{w_1,w_2}^{}\sum_{|\vec{d}|=d}^{}\sum_{k=0}^{m}\binom{m}{k} \int_{\fix_{\vec{d},\underline{k}}}^{}\phi_{12}^{k}(d_1Y_2+d_2Y_1)^{m-k} R(Y_1,Y_2)e^{-\frac{\theta_1+\theta_2-\phi_{12}}{Y_1+Y_2}}\prod_{i=1}^{2}e^{\theta_iz_i}h_i^{d_i-\bar{g}},
		\end{align*}
		where $z_i=\frac{P'(Y_i)}{P(Y_i)}-\frac{1}{x_i}$ and $h_i=\frac{x_i}{P(Y_i)}$ and \[R(Y_1,Y_2)=uQ(Y_1+Y_2, Y_1Y_2)\frac{(Y_1+Y_2)^{d-\bar{g}}}{(Y_1-Y_2)^{2\bar{g}}}. \]
		We pursue this calculation in Subsection \ref{sec:furth_inter}, in particular we use Proposition \ref{prop:f_int} to finish the proof.
		
	\end{proof}
	When $m=0$, we recover Theorem \ref{thm:r=2,sympl}. We specialize to the case $m=1$ to obtain a simple expression. 
	\begin{cor}
		Recall the definition of $T_{d,g}(N)$ from Theorem \ref{thm:r=2,sympl}. Let $Q$ be a homogeneous polynomial such that $\vd=m+\deg Q$, where $\deg Q$ is the weighted degree. Then 
		\begin{align*}
			\int_{[\IQ_d]^{\vir}}f_2Q(a_1,a_2)=&\frac{2}{N}\sum_{w_1,w_2}^{}\bigg(T_{d-1,g}(N) D\circ B(w_1,w_2)+\\ &\frac{1}{N}\frac{w_1w_2B(w_1,w_2)}{(w_1+w_2)}(T_{d-2,\bar{g}}(N)-NT_{d-1,\bar{g}}(N))\bigg)
		\end{align*}
		where $D \circ B(z_1,z_2)=\frac{z_1z_2}{2}\big(\frac{\partial}{\partial z_1}+\frac{\partial}{\partial z_2}\big)B(z_1,z_2)$ and the sum is taken over all the pairs of $N^\text{th}$ roots of unity $\{w_1,w_2\}$ with $w_1\ne \pm w_2$.

		In particular, when $d>g$ we get 
		\begin{align*}
			\int_{[\IQ_d]^{\vir}}f_2Q(a_1,a_2)=&\frac{2}{N}\bigg(1-\frac{1}{N}\bigg)^{g}\sum_{w_1,w_2}^{}\bigg( D\circ B(w_1,w_2)-\frac{w_1w_2B(w_1,w_2)}{(w_1+w_2)}\bigg).
		\end{align*}
	\end{cor}
	\begin{proof}
		Since $B$ is a homogeneous rational function in variables $Y_1$ and $Y_2$ of degree $Nd-1$, substituting $Y_1/w_1=Y_2/w_2=(1+q)^\frac{1}{N}$ gives a constant multiple of $(1+q)^{d-1/N}$. We use product rule to split the calculation.
		
		First we see that
		\begin{align}\label{eq:DeltaB}
			[q^d]T_{g}(t,Y_1,Y_2)\Delta B(Y_1,Y_2)\bigg|_{q_1=q_2=q}= \frac{2}{N}T_{d-1,g}(N) D\circ B(w_1,w_2),
		\end{align}
		since substituting $Y_1/w_1=Y_2/w_2=(1+q)^\frac{1}{N}$ in $\Delta B(Y_1,Y_2) $ gives us a constant times $q(1+q)^{d-1}$. The rest follows from the definition of $D$ and $T_{d,g}(N)$.
		
		Now we will find $[q^d]B(Y_1,Y_2)(\Delta+\mathfrak{d}_t)T_{g}(t,Y_1,Y_2)$. Let us define 
		\begin{equation*}
			T_g(q)=\bigg(1-\frac{q}{N(1+q)}\bigg)^g
		\end{equation*}
		for notational convenience. Note that 
		\begin{align*}
			\mathfrak{d}_tT_{g}(t,Y_1,Y_2)&=-(Y_1+Y_2)gT_{g-1}(t,Y_1,Y_2)(-2t\eta_1\eta_2)
		\end{align*}
		therefore
		\begin{equation*}
			\mathfrak{d}_tT_{g}(t,Y_1,Y_2)|_{t=1,q_1=q=q_2}=2g\frac{w_1w_2}{w_1+w_2}\frac{q^2T_{g-1}(q)}{N^2(1+q)^2}(1+q)^{\frac{1}{N}},
		\end{equation*}
		hence the the corresponding contribution is
		\begin{equation}\label{eq:DetlatT}
			[q^d]B(Y_1,Y_2)\mathfrak{d}_tT_{g}(t,Y_1,Y_2)|_{t=1,q_1=q=q_2}=\frac{2}{N^2}\frac{w_1w_2B(w_1,w_2)}{w_1+w_2}T_{d-2,g-1}(N).
		\end{equation}
		The other term simplifies as
		\begin{align}\label{eq:DeltaT}
			\Delta T_{g}(1,Y_1,Y_2)&=-gT_{g-1}(1,Y_1,Y_2)\big(q_1Y_2(\partial_1\eta_1+\partial_1\eta_2)+q_2Y_1(\partial_2\eta_1+\partial_2\eta_2)\big),
		\end{align}
		where we evaluate the partial derivatives
		\begin{align*}
			\partial_1\eta_1&=\bigg(\frac{1}{Y_1+Y_2}-\frac{P(Y_1)P''(Y_1)}{P'(Y_i)^2(Y_1+Y_2)}-  \frac{P(Y_i)}{P'(Y_1)(Y_1+Y_2)^2}\bigg)\partial_1Y_1\\
			\partial_1\eta_2&= -\frac{P(Y_2)}{P'(Y_2)(Y_1+Y_2)^2}\partial_1Y_1.
		\end{align*}
		Similar expressions hold for  $\partial_2\eta_1$ and $\partial_2\eta_2$. Note that we also know that $\partial_iY_i= \frac{1}{NY_i^{N-1}}=\frac{1}{P'(Y_i)}$. Using this we find the following identities:
		\begin{align*}
			\frac{q_1Y_2}{(Y_1+Y_2)P'(Y_1)}+\frac{q_2Y_1}{(Y_1+Y_2)P'(Y_2)}\bigg|_q&= \frac{2}{N}\frac{w_1w_2}{(w_1+w_2)}\frac{q(1+q)^{\frac{1}{N}}}{(1+q)}\\
			\frac{q_1Y_2P(Y_1)P''(Y_1)}{(Y_1+Y_2)P'(Y_1)^3}+\frac{q_2Y_1P(Y_2)P''(Y_2)}{(Y_1+Y_2)P'(Y_1)^3}\bigg|_q&= \frac{2(N-1)}{N^2}\frac{w_1w_2}{(w_1+w_2)}\frac{q^2(1+q)^{\frac{1}{N}}}{(1+q)^2}\\
			\frac{q_1Y_2P(Y_1)}{(Y_1+Y_2)^2P'(Y_1)^2}+\frac{q_2Y_1P(Y_2)}{(Y_1+Y_2)^2P'(Y_2)^2}\bigg|_q&= \frac{1}{N^2}\frac{w_1w_2}{(w_1+w_2)}\frac{q^2(1+q)^{\frac{1}{N}}}{(1+q)^2}\\
			\frac{1}{Y_1+Y_2}\bigg(\frac{q_1Y_2P(Y_2)}{P'(Y_2)P'(Y_1)}+\frac{q_2Y_1P(Y_1)}{P'(Y_1)P'(Y_2)}\bigg)\bigg|_q&= \frac{1}{N^2}\frac{w_1w_2}{(w_1+w_2)}\frac{q^2(1+q)^{\frac{1}{N}}}{(1+q)^2}.
		\end{align*}
		Substituting the above expressions back in \eqref{eq:DeltaT}, we obtain
		\begin{align*}
			\Delta T_{g}(1,Y_1,Y_2)\bigg|_{q_1=q_2=q}=gT_{g-1}(q)\frac{w_1w_2}{w_1+w_2}\frac{2}{N}\frac{-q}{(1+q)^2}(1+q)^{\frac{1}{N}}.
		\end{align*}
		Therefore
		\begin{align}\label{eq:DeltaTq1=q2}
			[q^d]B(Y_1,Y_2)\Delta T_{g}(1,Y_1,Y_2)|_{,q_1=q=q_2}=\frac{-2}{N}\frac{w_1w_2B(w_1,w_2)}{w_1+w_2}T_{d-1,g-1}(N).
		\end{align}
		We get the required expression by summing \eqref{eq:DeltaB}, \eqref{eq:DetlatT} and \eqref{eq:DeltaTq1=q2}.
	\end{proof}
	
		\subsection{Further calculations over fixed loci}\label{sec:furth_inter}
	The following results are crucially used to obtain Theorem \ref{thm:f-classes}. They are analogue of Proposition \ref{prop:theta_sum_int} and \ref{prop:fixed_loci_calculation}.
	
	\begin{prop}\label{prop:theta_int_f_class}
		Let $R$ be a homogeneous polynomial with weighted degree $Nd-2\bar{g}(N-1)-p-u$. 
		Let $R(Y_1,Y_2)$ be a homogeneous rational function of degree $s=Nd-2\bar{g}(N-1)$. We borrow the notation $X_{\vec{d}}$, $Y_i$, $P(Y)$, $B(Y)$, $h_i$ and $z_i$ from Proposition \ref{prop:theta_sum_int}. Then 
		\begin{align*}
			&\int_{X_{\vec{d}}}^{}(d_1Y_2+d_2Y_1)^uR(Y_1,Y_2)\prod_{i=1}^{2}\frac{\theta_i^{p_i}}{p_i!}e^{\theta_iz_i}h_i^{d_i-\bar{g}}
			\\
			&= [q_1^{d_1}q_2^{d_2}]\Delta^u \bigg(R(Y_1,Y_2)\prod_{i=1}^{2}\binom{g}{p_i} \frac{B(Y_i)^{g-p_i}P(Y_i)^{p_i}}{P'(Y_i)}\bigg)
		\end{align*}
		where $Y_i=w_i(1+q_i)^{\frac{1}{N}}$ as a power series in $q_i$ on the right hand side.
	\end{prop}
	\begin{proof}
		Let $g(x)=\sum a_dx^d$. The generating functions of the form $f(x)=\sum d^ka_dx^d$ can be evaluated as \[f(x)=\bigg(x\frac{\partial}{\partial x}\bigg)^kg(x). \]This holds true for multivariate generating functions (by using partial derivatives). Using the proof of Proposition \ref{prop:theta_sum_int}, specifically equation \ref{eq:power_series_version_theta}, we get the required expression.
	\end{proof}
	\begin{prop}\label{prop:f_int}
		The following identity holds
		\begin{align*}
			&\int_{X_{\vec{d}}}^{}\phi_{12}^k(d_1Y_2+d_2Y_1)^{m-k} R(Y_1,Y_2)e^{-\frac{\theta_1+\theta_2-\phi_{12}}{Y_1+Y_2}}\prod_{i=1}^{2}e^{\theta_iz_i}h_i^{d_i-\bar{g}}
			\\
			&= [q_1^{d_1}q_2^{d_2}]\Delta^{m-k}\mathfrak{d}_t^k F_t(Y_1,Y_2)\bigg|_{t=1}
		\end{align*}
		where $\eta_i=\frac{P(Y_i)}{B(Y_i)(Y_1+Y_2)}$, $\mathfrak{d}_t=-(Y_1+Y_2)\partial_t$ and
		\begin{align*}
			F_t(Y_1,Y_2)&=  R(Y_1,Y_2)\prod_{i=1}^{2} \frac{B(Y_i)^{g}}{P'(Y_i)}	
			\bigg(\prod_{i=1}^{2} (1-\eta_i)-\prod_{i=1}^{2}t^2\eta_i\bigg)^g.
		\end{align*}
	\end{prop}
	\begin{proof}
		Using Proposition \ref{prop:phi_int} we may replace even powers of $\phi_{12}$ with suitable expression in $\theta_i$'s. Therefore we can make the following replacement 
		\begin{align*}
			\phi_{12}^ke^{-\frac{\theta_1+\theta_2- \phi_{12}}{Y_1+Y_2}}
			&\to\sum_{p=0}^{\infty}\frac{(-1)^{p+\ell}}{p!(Y_1+Y_2)^p}\bigg(\sum_{\substack{\ell+r+s=p\\\ell\equiv k\mod 2}}^{}\binom{p}{\ell,r,s}\theta_1^{r}\theta_2^{s}\phi_{12}^{\ell+k}  \bigg)\\
			&\to\sum_{p=0}^{\infty}\sum_{\substack{\ell+r+s=p\\\ell\equiv k\mod 2}}^{}\frac{(-1)^{p+k-\frac{\ell+k}{2}}}{p!}\binom{p}{\ell,r,s}\binom{\ell+k}{\frac{\ell+k}{2}}\binom{g}{\frac{\ell+k}{2}}^{-1}\frac{\theta_1^{r+\frac{\ell+k}{2}}\theta_2^{s+\frac{\ell+k}{2}}}{(Y_1+Y_2)^p} 
		\end{align*}
		We use Proposition \ref{prop:theta_int_f_class} and binomial identities to obtain that the required expression  is
		\begin{align*}
			\sum_{p=0}^{\infty}\sum_{\substack{\ell+r+s=p\\\ell\equiv k\mod 2}}^{}(-1)^{p+k-\frac{\ell+k}{2}}\binom{p}{\ell,r,s}\frac{(p+k)!}{p!}\binom{p+k}{r+\frac{\ell+k}{2}}^{-1}\binom{\ell+k}{\frac{\ell+k}{2}}\binom{g}{\frac{\ell+k}{2}}^{-1}\\
			\cdot \binom{g}{r+\frac{\ell+k}{2}}\binom{g}{s+\frac{\ell+k}{2}}[q_1^{d_1}q_2^{d_2}]\Delta^{m-k} \frac{J(Y_1,Y_2)}{(Y_1+Y_2)^p} \bigg(\frac{Y_1}{h(Y_1)}\bigg)^{r+\frac{\ell+k}{2}}\bigg(\frac{Y_2}{h(Y_2)}\bigg)^{s+\frac{\ell+k}{2}},
		\end{align*}
		where $h(Y_i)=Y_iB(Y_i)/P(Y_i)$ and \[J(Y_1,Y_2)= R(Y_1,Y_2)\prod_{i=1}^{2} \frac{B(Y_i)^{g}}{P'(Y_i)}. \]
		
		The binomial factor simplifies to give us
		\begin{align*}
			[q_1^{d_1}q_2^{d_2}]\Delta^u\sum_{p=0}^{\infty}&\sum_{\substack{\ell+r+s=p\\2| \ell- k}}^{}(-1)^{\frac{\ell+k}{2}}\frac{(k+\ell)!}{\ell!}\binom{g}{\frac{\ell+k}{2}}\binom{g-\frac{\ell+k}{2}}{r} \binom{g-\frac{\ell+k}{2}}{s}\\
			&\cdot\frac{J(Y_1,Y_2)}{(Y_1+Y_2)^p} \bigg(\frac{-Y_1}{h(Y_1)}\bigg)^{r+\frac{\ell+k}{2}}\bigg(\frac{-Y_2}{h(Y_2)}\bigg)^{s+\frac{\ell+k}{2}}
		\end{align*}
		We sum over $r$ and $s$ keeping $\ell $ fixed after pulling out the terms independent of $r,s$ and $\ell$ to obtain 
		\begin{align*}
			[q_1^{d_1}q_2^{d_2}]\Delta^{m-k}(-1)^k(Y_1+Y_2)^k J(Y_1,Y_2)
			\sum_{2|(\ell- k)}^{}\frac{(k+\ell)!}{\ell!}\binom{g}{\frac{\ell+k}{2}}(-1)^{\frac{\ell+k}{2}}
			\\
			\cdot
			\prod_{i=1}^{2}(-\eta_i)^{\frac{\ell+k}{2}}
			(1-\eta_i )^{g-\frac{\ell+k}{2}}	
		\end{align*}
		The result follows by noting that 
		\begin{align*}
			\sum_{2|(\ell- k)}^{}\frac{(k+\ell)!}{\ell!}\binom{g}{\frac{\ell+k}{2}}(-1)^{\frac{\ell+k}{2}}
			\prod_{i=1}^{2}(-\eta_i)^{\frac{\ell+k}{2}}
			(1-\eta_i )^{g-\frac{\ell+k}{2}}= \partial_t^k \bigg(\prod_{i=1}^{2} (1-\eta_i)-\prod_{i=1}^{2}t^2\eta_i\bigg)^g\bigg|_{t=1}.
		\end{align*}
	\end{proof}

\section{Virtual Euler characteristics}\label{sec:vir_Euler_char}

The Euler characteristic of the symmetric product of curves is given by the well known formula
	\begin{align*}
		e(C^{[d]}) 	= [q^{d}](1-q)^{2g-2}.
	\end{align*}

	Let $\vec{d}=(d_1,\dots,d_r)$ and $X_{\vec{d}}=C^{[d_1]}\times\cdots\times C^{[d_r]}$. Then the multiplicative property of Euler characteristic implies  
	\begin{align*}
		\sum_{|\vec{d}|=d}e(X_{\vec{d}}) 	= [q^{d}](1-q)^{r(2g-2)}.
	\end{align*}

	Let $\IQ_{d}$ be the symplectic isotropic Quot scheme with $N=2n$. The fixed loci under the $\mathbb{C}^*$ action described in Section \ref{sec:fixed_loci}. The localization formula give us explicit expression for the Euler characteristics:
	\begin{equation*}
		\sum_{d=0}^{\infty}e(\IQ_d)q^d=2^r\binom{n}{r}(1-q)^{r(2g-2)}.
	\end{equation*}
	
	Since the isotropic Quot scheme are not necessarily smooth, the virtual Euler characteristic $e^{\vir}(\IQ_{d})$  may not coincide with the topological Euler characteristic. Define the formal power series
	\[A_{N,r,g}(q)=\sum_{d=0}^{\infty}e^{\vir}(\IQ_{d})q^d.\] 

The virtual localization formula gives
	\begin{align*}
		e^{\vir}(\IQ_d)=\sum_{d_1+d_2=d}^{}\sum_{w_1,w_2}\int_{\fix_{\vec{d}, \underline{k}}}^{} c(\fix_{\vec{d}, \underline{k}})\frac{c_{\mathbb{C}^*}(\nbun)}{e_{\mathbb{C}^*}(\nbun)}.
	\end{align*}
	We know how to evaluate the above integral (see Section \ref{sec:Euler_class_sympl}), but the details are computationally challenging. We do not a have a closed form expression or a conjecture for $A_{N,r,g}(q)$.
	
Over $\mathbb{P}^1$, we find a finite number of values using computers. We used Sagemath \cite{sagemath} for these calculations:
	\begin{align*}
		A_{4,2,0}(q)=&4 + 16q + 32q^2 + 112q^3 + (-396)q^4 + 6800q^5 + (-85856)q^6 + 1122544q^7+ \\& (-14660608)q^8  + 192011264q^9 + (-2520726176)q^{10} + 33164547968q^{11}+\cdots \\
		A_{6,2,0}(q)= &12+48q+96q^2+228q^3-3246q^4+\cdots \\
		A_{8,2,0}(q)=& 24+96q+192q^2+464q^3+\cdots
	\end{align*}
	We observe that $e^{\vir}(\IQ_d)$ differs from the topological Euler characteristic when $d\ge2$, which indicates that $\IQ_{d}$ is not smooth. When $d=0,1$, the space $\IQ_{d}$ is always smooth.

	\bibliographystyle{alphnum}
	\bibliography{Iso_quot}

\end{document}